\RequirePackage{amsthm}
\documentclass[sn-mathphys]{sn-jnl}
\usepackage{graphicx}
\usepackage{multirow}
\usepackage{amsmath,amssymb,amsfonts}
\usepackage{amsthm}
\usepackage{mathrsfs}
\usepackage[title]{appendix}
\usepackage{xcolor}
\usepackage{textcomp}
\usepackage{manyfoot}
\usepackage{booktabs}
\usepackage{algorithm}
\usepackage{algorithmicx}
\usepackage{algpseudocode}
\usepackage{listings}
\usepackage{hyperref}
\usepackage[outdir=./]{epstopdf}
\usepackage[colorinlistoftodos]{todonotes}
\newread\tmp

\theoremstyle{thmstyleone}
\newtheorem{theorem}{Theorem}[section]

\theoremstyle{thmstyletwo}

\newtheorem{remark}[theorem]{Remark}

\theoremstyle{thmstylethree}
\newtheorem{definition}[theorem]{Definition}

\begin{document}

\title[Article Title]{Topologically Interlocking Blocks inside the Tetroctahedrille}

\author*[1]{\fnm{Reymond} \sur{Akpanya} \orcid{https://orcid.org/0009-0009-0195-992X}}\email{akpanya@art.rwth-aachen.de}

\author[1]{\fnm{Tom} \sur{Goertzen} \orcid{https://orcid.org/0009-0003-3399-3416}}\email{goertzen@art.rwth-aachen.de}

\author[1]{\fnm{Alice C.} \sur{Niemeyer} \orcid{https://orcid.org/0000-0002-2163-3240}}\email{niemeyer@art.rwth-aachen.de}

\affil[1]{\orgdiv{Chair of Algebra and Representation Theory}, \orgname{RWTH Aachen University}, \orgaddress{\street{Pontdriesch 10-12}, \city{Aachen}, \postcode{52062}, \country{Germany}}}

\abstract{A topological interlocking assembly consists of rigid blocks together with a fixed frame, such that any subset of blocks is kinematically constrained and therefore cannot be removed from the assembly. 
In this paper we pursue a modular approach to construct (non-convex) interlocking blocks by combining finitely many tetrahedra and octahedra. This gives rise to polyhedra whose vertices can be described by the tetrahedral-octahedral~honeycomb, also known as tetroctahedrille. 
We show that the resulting interlocking blocks are very versatile and allow many possibilities to form topological interlocking assemblies consisting of copies of a single block.  We formulate a generalised construction of some of the introduced blocks to construct families of topological interlocking blocks. Moreover, we demonstrate a geometric application by using the tetroctahedrille to approximate given geometric objects. Finally, we show that given topological interlocking assemblies can be deformed continuously in order to obtain new topological interlocking assemblies.
}


\maketitle
\section{Introduction}\label{section:Introduction}
Since ancient times humans have strived to build constructions which are well adapted to their needs and withstand the demands placed upon them by time and the environment. As early as in the neolithic age, large stones were stacked to form dry stone walls in various countries around the world. Arches have long been guided by the stereotomic principle to use the weight of a stone and the neighbouring 
stones to hold it in place. In 1699, Joseph Abeille presented his design of an Abeille tile to the French Acad\'{e}mie des Sciences \citep{gallon_machines_1735}. This tile could be assembled forming a flat structure in such a way that each tile was prevented from moving by its neighbouring tiles. This idea forms the basis of the concept of topological interlocking. Other blocks have been found by Truchet and Frézier, and possess similar assembly rules as the Abeille tile \citep{frezier_theorie_1738}.

In this paper, a new method to create topologically interlocking blocks is presented. Starting from a three-dimensional lattice, whose points are the integer combinations of three basis vectors, we construct blocks with interlocking properties whose vertices are points of the given lattice. 
Given an assembly of blocks that lie in a given three-dimensional lattice, the position of the vertices of any block in the assembly can easily be computed using only translations and the symmetries of the given lattice. 
Thus, describing an interlocking block and the corresponding assemblies that lie in a lattice can be easily achieved.
The lattice of interest in this paper is the face-centred cubic lattice also known as the tetrahedral-octahedral~honeycomb or tetroctahedrille. The lattice, or more precisely the honeycomb corresponding to the lattice, can be created by alternating octahedra and tetrahedra to create a space filling. Here, we exploit this lattice to construct blocks that allow topological interlocking assemblies.
Known examples of topological interlocking assemblies, namely the tetrahedra interlocking and the octahedra interlocking, can be placed in this lattice by rigid motions, so that the vertices of the assembled blocks can be described by the lattice points.

\paragraph{Literature Review}

Many different blocks that allow topological interlocking assemblies are known. For instance, \cite{glickman_g-block_1984} describes an interlocking of regular tetrahedra  as a pavement system and \cite{dyskin_new_2001} introduce it as an example of a new material design concept. 
Moreover, \cite{dyskin_topological_2003} describe topological interlocking assemblies for all other Platonic solids and \cite{dyskin_fracture_2003} introduce a family of non-convex blocks with non-planar surface, called \emph{osteomorphic  blocks}. These blocks are very versatile and facilitate different design strategies.

Many researchers have been working on methods to design new interlocking blocks. For instance, \cite{kanel-belov_interlocking_2010} present a method to construct convex polyhedra that serve as interlocking blocks. Generalisations of this construction are exploited by \cite{wang_design_2019} and \cite{weizmann_topological_2016} to create candidates for curved interlocking assemblies consisting of convex blocks.

The interplay between plane tessellations and space-filling structures that lead to various topological interlocking assemblies with convex bodies is described by \cite{viana_topological_2018}.
Further, \cite{viana_towards_2021} gives a method for a systematic approach of creating topological interlocking assemblies with convex blocks and points out open problems in the context of topological interlocking assemblies.

\cite{subramanian_delaunay_2019,akleman_generalized_2020} introduce a general method for creating TI assemblies based on varying Voronoi Domains. \cite{topological22} exploit an Escher-like approach to obtain TIA by deforming (non-convex) fundamental domains of wallpaper symmetries. Moreover, \cite{spherical} investigate a group-theoretic approach to construct spherical TIA with identical parts.

Topological interlocking assemblies yield a concept with many applications in architecture and related fields. 
For instance, manipulating topological interlocking assemblies that form flat vaults yields various possibilities to create topological interlocking assemblies realising curved structures \citep{fallacara_topological_2019}.  The process of finding, analysing and manufacturing a structure based on the topological interlocking assembly of cubes is described by \cite{lecci_design_2021}. 
\cite{tessmann_topological_2012,tessmann_extremely_2013} present several topological interlocking assemblies, both flat and curved, obtained by various methods such as deformations and modifications of surfaces. Further, \cite{piekarski_floor_2020} reviews and gives new methods to generate floor slabs consisting of topological assemblies. For instance, Piekarski discusses the possibility of assembling inwards starting from the frame by introducing certain keystones which are obtained as divisions of interlocking blocks. Application of osteomorphic type blocks for the design and construction of a sustainable masonry system allowing a versatile approach are discussed by \cite{harsono_integration_2023}. More applications to masonry structures are described by \cite{moreno_gata_designing_2019}.

For a recent overview on the research progress on TI assemblies, see \citep{dyskin_topological_2019,estrin_design_2021}.

\paragraph{Structure of the Paper}

We first present the tetroctahedrille, a three-dimensional lattice that gives rise to a quasi regular space filling consisting of octahedra and tetrahedra, and further examples of polyhedra that can be found within this lattice, see Section~\ref{section:Tetroctahedrille}. In Section~\ref{sections:approxiamtions}, we show that any geometric object can be approximated inside the tetroctahedrille. 
Next, we give an abbreviated definition of a topological interlocking assembly in Section~\ref{section:TopologicalInterlocking}. We then discuss interlocking blocks that can be constructed by exploiting the lattice of interest in Section~\ref{section:NewBlocks}. Further, we describe corresponding assemblies for every presented interlocking block. Last, we elaborate on possible deformations that can be applied to the presented blocks to create new blocks preserving the interlocking property, see Section~\ref{section:Modyfying}.

\section{Tetroctahedrille}\label{section:Tetroctahedrille}
Exploiting a modular approach to construct new polyhedra by combining finitely many polyhedra, such as the Platonic solids, can yield geometric shapes that allow various applications \citep{stuttgen_modular_2023}. In this paper, we present polyhedra that are computed by using a finite number of octahedra and tetrahedra and elaborate on their geometric properties.
This section covers the necessary notions and details of the structures that arise from combining tetrahedra and octahedra.

\begin{figure}[H]
\begin{minipage}{.3\textwidth}
    \centering
    \includegraphics[height=4.cm]{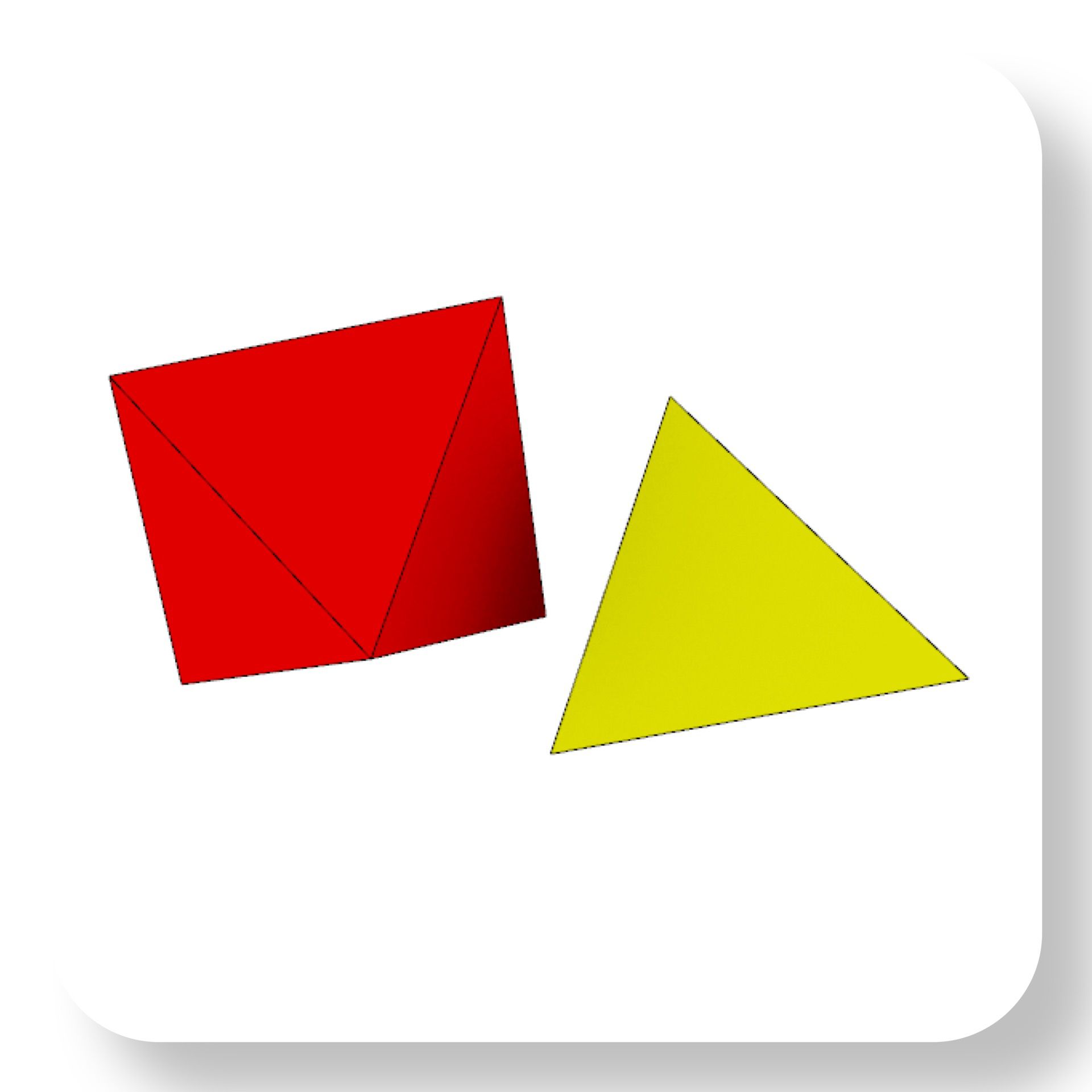}
    \label{octtet1}
\end{minipage}
\begin{minipage}{0.5cm}
    
\end{minipage}
\begin{minipage}{.3\textwidth}
    \centering
    \includegraphics[height=4.cm]{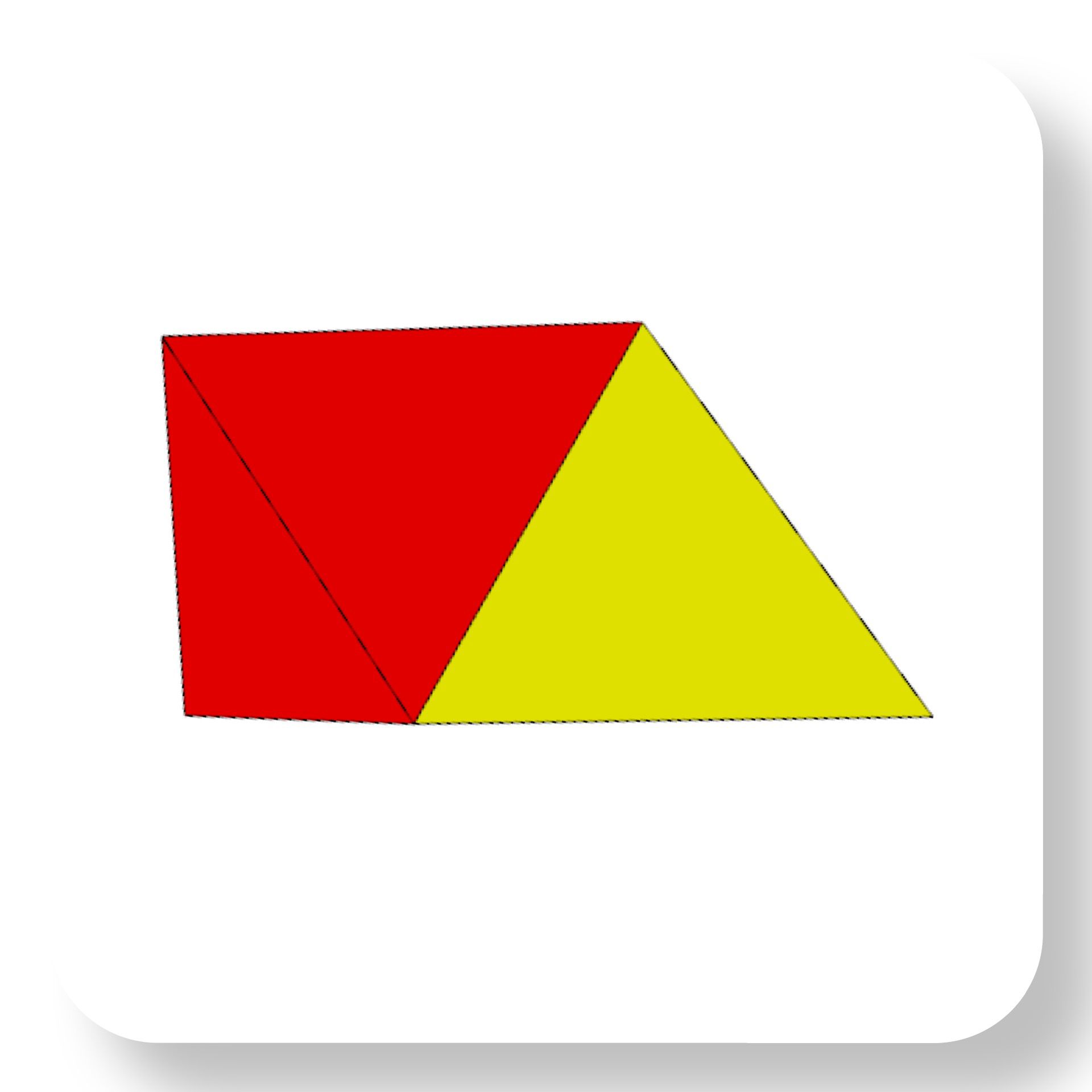}
    \label{octtet2}
\end{minipage}
\begin{minipage}{0.5cm}
    
\end{minipage}
\begin{minipage}{.3\textwidth}
    \centering
    \includegraphics[height=4.cm]{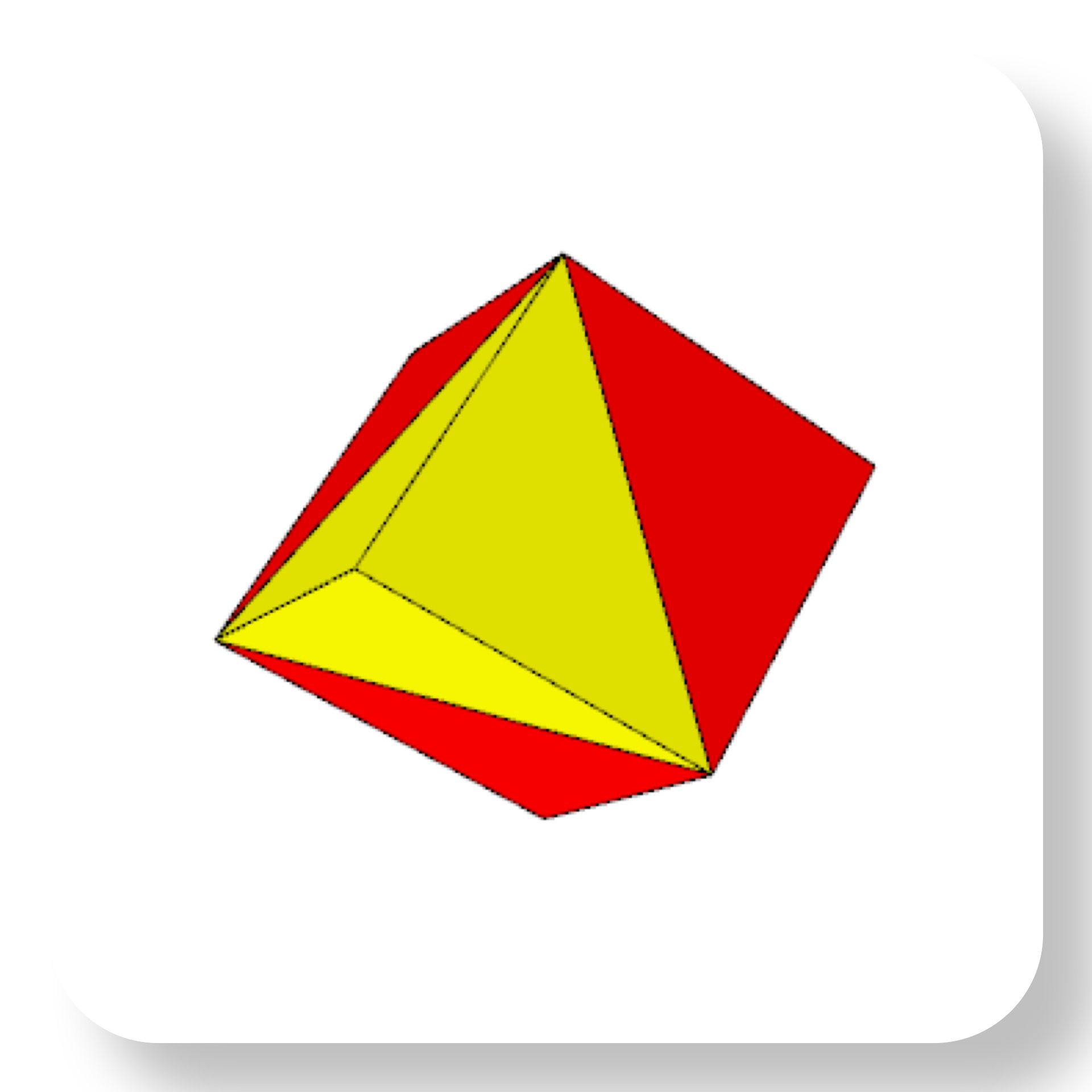}
    \label{octtet3}
\end{minipage}
\caption{Attaching a tetrahedron to an octahedron so that two faces are identified}
\label{octet}
\end{figure}

The fact, that copies of tetrahedra and octahedra can be exploited to construct structures such as lattices and space fillings, can be seen by the following observation:
\begin{remark}
A rigid regular tetrahedron and a rigid octahedron with all edge lengths $\sqrt{2}$ can be placed in three-dimensional space so that a face $F_T$ of the tetrahedron coincides with a face $F_O$ of the octahedron. Moreover, each edge of $F_T$ coincides with exactly one edge of $F_O$. The resulting dihedral angles at these edges add up to $\pi$, see Figure~\ref{octet}. 
    
\end{remark}

We place polyhedra, which are copies of rigid regular tetrahedra and octahedra with all edges of same length in $\mathbb{R}^3$, so that 
\begin{enumerate}
    \item[(R1)] if two polyhedra have a non-empty intersection, then this intersection is either a vertex, an edge or a face of both polyhedra,
    \item[(R2)] the intersection of two polyhedra is a face if and only if one polyhedron is a tetrahedron and the other polyhedron is an octahedron. 
\end{enumerate}
If, in addition, each face is the intersection of a tetrahedron and an octahedron, this arrangement
gives rise to a quasi regular space filling, i.e. following the formulated rules allows to fill the Euclidean space with regular tetrahedra and octahedra without gaps.
In the literature, this space filling is often referred to as \textit{tetroctahedrille} or \textit{octahedral~tetrahedral~honeycomb} and has various applications in scientific fields such as crystallography and architecture. 

An alternative way to define the tetroctahedrille is to describe the vertices of the individual tetrahedra and octahedra in the space filling as points of a three-dimensional lattice, namely the \emph{face-centred~cubic~lattice}.
The tetrahedra and octahedra with edge lengths $\sqrt{2}$ in this honeycomb can be described as subsets of the following set of integer points:
\[
\langle e_1+e_2,e_1+e_3,e_2+e_3\rangle_\mathbb{Z}, 
\]
where $\{e_1,e_2,e_3\}$ denotes the standard basis of the Euclidean 3-space. 
In this paper we will present various examples of polyhedra that can be cut out from the tetroctahedrille or, alternatively, that can be constructed by combining finitely many tetrahedra and octahedra following the above rules.
For instance, by satisfying the above rules, we can assemble a tetrahedron with edge lengths $2\sqrt{2}$ from 1 octahedron and 3 tetrahedra with edge lengths $\sqrt{2}$ and an octahedron with edge length $2\sqrt{2}$ from 6 octahedra and 4 tetrahedra with edge lengths $\sqrt{2}$, see Figure~\ref{tetra_oct_explode}. 
For a more general introduction to this honeycomb, we refer to \citep{conway}.

\begin{figure}[H]
\begin{minipage}{.3\textwidth}
    \centering
    \includegraphics[height=4cm]{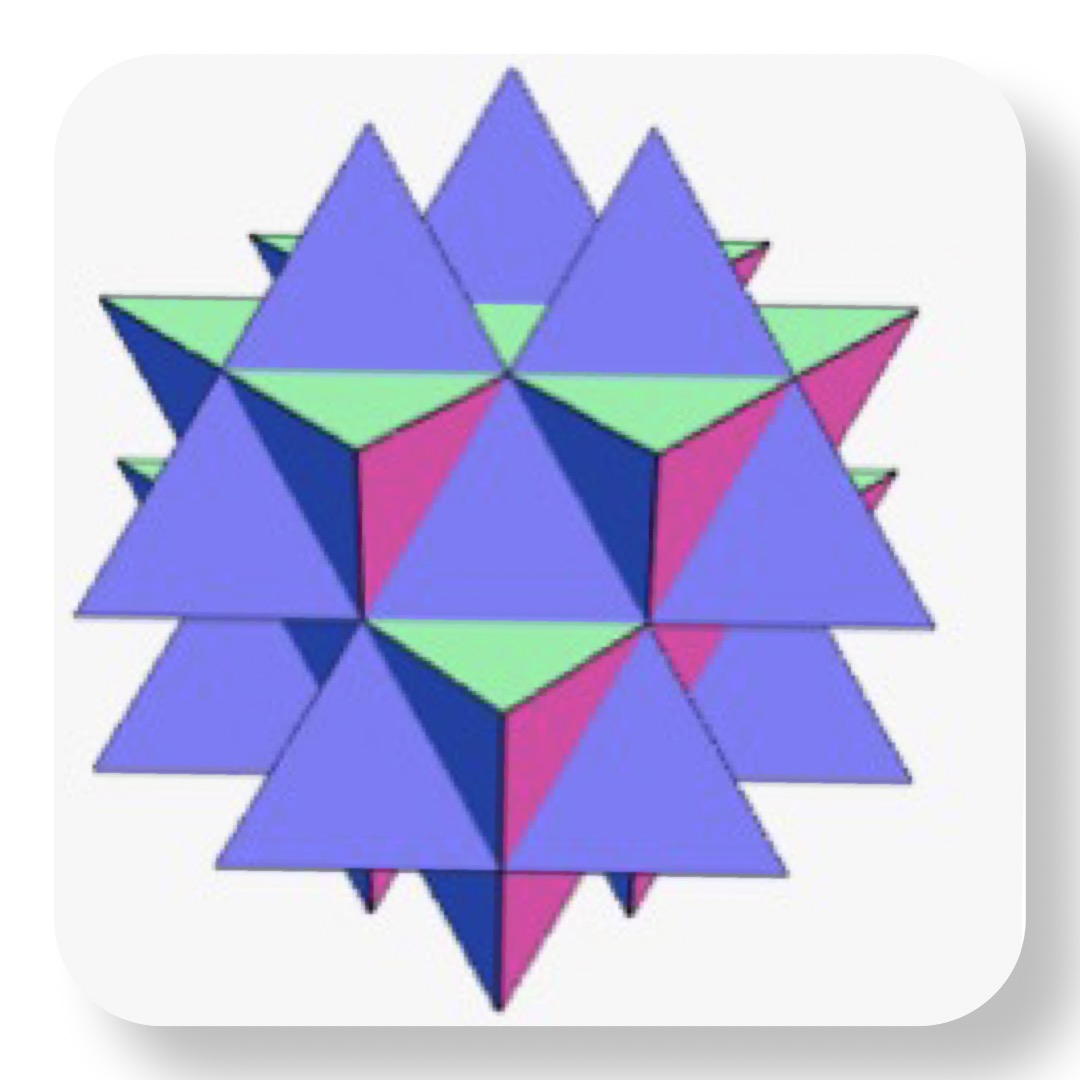}
\end{minipage}
\begin{minipage}{0.5cm}
    
\end{minipage}
\begin{minipage}{.3\textwidth}
    \centering
    \includegraphics[height=4cm]{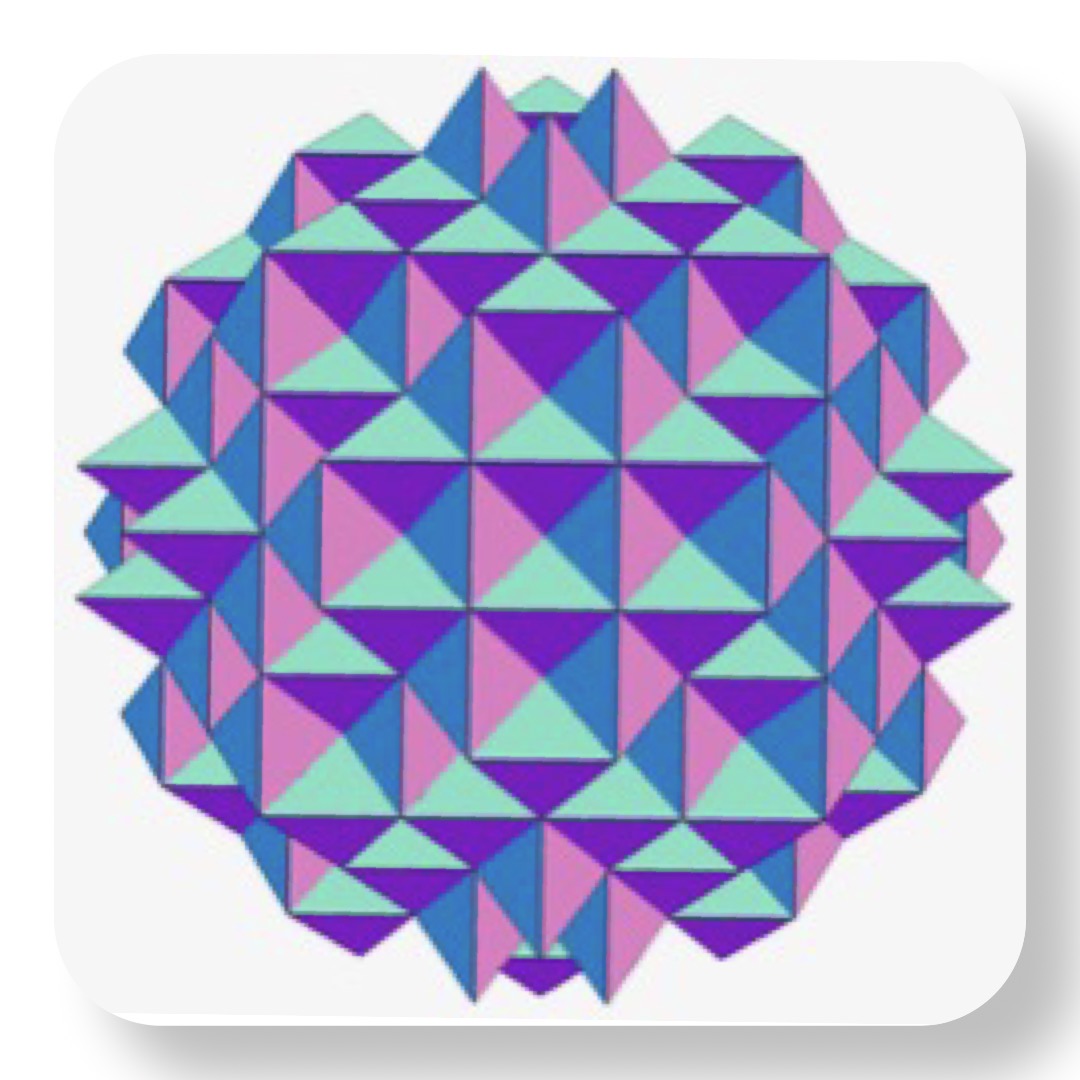}
\end{minipage}
\begin{minipage}{0.5cm}
    
\end{minipage}
\begin{minipage}{.3\textwidth}
    \centering
    \includegraphics[height=4cm]{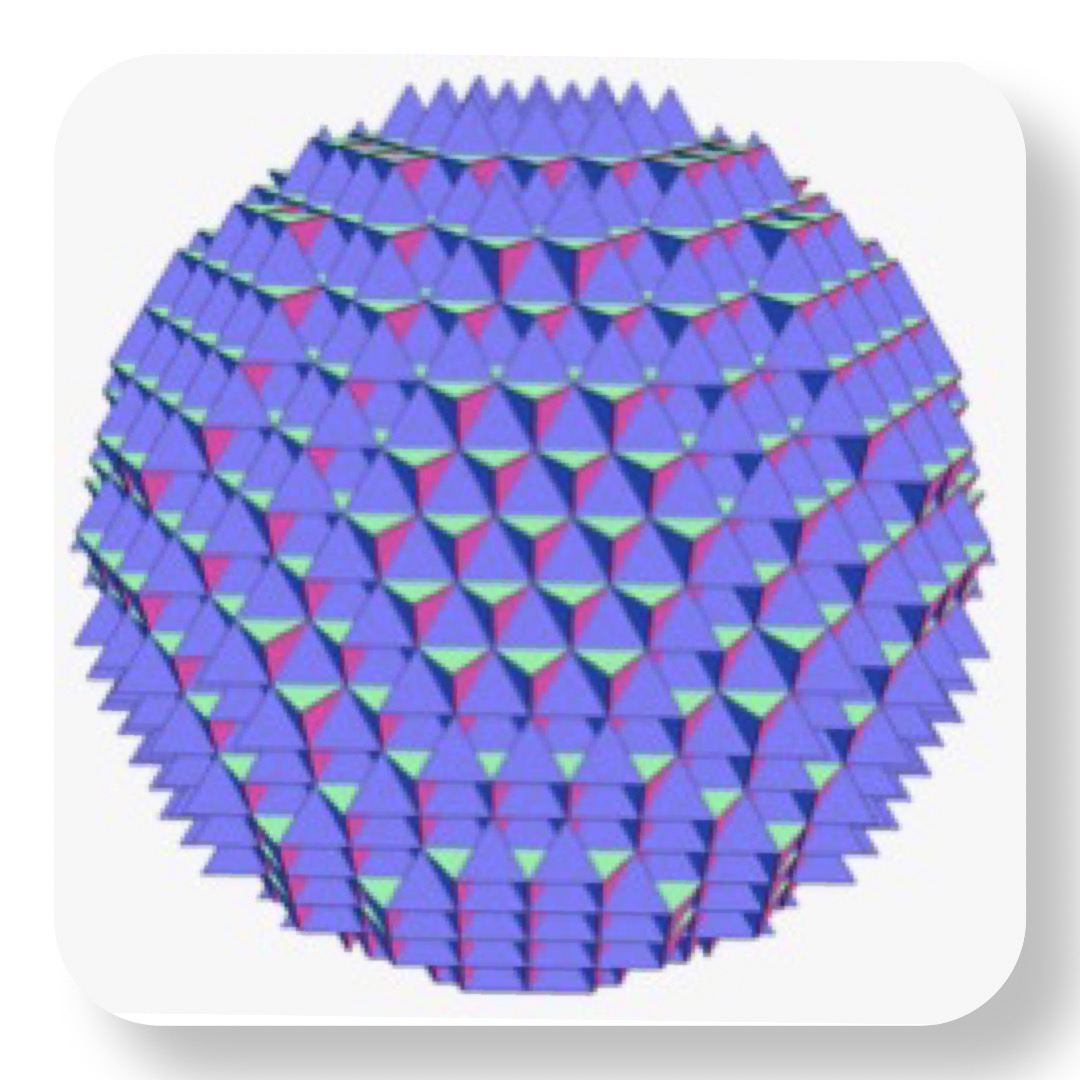}
\end{minipage}
\caption{Approximations of the unit sphere with scaled tetrahedra and octahedra }
\label{approximation}
\end{figure}

Since regular tetrahedra and octahedra are convex, they can be described by their sets of incident vertices.
We therefore describe the tetrahedra in the tetroctahedrille by the sets of their four vertices and the octahedra by the sets of their six vertices.
Thus, every tetrahedron and octahedron in this honeycomb can be described as follows: 
\begin{remark}
\label{v1v2v3}
 Let $0$ be defined by $0:=(0,0,0)^t$ and $v_1,v_2,v_3$ be three-dimensional coordinates in $\mathbb{R}^3$ given by: 
    \begin{align*}
        v_1:=e_2+e_3,
        v_2:=e_1+e_3,
        v_3:=e_1+e_2.
    \end{align*}
Hence, it is possible to describe every octahedron 
in the tetroctahedrille by a translation of the following octahedron:
\[
\{v_1,v_2,v_3,v_1+v_2,v_1+v_3,v_2+v_3\}. 
\]
Further, the tetrahedra in this honeycomb can be described by applying translations to one of the following tetrahedra:
\begin{align*}
   & T_1=\{0,v_1,v_2,v_3\},\\
   & T_2=\{v_2,v_3,v_2+v_3,v_2+v3-v_1\},\\
   & T_3=\{v_1,v_2,v_1+v_2,v_1+v_2-v_3\},\\
   & T_4=\{v_2,v_1+v_2,v_2+v_3,2v_2\}.\\
\end{align*}
We denote a polyhedron in the tetroctahedrille that is obtained by applying a translation by a vector $v\in \mathbb{R}^3$ to a polyhedron $X$ by $v+X$. For instance, the tetrahedron $v_1+T_1$, that arises from applying the translating by the vector $v_1$ to the polyhedron $T_1$, can be represented by the following set of coordinates:
\[
v_1+T_1=\{v_1,2v_1,v_1+v_2,v_1+v_3\}.
\]
\end{remark}
One application of the tetrahedral-octahedral~honeycomb is the approximation of surfaces of given three-dimensional bodies. An intuitive approach to achieve this approximation is to construct a polyhedral complex $P$ such that 
\begin{enumerate}
    \item $P$ is fully contained in the tetroctahedrille,
    \item $P$ is fully contained in the given three-dimensional body,
    \item $P$ is maximal in the sense that we cannot add any tetrahedra or octahedra to obtain a polygonal complex that satisfies 1. and 2. 
\end{enumerate}
 By scaling the tetrahedra and octahedra in the tetroctahedrille or equivalently scaling the given three-dimensional body, we can obtain more exact approximations of the given body. Note, that the triangulated surface that results from the above approximation consists of equilateral triangles with the same edge lengths.  
Figure~\ref{approximation} demonstrates this approximation on the unit sphere, where the common edge length of the tetrahedra and octahedra varies.
More details on the approximation of surfaces by exploiting the tetroctahedrille along with more complex examples can be found in Section~\ref{sections:approxiamtions}. 

\begin{figure}[H]
\begin{minipage}{.3\textwidth}
    \centering
    \includegraphics[height=4cm]{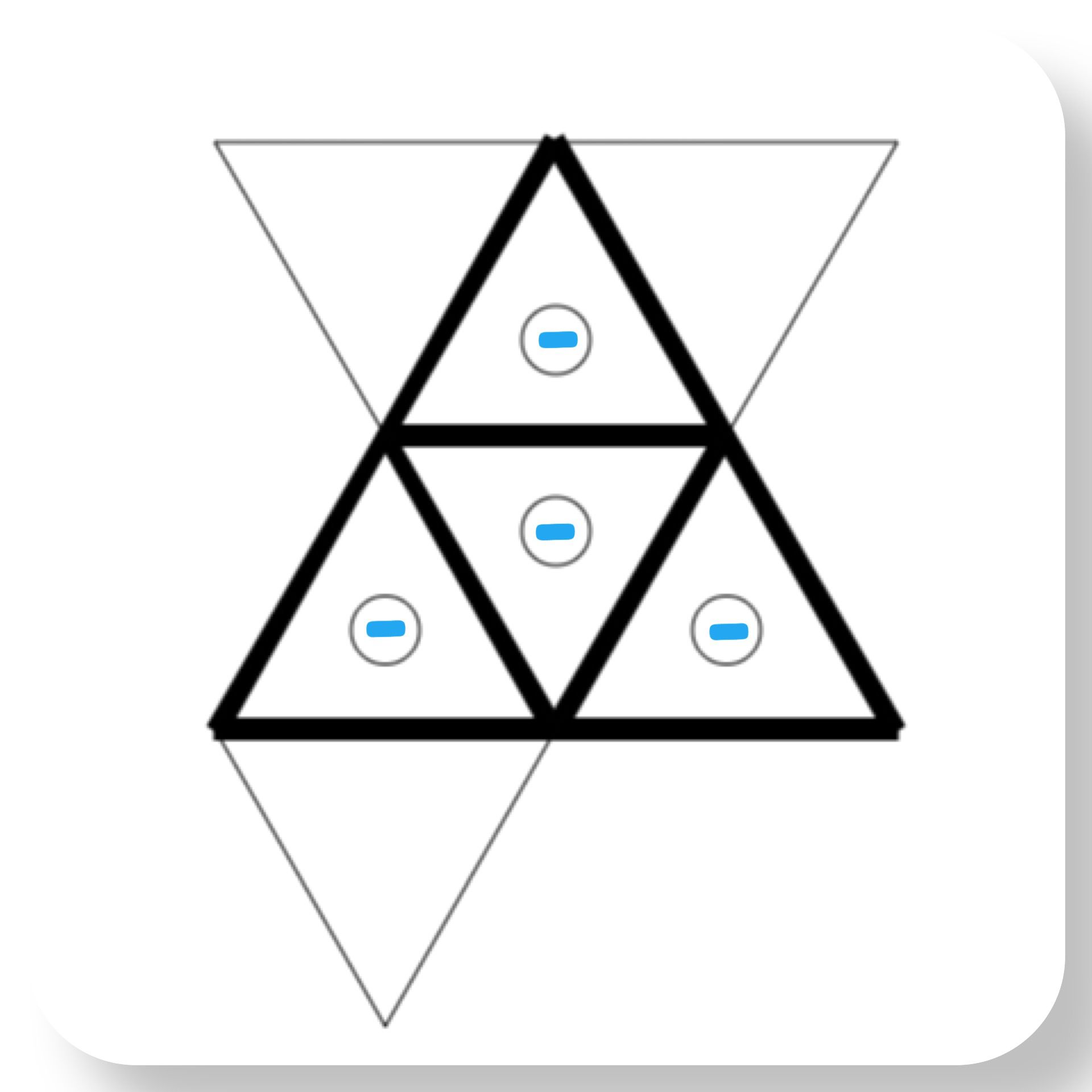}
\end{minipage}
\begin{minipage}{0.5cm}
    
\end{minipage}
\begin{minipage}{.3\textwidth}
    \centering
    \includegraphics[height=4cm]{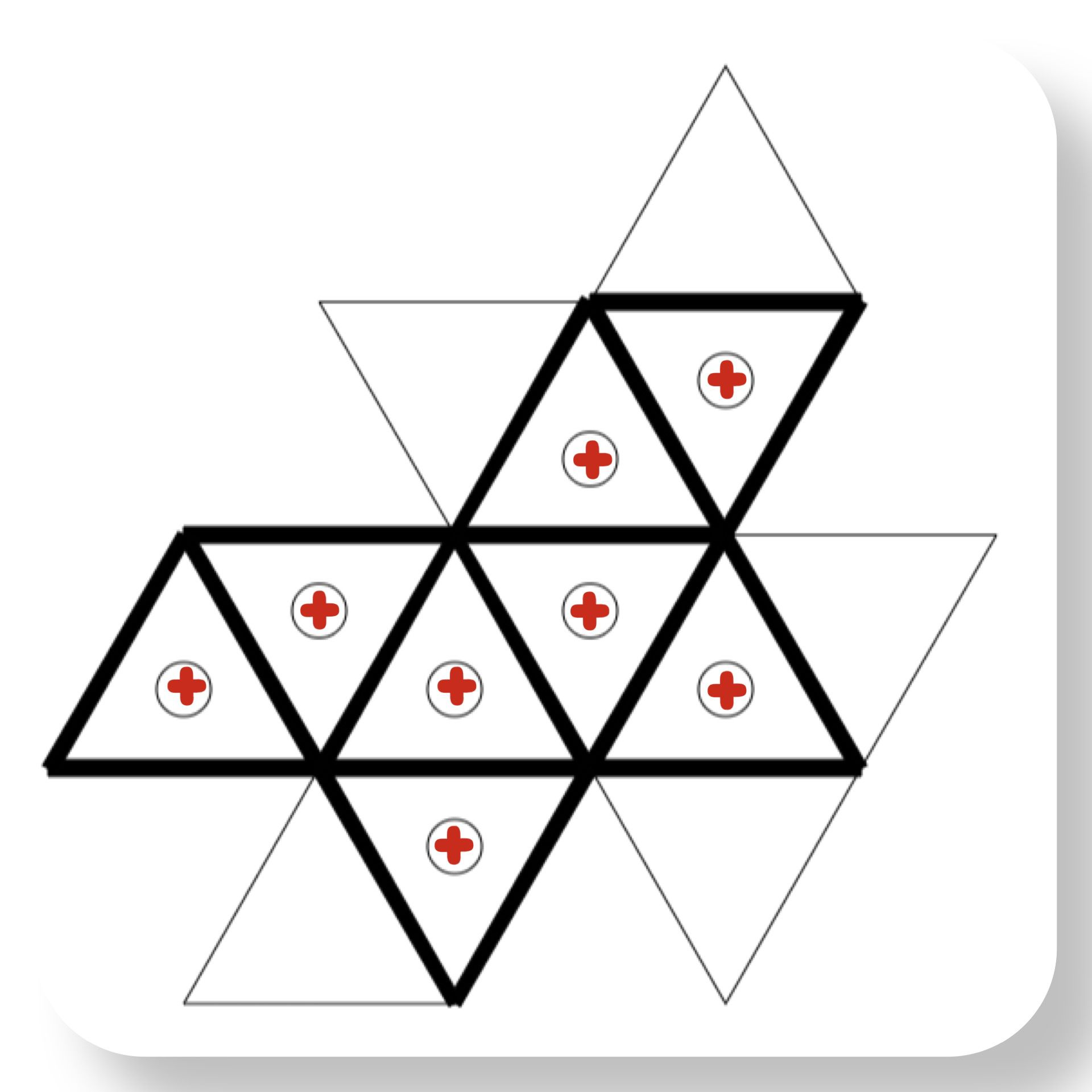}
\end{minipage}
\begin{minipage}{0.5cm}
    
\end{minipage}
\begin{minipage}{.3\textwidth}
    \centering
    \includegraphics[height=4cm]{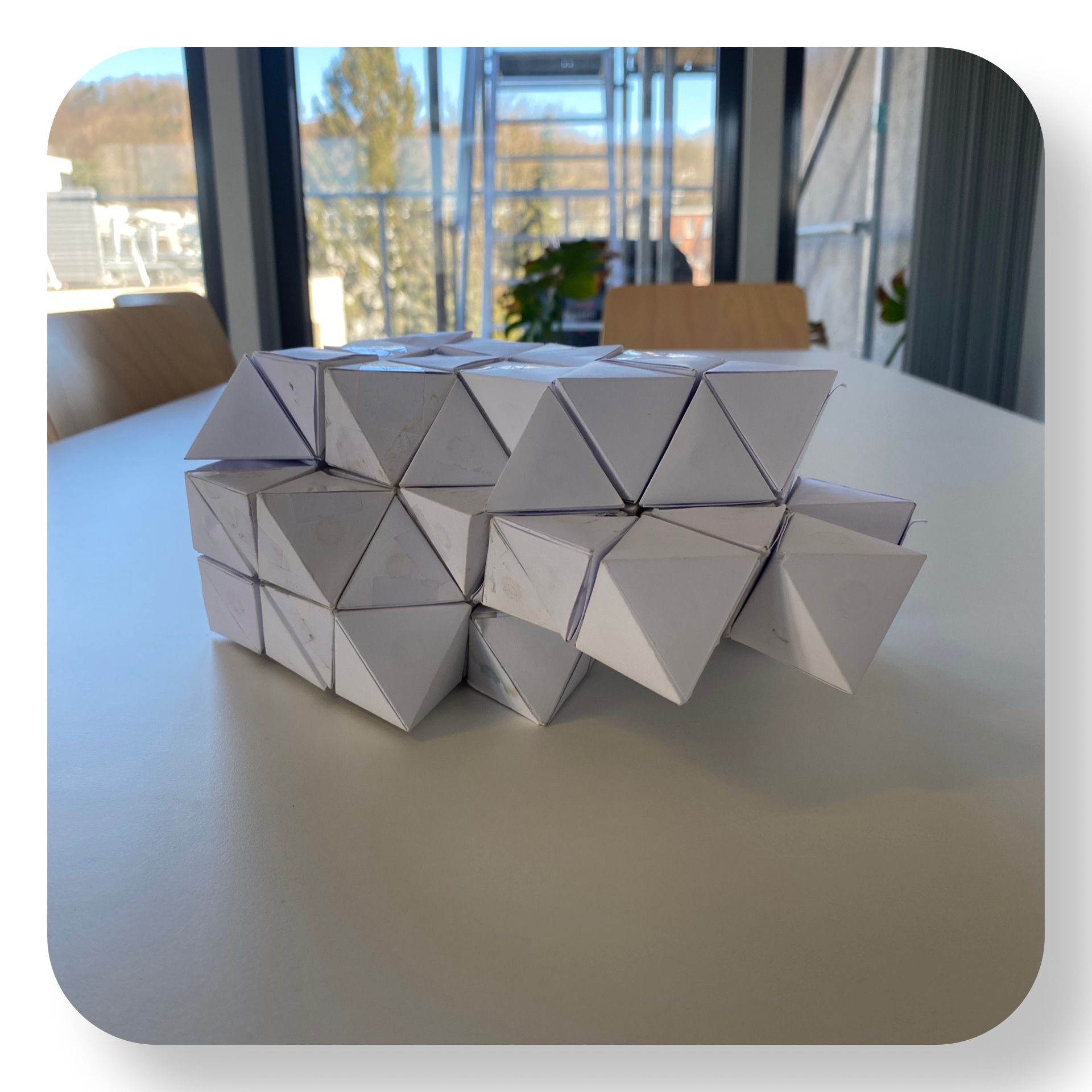}
\end{minipage}
\caption{(Left, middle) Folding plans for the tetrahedron and octahedron with $+,-$ signs according to the poling; (right) a block constructed by assembling tetrahedra and tetrahedra with inserted magnets}
\label{magnets}
\end{figure}

In order to explore the properties and possibilities of the given honeycomb further, we can build octahedra and tetrahedra in real life.
For instance, we use magnetic paper models incorporating the rules of the tetrahedral-octahedral~honeycomb to construct polyhedra that are contained in the honeycomb. We glue magnets to the different tetrahedra and octahedra in such a way that only tetrahedra and octahedra can touch.
This approach leads to folding plans with circles corresponding to the magnets and $+,-$ signs according to the polarity, see Figure~\ref{magnets}.

Such models can be  used to construct large structures contained inside the tetrahedral-octahedral~honeycomb without being concerned about structural failure, since magnetic forces can stabilise this construction. \cite{bridges23} similarly use magnets to understand the self-assembling behaviour of copies of the Versatile Block, another topological interlocking block.

\section{Approximating Objects using the Tetroctahedrille }\label{sections:approxiamtions}

As described above, the tetroctahedrille fills the three-dimensional Euclidean space. In this section, we exploit this fact to construct approximations of geometric objects that solely consist of tetrahedra and octahedra. This can then be used for constructing TI assemblies which approximate the given geometric object.

\begin{figure}[H]
\centering
\begin{minipage}{.45\textwidth}
    \centering
    \includegraphics[height=3.5cm]{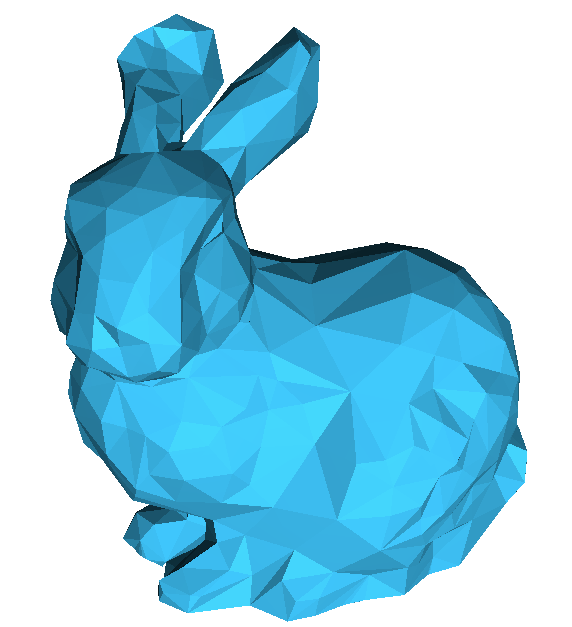}
\end{minipage}
\begin{minipage}{.45\textwidth}
    \centering
    \includegraphics[height=3.5cm]{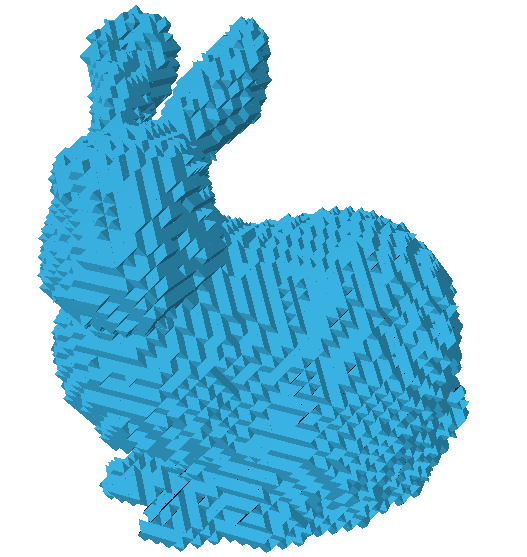}
\end{minipage}
\caption{An approximation of the Stanford Bunny \url{http://graphics.stanford.edu/data/3Dscanrep/} (left) and an approximation generated using the tetroctahedrille (right)}
\label{stanfordbunny}
\end{figure}

 A given geometric object such as the Stanford Bunny, shown in Figure \ref{stanfordbunny}, can be approximated using the tetroctahedrille. For this we assume that we receive a triangulated surface, for instance in form of an STL-file.
In order to obtain the desired approximation, we proceed as follows:
\begin{enumerate}
    \item Simplify the original model to guarantee cost-efficient computing, i.e. reduce the number of triangles;
    \item For each triangle compute finitely many convex combinations of its vertices;
    \item For each convex combination compute the corresponding cell in the tetroctahedrille;
    \item Unify all cells to receive an approximation of the original data within the tetroctahedrille.
\end{enumerate}

\begin{figure}[H]
\centering
\begin{minipage}{.45\textwidth}
    \centering
    \includegraphics[height=3.5cm]{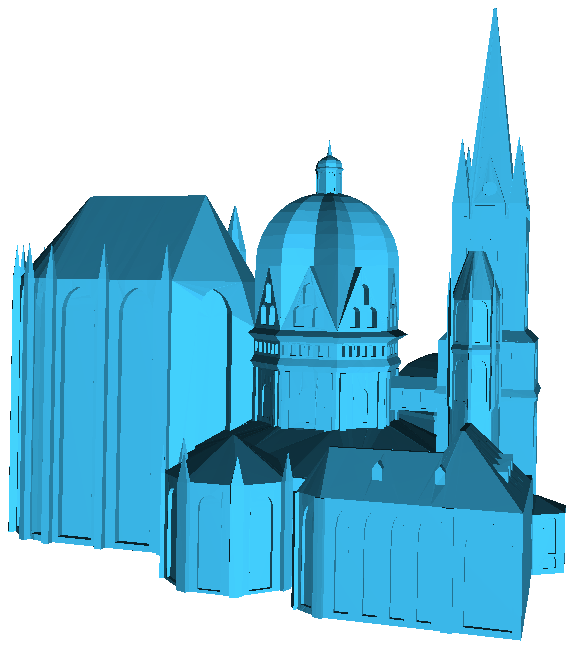}
\end{minipage}
\begin{minipage}{.45\textwidth}
    \centering
    \includegraphics[height=3.5cm]{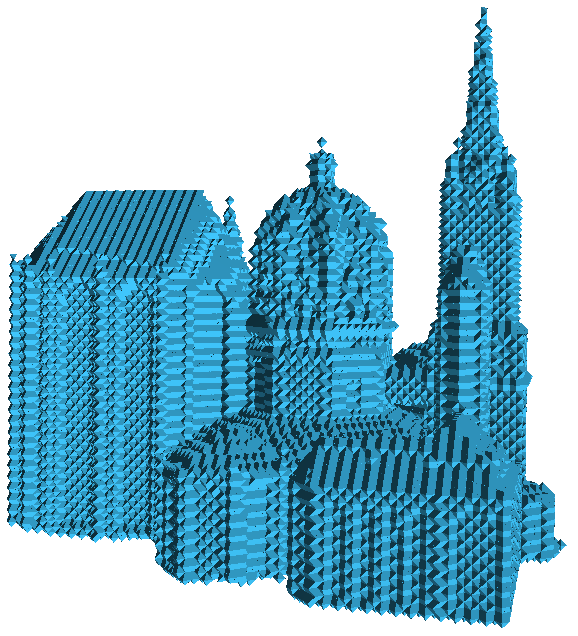}
\end{minipage}
\caption{An STL-model of the cathedral in Aachen given by Philipp Hoffmann available at \url{https://www.thingiverse.com/thing:4820139} (left) and an approximation generated using the tetroctahedrille (right)}
\label{AachenCathedral}
\end{figure}

We can generate arbitrarily exact approximations by scaling the model, while keeping the size of tetrahedra and octahedra constant. In Figure \ref{stanfordbunny} we present an approximation of the Stanford bunny using the tetroctahedrille.
As another example we generate an approximation of the cathedral in Aachen using the ideas described above. The resulting approximation can be seen in Figure \ref{AachenCathedral}.

This leads to an alternative approach for designing
assemblies with rigid parts. \cite{WangAssembliesRigidParts} provide a recent overview of such computational methods.

\section{Definition of Topological Interlocking Assemblies}\label{section:TopologicalInterlocking}

A topological interlocking assembly consists of independent blocks, some of them fixed, such that no group of blocks can be moved. The fixed blocks are called the \emph{frame of the assembly}. Note, that the blocks of the assembly do not have to be congruent, i.e. identical copies of one block.

\begin{figure}[H]
\begin{minipage}{.5\textwidth}
    \centering
    \includegraphics[height=5.cm]{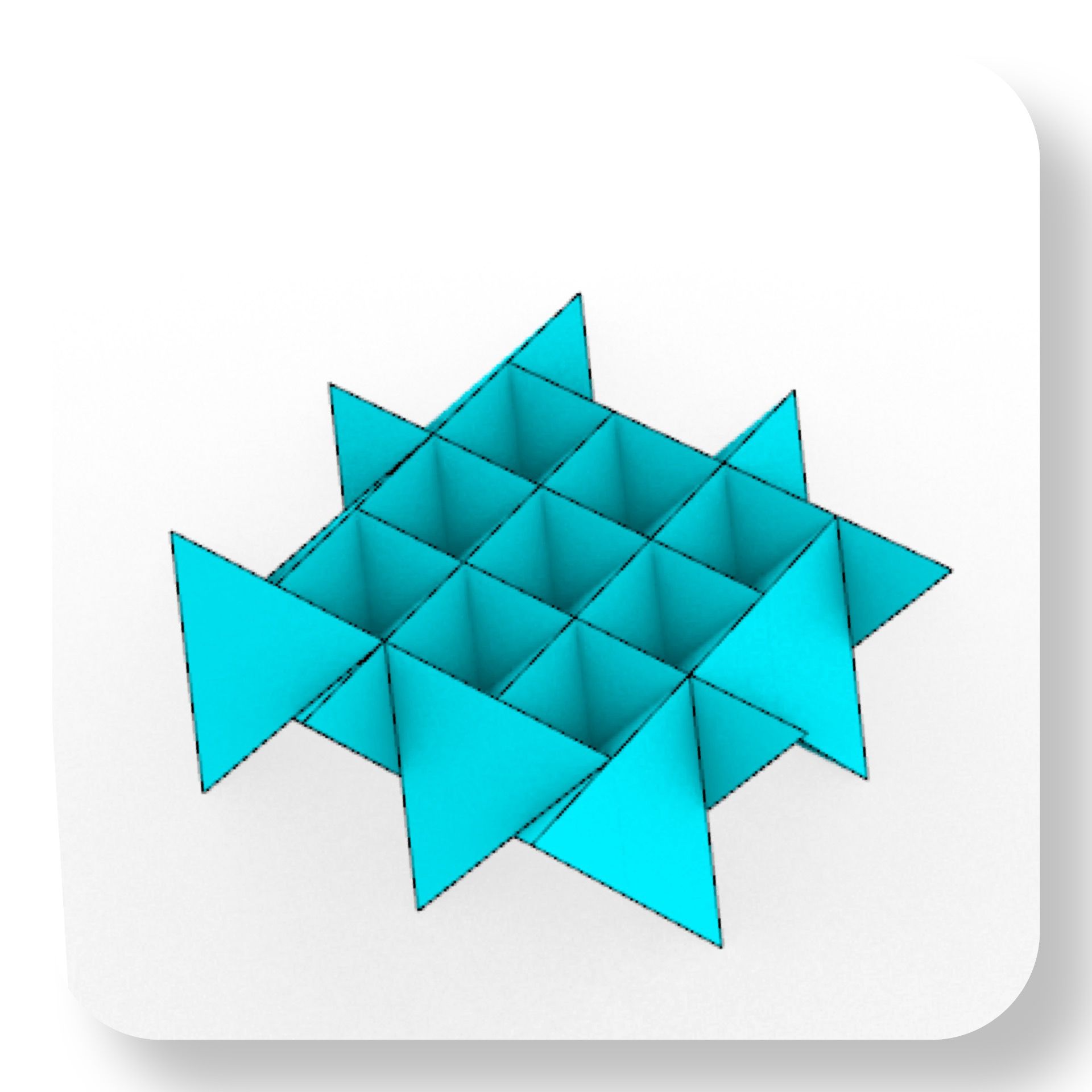}
    \label{TetrahedraInterlocking}
\end{minipage}
\begin{minipage}{.5\textwidth}
    \centering
    \includegraphics[height=5.cm]{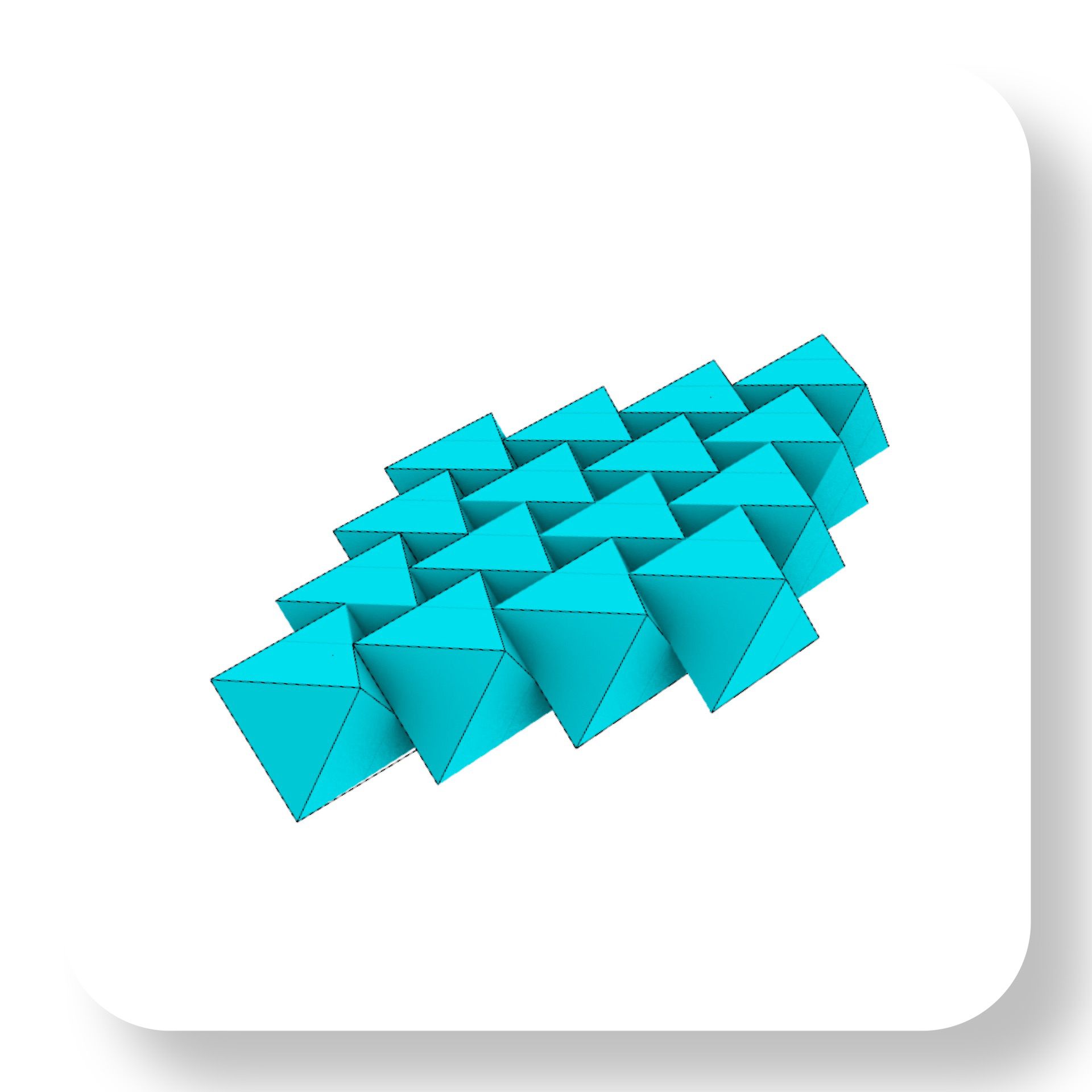}
    \label{OctahedraInterlocking}
\end{minipage}
\caption{Tetrahedra- and octahedra-interlocking assemblies, see \citep{dyskin_topological_2003}}
\label{tetra_oct_interlocking}
\end{figure}

For the purpose of this paper, a \emph{block} is a three-dimensional subset $X\subset \mathbb{R}^3$ with polyhedral boundary  denoted by $\partial X$. For simplicity, we only focus on topological interlocking assemblies that consist of identical blocks. Examples of such topological interlocking assemblies are the tetrahedra- and octahedra-interlocking assemblies, see Figure~\ref{tetra_oct_interlocking}.

An \emph{assembly of blocks} is a family $(X_i)_{i\in I}$ of blocks, where $I$ is a countable but not necessarily finite index set such that any two blocks only meet at their boundary, i.e. $$X_i \cap X_j = \partial X_i \cap \partial X_j.$$ For the assemblies in Figure~\ref{tetra_oct_interlocking}, the frame is chosen to consist of the outer blocks. However, assemblies with isolated frame blocks can also be constructed, see \citep{key} for instance. Here, the authors describe assemblies in the context of interlocking puzzles whose frames consist of single blocks, called keys. 

\cite{topological22} formulate a mathematical definition of TI assemblies and Definition \ref{def:interlocking} below is a version adjusted for the context of this paper.
Describing the motion of a block mathematically requires the definition of continuous motions.

\begin{definition}\label{def:interlocking}
A continuous (differentiable) motion is a continuous (differentiable) map $\gamma:[0,1]\to SE(3)$ such that $\gamma(0)=\mathbb{I}$, where $SE(3)$ is the group of rigid Euclidean motions composed of rotations and translation and $\mathbb{I}$ is the identity element.
\end{definition}
\noindent
With the notion of a continuous motion, we can present the definition of a topological interlocking.
\begin{definition}[\cite{topological22}]
Let $(X_i)_{i\in I}$ be an assembly of blocks for a countable index set $I$ and let $J\subset I$. The pair $((X_i)_{i\in I},J)$ is a \emph{topological interlocking assembly} (short TI assembly) with \emph{frame} $J$ if for all finite non-empty subsets $\emptyset \neq S\subset I\setminus J$ and for all non-trivial piece-wise differentiable motions $(\gamma_i)_{i\in S}$ there exists $t\in [0,1]$ such that 
$$ Y_i^t = \begin{cases} 
          X_i & i \in I\setminus S \\
          \gamma_i(t)(X_i) & i \in S 
       \end{cases}
    $$
    is not an assembly of blocks, i.e.\ there exists $i,j \in I$ such that $Y_i^t\cap Y_j^t\neq \partial Y_i^t\cap \partial Y_j^t$ meaning any small perturbations (movements) of blocks lead to intersections.
    
\end{definition}
Proving that an assembly of blocks is a TI assembly turns out to be a high-dimensional problem as motions of any group of blocks have to be considered simultaneously. 

\cite{wang_design_2019} show that proving that an assembly of blocks is a TI assembly can be done by solving a linear optimisation problem, whose solutions describe all infinitesimal movements of blocks that do not lead to intersections.

\begin{figure}[H]
\begin{minipage}{.5\textwidth}
    \centering
    \includegraphics[height=5cm]{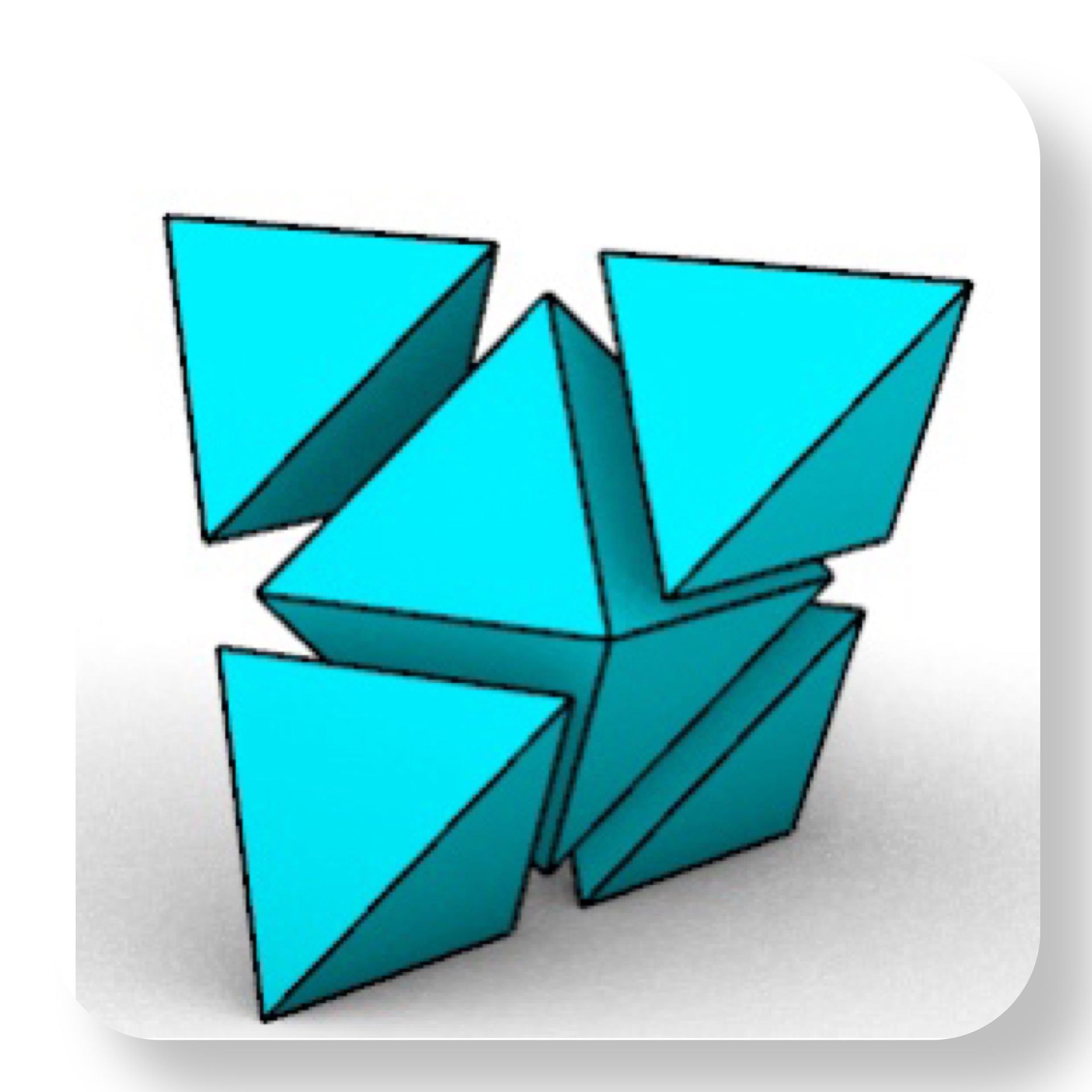}
    \label{TetrahedronExplode}
\end{minipage}
\begin{minipage}{.5\textwidth}
    \centering
    \includegraphics[height=5cm]{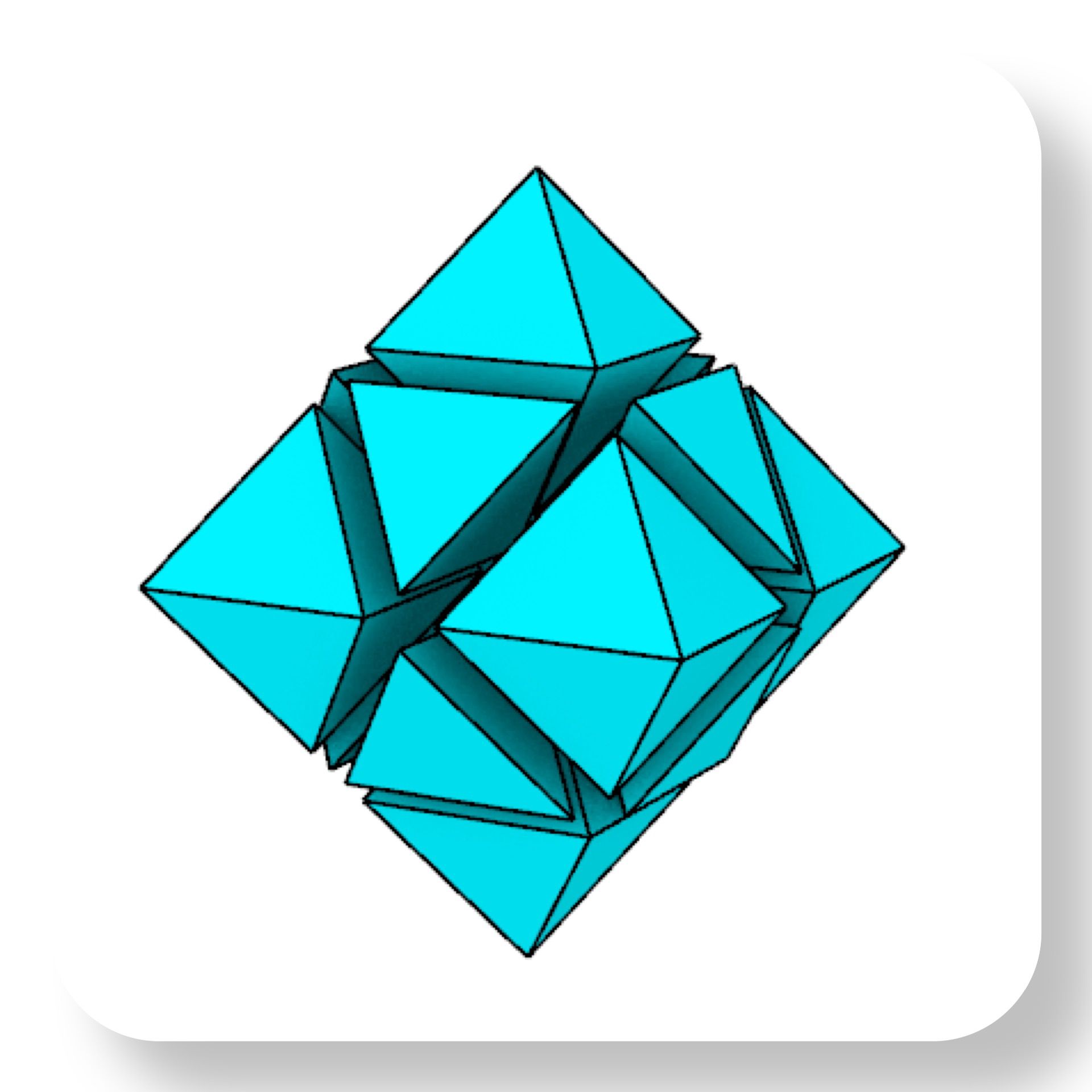}
    \label{OctahedronExplode}
\end{minipage}
\caption{Exploded view of tetrahedron and octahedron inside the tetrahedral-octahedral~honeycomb}
\label{tetra_oct_explode}
\end{figure}

\begin{figure}[H]
\begin{minipage}{.5\textwidth}
    \centering
    \includegraphics[height=5cm]{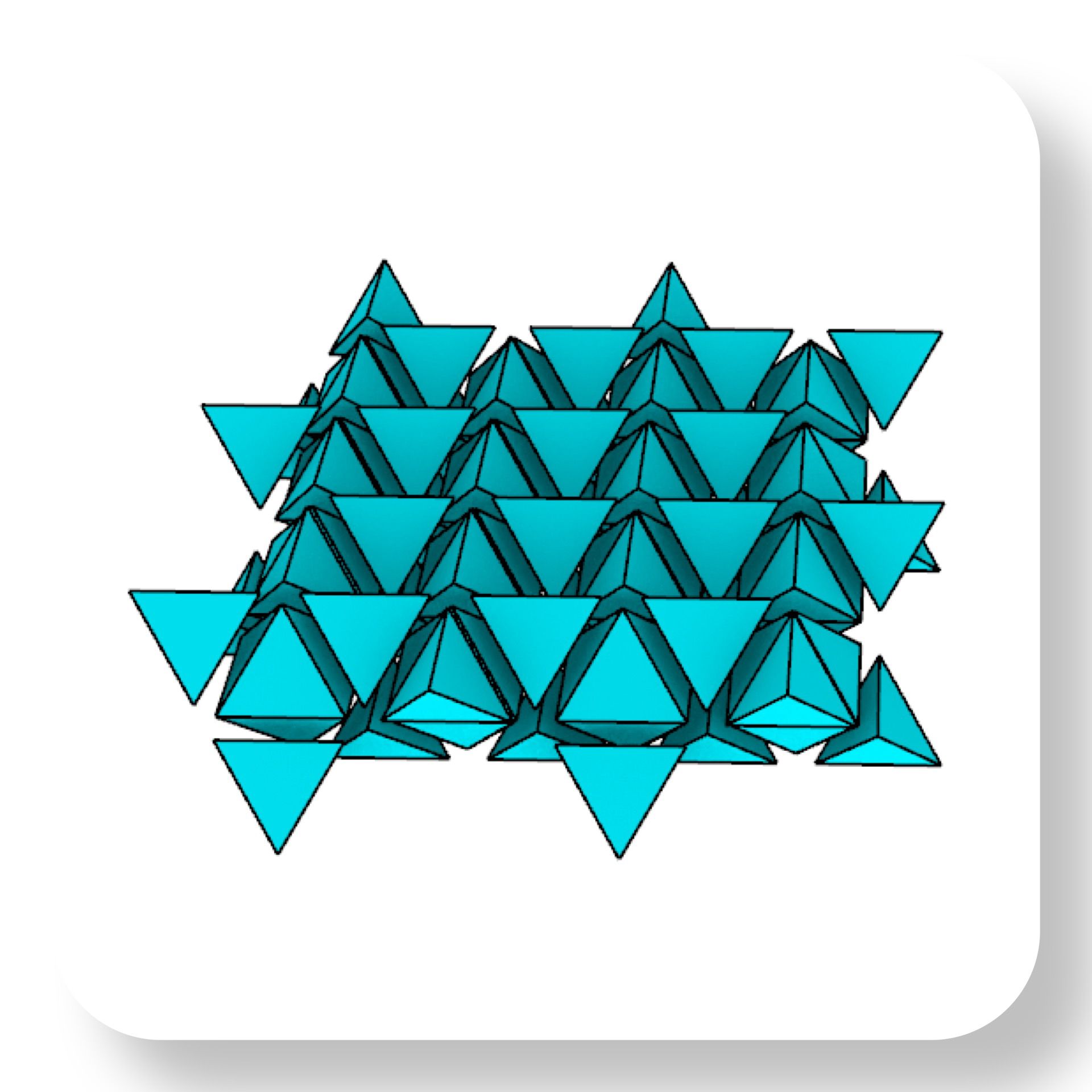}
    \label{TetrahedronExplodeInt}
\end{minipage}
\begin{minipage}{.5\textwidth}
    \centering
    \includegraphics[height=5cm]{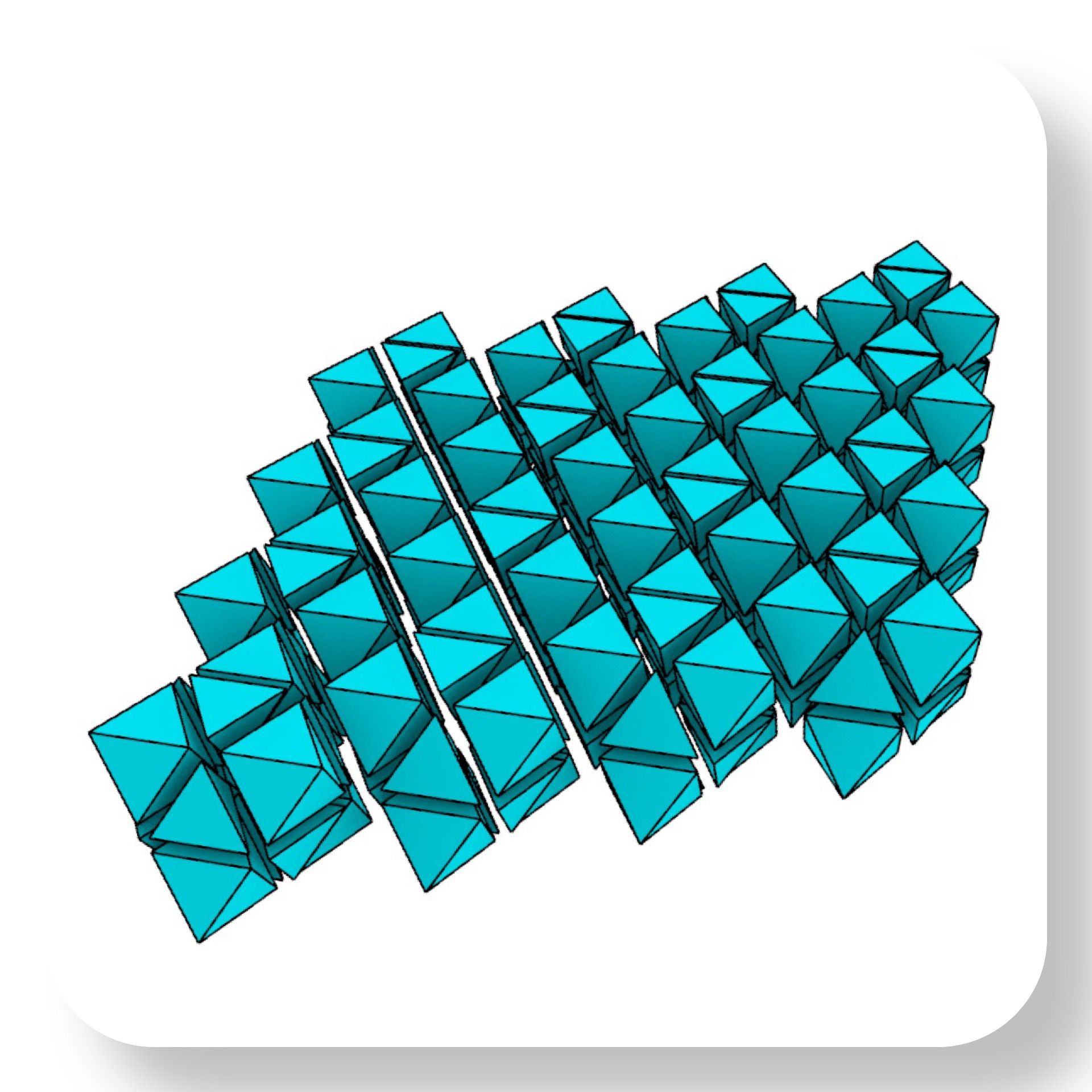}
    \label{OctahedronExplodeInt}
\end{minipage}
\caption{Exploded view of tetrahedra and octahedra TI assemblies inside the tetrahedral-octahedral~honeycomb}
\label{tetra_oct_interlocking_explode}
\end{figure}

Both the tetrahedra and octahedra TI assemblies, shown in Figure~\ref{tetra_oct_interlocking}, can be viewed as part of the tetrahedral-octahedral~honeycomb as ``larger'' tetrahedra and octahedra can be constructed from copies of tetrahedra and octahedra, see Figure~\ref{tetra_oct_explode}. We observe that the assembly rule in the assemblies in Figure~\ref{tetra_oct_interlocking_explode} then follows the logic of the tetrahedral-octahedral~honeycomb.

\section{Interlocking Blocks in the Tetroctahedrille}
\label{section:NewBlocks}
In this section, we first introduce examples of blocks with interlocking properties, that are constructed by exploiting the tetroctahedrille. Next, we give possible topological interlocking assemblies for each presented block.

In order to simplify the construction of the topological interlocking blocks that can be created by arranging copies of tetrahedra and octahedra complying with rules R1 and R2, we define the tetrahedral-octahedral~decomposition of a three-dimensional body.
\begin{definition}
    Let $X$ be a block in $\mathbb{R}^3$. Furthermore, let $m,n$ be natural numbers.  If there exist copies of regular tetrahedra $T_1,\ldots,T_m$ and regular octahedra $O_1,\ldots, O_n$ with all edges of same length so that 
    \[
    X=T_1\cup\ldots\cup T_m \cup O_1\cup\ldots\cup O_m,  
    \]
    where any two of the above copies of the Platonic solid only intersect in their boundary, 
    then $(T_1,\ldots,T_m,O_1,\ldots,O_n)$ is called a \emph{tetrahedral-octahedral~decomposition} of $X.$
\end{definition}
The polyhedra that are the focus of this investigation, i.e. the polyhedra that can be cut out from the tetroctahedrille, have a tetrahedral-octahedral~decomposition. 
For example, up to isometry, a tetrahedral-octahedral~decomposition of the right most block shown in Figure~\ref{octet} is  given by 
\[
(\{0,v_1,v_2,v_3\},\{v_1,v_2,v_3,v_1+v_2,v_1+v_3,v_2+v_3\}).
\]

Note, the tetroctahedrille yields a rich combinatorial structure, since the incidence structure of the polyhedra that are contained in the honeycomb can be examined further, see Remark~\ref{simplSurf}.
\begin{remark} \label{simplSurf}
Since it is possible to triangulate a block that admits an octahedral tetrahedral decomposition so that all faces of the triangulation form equilateral triangles with the same edge lengths, the incidence structure of vertices, edges and faces of the blocks can be the focus of further investigations. Functions and implementations to study and manipulate the combinatorics of the incidence structure of the described blocks can be found in the GAP4 package SimplicialSurfaces \citep{SimpSurf21}.
\end{remark}

\subsection{Construction of the Blocks}
Here, we exploit the notion of a tetrahedral-octahedral~decomposition to construct blocks that allow topological interlocking assemblies and further present some additional properties. Note, we will use the tetrahedra and octahedra presented in Remark~\ref{v1v2v3} for these constructions.
\subsubsection{The Kitten}
\begin{figure}[H]
\begin{minipage}{.3\textwidth}
    \centering
    \includegraphics[height=4cm,width=4cm]{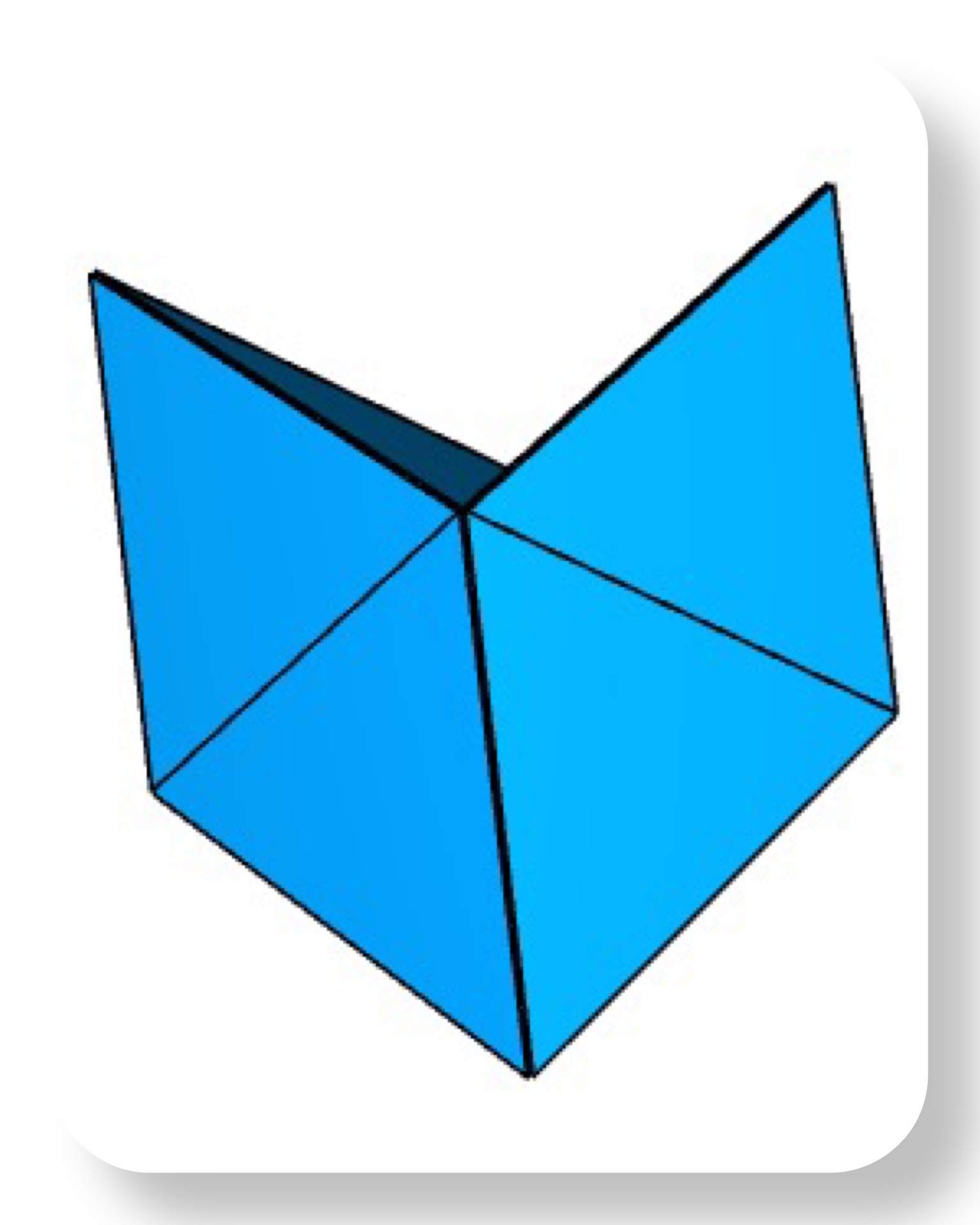}
\end{minipage}
\begin{minipage}{0.5cm}
    
\end{minipage}
\begin{minipage}{.3\textwidth}
    \centering
    \includegraphics[height=4cm,width=4cm]{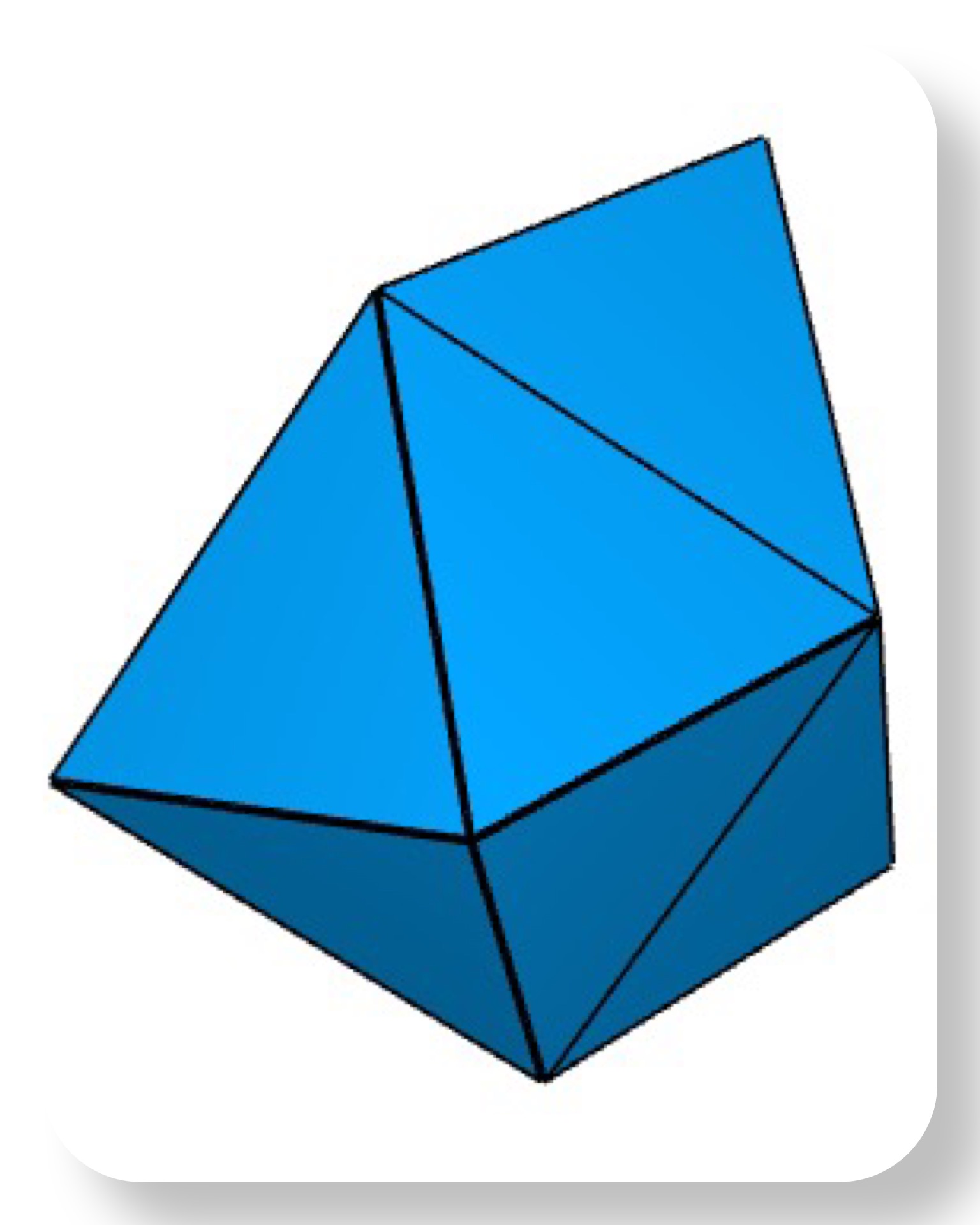}
\end{minipage}
\begin{minipage}{0.5cm}
    
\end{minipage}
\begin{minipage}{.3\textwidth}
    \centering
    \includegraphics[height=4cm,width=4cm]{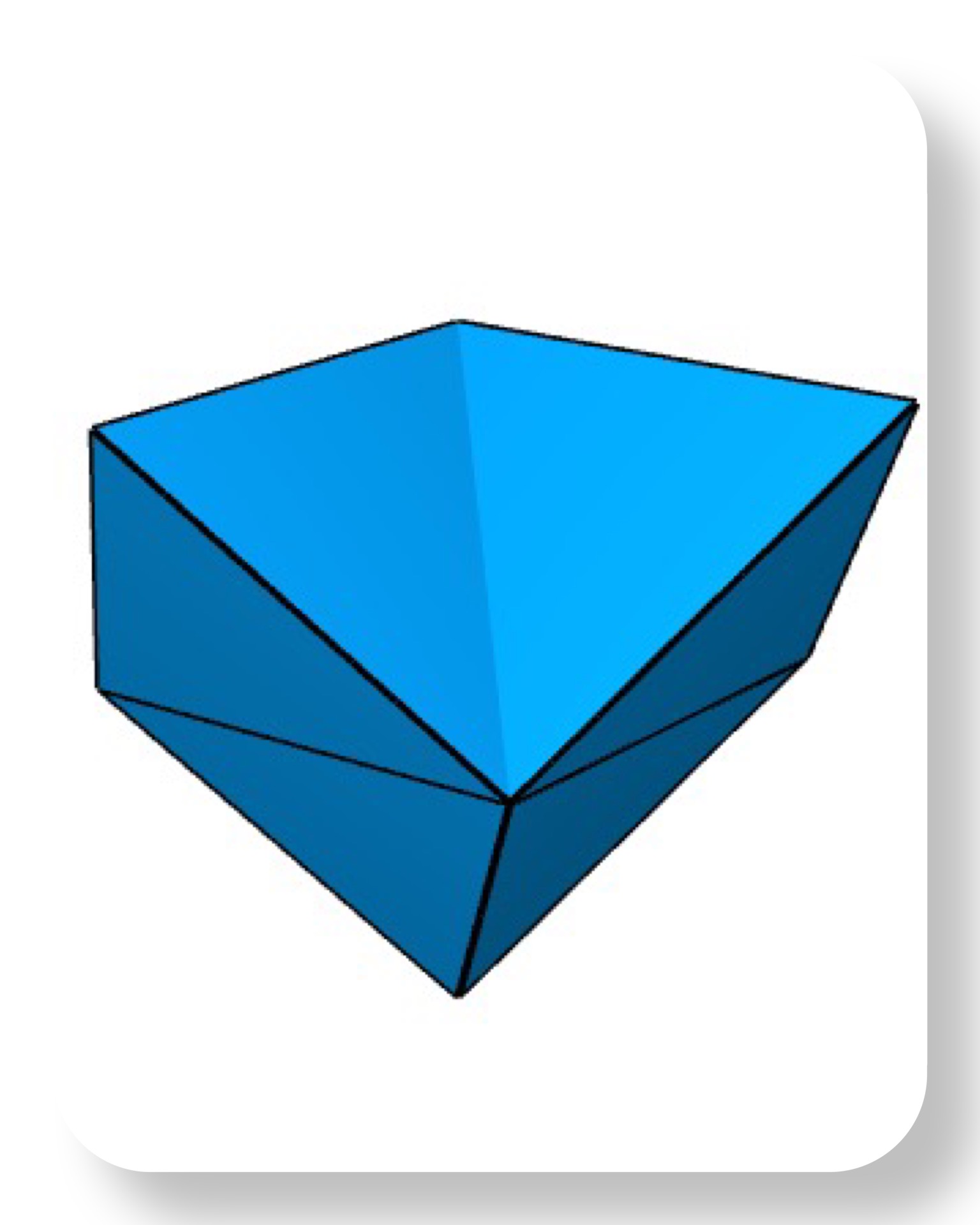}
\end{minipage}
\caption{Various views of the constructed kitten, consisting of one octahedron and two tetrahedra}
\label{kitten}
\end{figure}

In order to construct the first block, two tetrahedra are attached to an octahedron so that all three polyhedra  share a common edge, see Figure~\ref{kitten}.
 
\begin{figure}[H]
    \centering
    \includegraphics[scale=1,viewport=11cm 21.5cm 0cm 26.5cm]{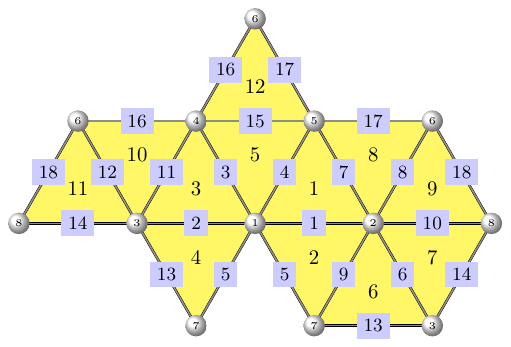}
\caption{A folding plan of the kitten}
\label{foldingplankitten}
\end{figure}

\begin{figure}[H]
\begin{minipage}{.5\textwidth}
    \centering
    \includegraphics[height=5cm,width=5cm]{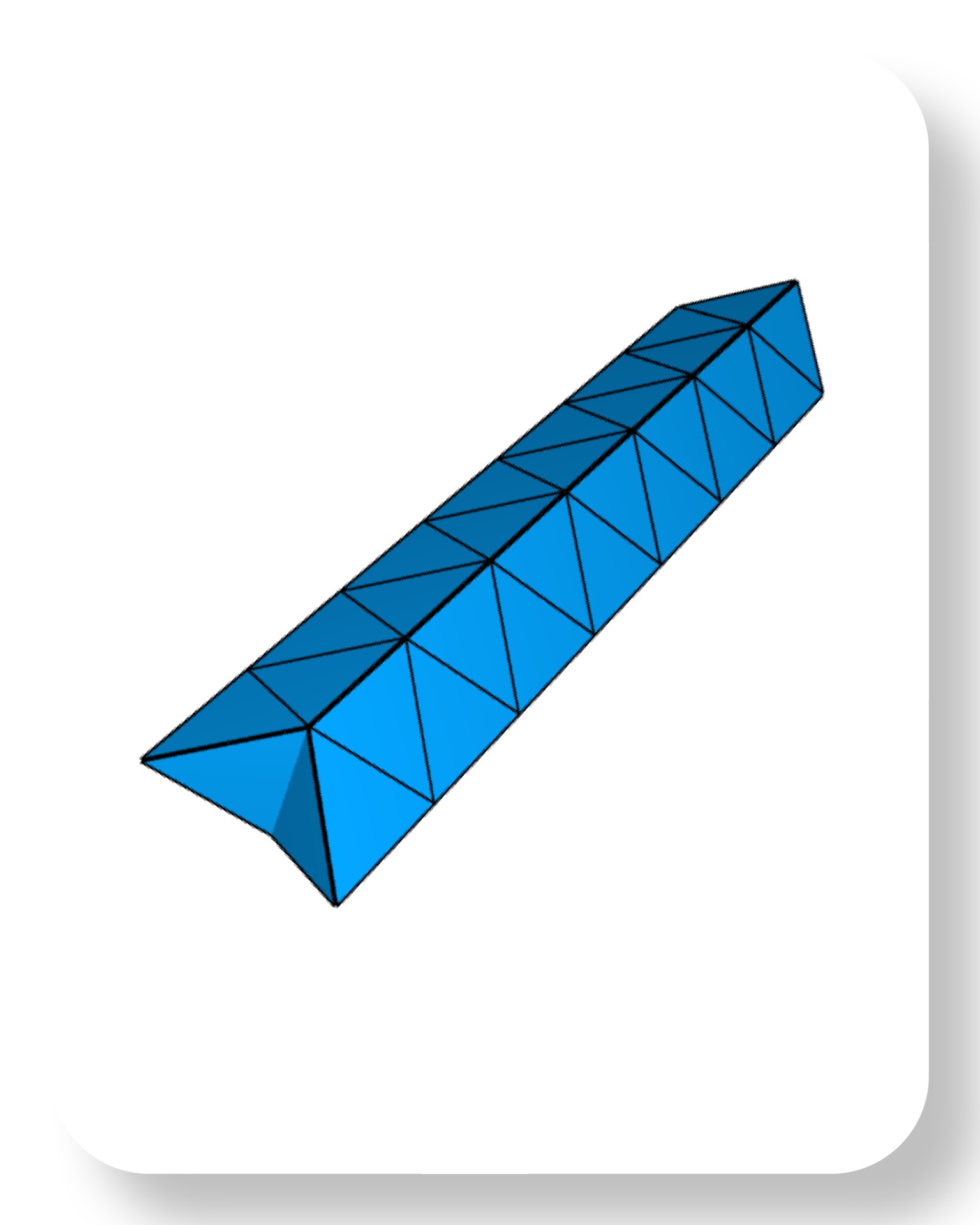}
\end{minipage}
\begin{minipage}{0.5cm}
    
\end{minipage}
\begin{minipage}{.5\textwidth}
    \centering
    \includegraphics[height=5cm,width=5cm]{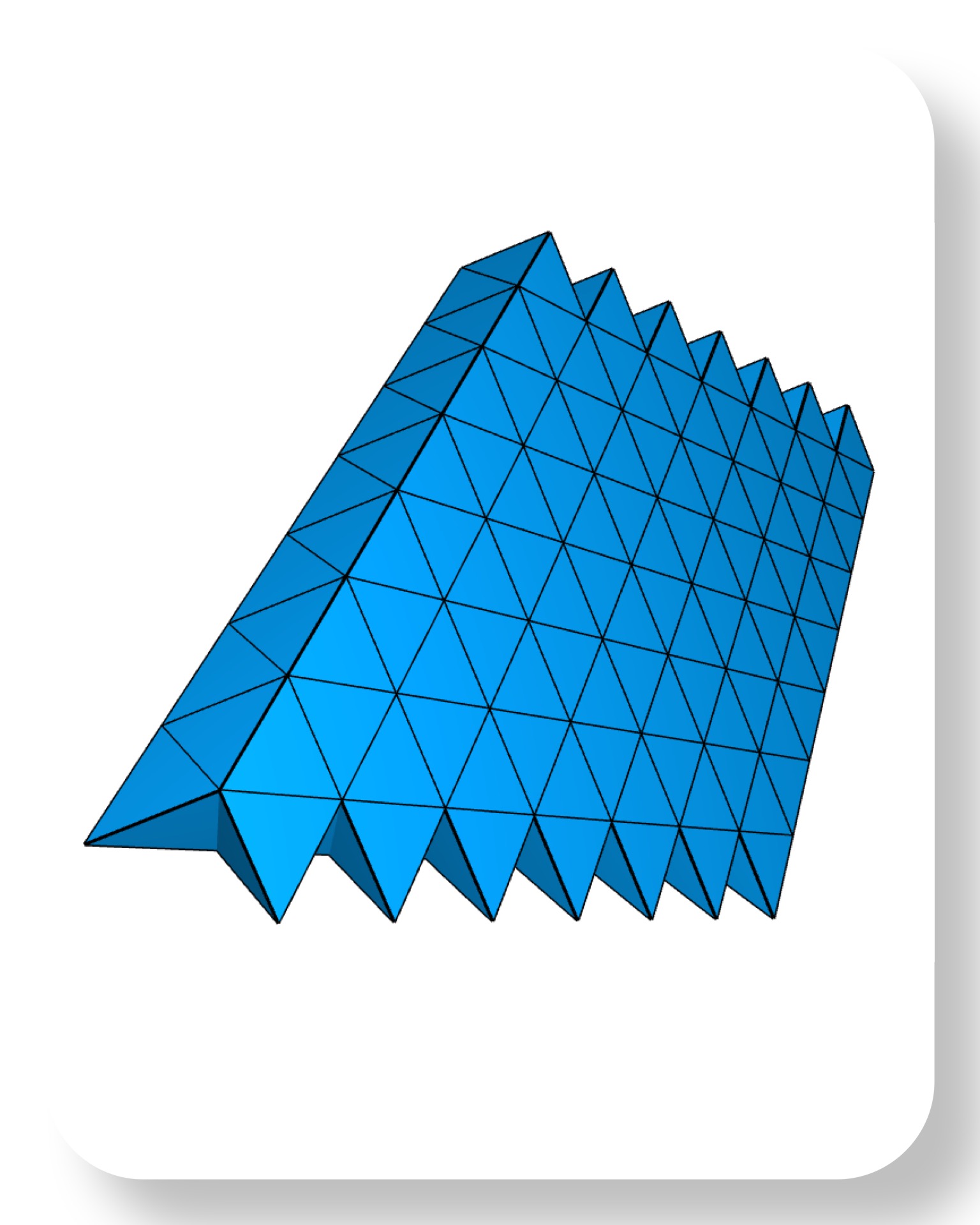}
\end{minipage}
\caption{Assembling copies of the kitten to form a strip (left) and a plane (right)}
\label{kittenspacefilling}
\end{figure}

This block is called \emph{kitten}. Up to isometry, the tetrahedral-octahedral~decomposition of the kitten is given by $(T_1,T_2,O).$
This block can be described by the following ordered list of coordinates in $\mathbb{R}^3$ representing the eight vertices of the polyhedron:
 \begin{align*}
\left[ \begin{pmatrix} 0\\1\\1\end{pmatrix},
\begin{pmatrix} 1\\0\\1\end{pmatrix},
\begin{pmatrix} 1\\1\\0\end{pmatrix},
\begin{pmatrix} 1\\2\\1\end{pmatrix},
\begin{pmatrix} 1\\1\\2\end{pmatrix},
\begin{pmatrix} 2\\1\\1\end{pmatrix},
\begin{pmatrix} 0\\0\\0\end{pmatrix},
\begin{pmatrix} 2\\0\\0\end{pmatrix}\right].  
 \end{align*}
For the purpose of this paper, the vertices of the resulting triangulation are referred to by their position in the list of coordinates. 
Moreover, the underlying incidence structure of the triangulation shown in Figure~\ref{kitten} is determined by the following list of faces given by their incident vertices:
\begin{align*}
[ [ 1, 2, 5 ], [ 1, 2, 7 ], [ 1, 3, 4 ], [ 1, 3, 7 ], [ 1, 4, 5 ], 
  [ 2, 3, 7 ],\\
  [ 2, 3, 8 ], [ 2, 5, 6 ], [ 2, 6, 8 ], [ 3, 4, 6 ], 
  [ 3, 6, 8 ], [ 4, 5, 6 ] ].
\end{align*}
The symmetry group of the block and the automorphism group of the underlying incidence geometry of the given triangulation are both isomorphic to $C_2 \times C_2 $.
Figure~\ref{foldingplankitten} shows a folding plan of the corresponding triangulation of this block. 
The kitten is an example of a polyhedron that allows a space-filling, i.e. we can fill the three-dimensional space without gaps by assembling copies of this block. This can be seen by positioning copies of the kitten by using the translation $v_1$ to create a strip and then assembling copies of the resulting strip to create the desired space filling, see Figure \ref{kittenspacefilling}.
\subsubsection{The Cushion}
Next, we construct \emph{the cushion}. This block can be constructed by 
attaching four tetrahedra to an octahedron so that for each tetrahedron there is exactly one tetrahedron sharing a common edge with the first tetrahedron.
Then, up to isometry, the tetrahedral-octahedral~decomposition of $X$ is given by the following subsets of the face-centred cubic lattice:
    \begin{align*}
&(T_1,T_2,v_1+T_1,v_1+T_2,O).
    \end{align*}

\begin{figure}[H]
\begin{minipage}{.3\textwidth}
    \centering
    \includegraphics[height=4cm]{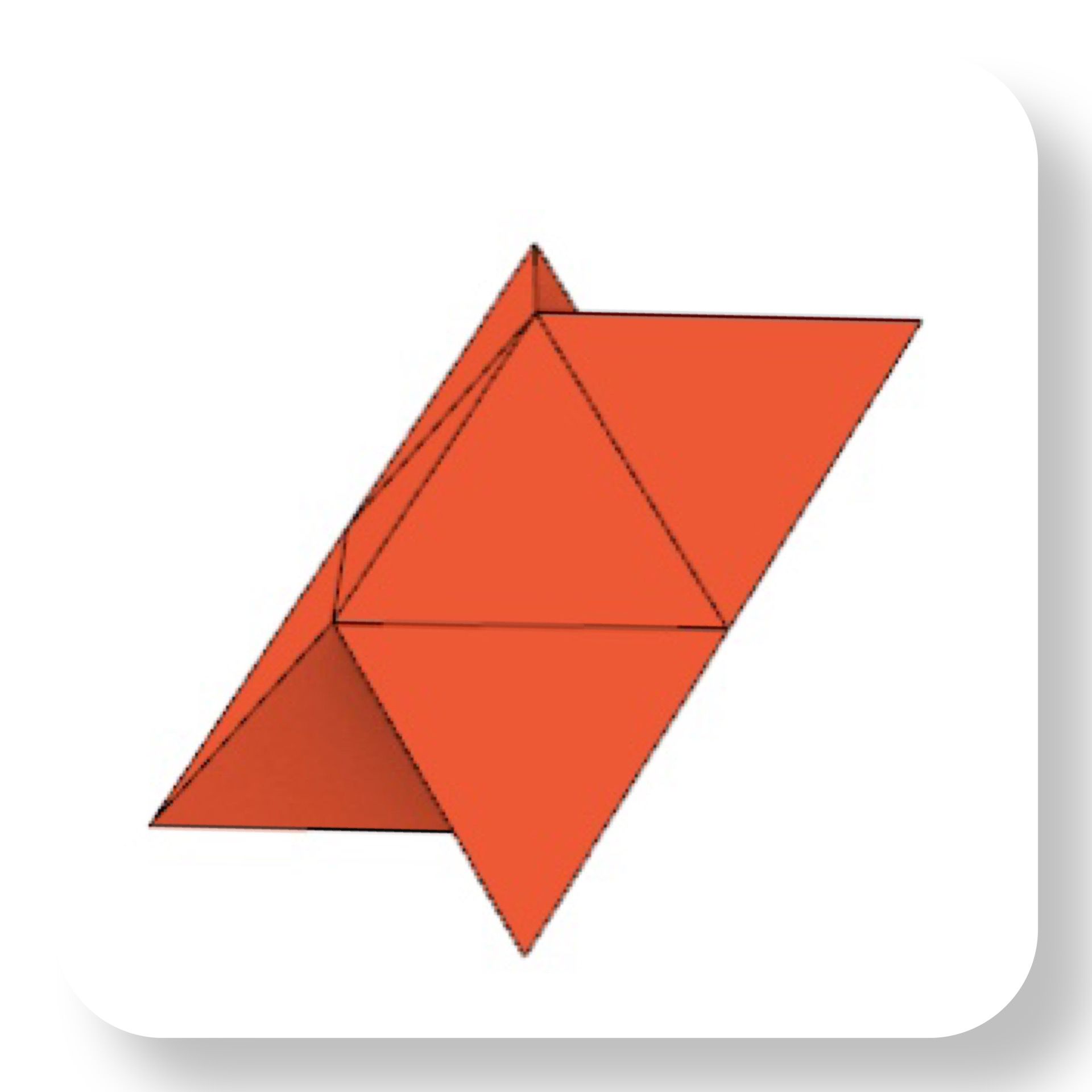}
\end{minipage}
\begin{minipage}{0.5cm}
    
\end{minipage}
\begin{minipage}{.3\textwidth}
    \centering
    \includegraphics[height=4cm]{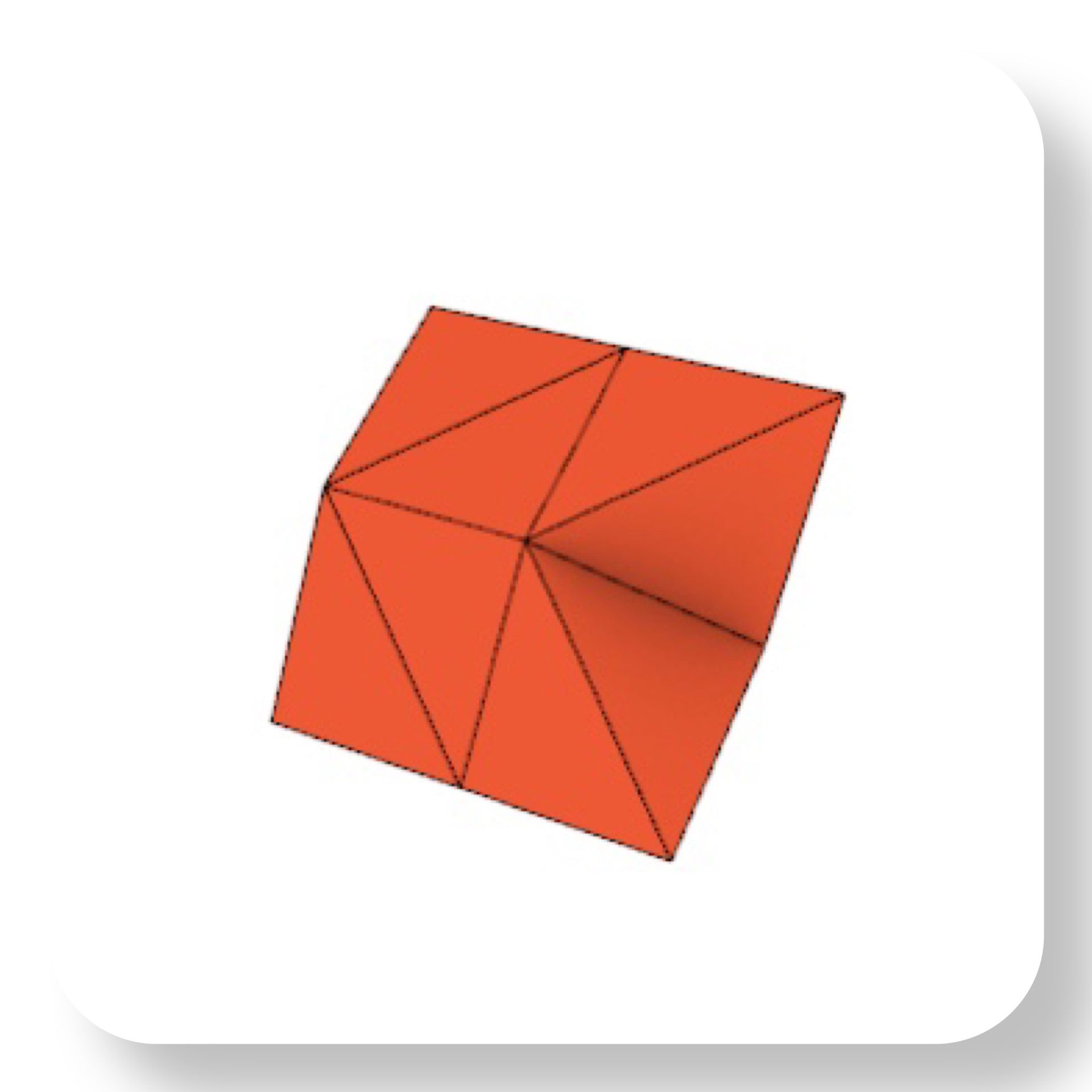}
\end{minipage}
\begin{minipage}{0.5cm}
    
\end{minipage}
\begin{minipage}{.3\textwidth}
    \centering
    \includegraphics[height=4cm]{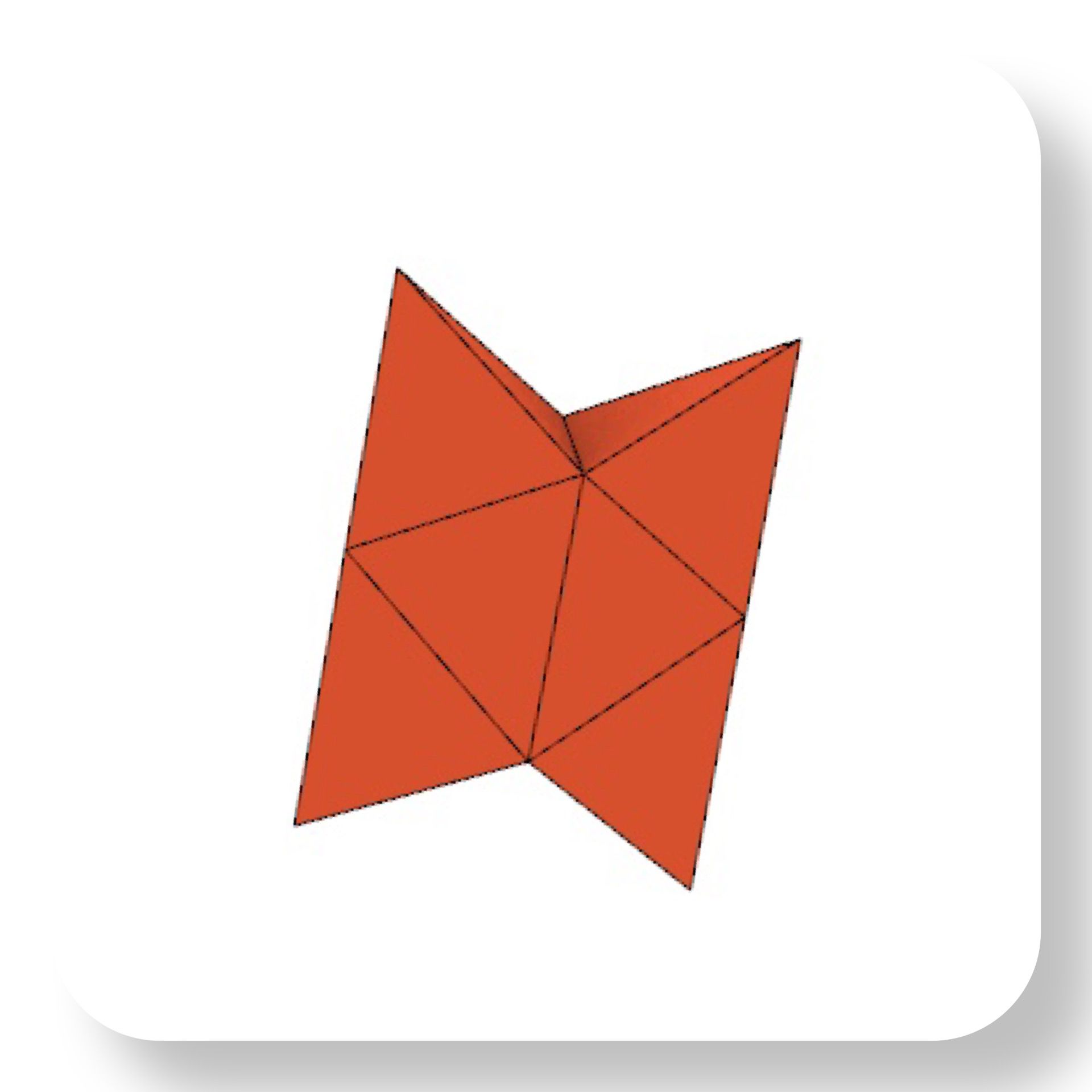}
\end{minipage}
\caption{Various views of the constructed cushion consisting of one octahedron and four tetrahedra}
\label{cushion}
\end{figure}

Here, we present the ordered list of vertex coordinates of the block shown in Figure~\ref{cushion}:
 \begin{align*}
\left[ \begin{pmatrix} 0\\1\\1\end{pmatrix},
\begin{pmatrix} 1\\0\\1\end{pmatrix},
\begin{pmatrix} 1\\1\\0\end{pmatrix},
\begin{pmatrix} 1\\2\\1\end{pmatrix},
\begin{pmatrix} 1\\1\\2\end{pmatrix},
\begin{pmatrix} 2\\1\\1\end{pmatrix},
\begin{pmatrix} 0\\0\\0\end{pmatrix},
\begin{pmatrix} 2\\0\\0\end{pmatrix},
\begin{pmatrix} 0\\2\\2\end{pmatrix},
\begin{pmatrix} 2\\2\\2\end{pmatrix}\right].  
 \end{align*}
Furthermore, the vertices of faces of the corresponding triangulation are given by 
\begin{align*}
[ &[ 1, 2, 5 ], [ 1, 2, 7 ], [ 1, 3, 4 ], [ 1, 3, 7 ], [ 1, 4, 9 ], 
  [ 1, 5, 9 ],
  [ 2, 3, 7 ], [ 2, 3, 8 ],\\ &[ 2, 5, 6 ], [ 2, 6, 8 ], 
  [ 3, 4, 6 ], [ 3, 6, 8 ], 
  [ 4, 5, 9 ], [ 4, 5, 10 ], [ 4, 6, 10 ], 
  [ 5, 6, 10 ] ].
\end{align*}

\begin{figure}[H]
    \centering    \includegraphics[scale=0.9,viewport=12cm 23cm 0cm 27cm]{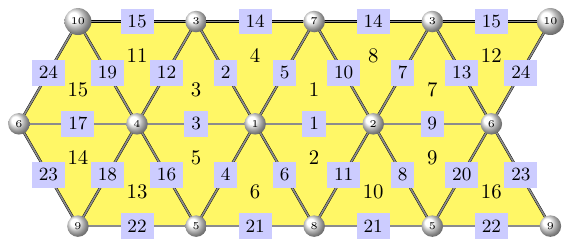}
\caption{A folding plan of a triangulation of the cushion}
\label{foldingplancushion}

\end{figure}

\begin{figure}[H]
\begin{minipage}{.3\textwidth}
    \centering
    \includegraphics[height=4cm]{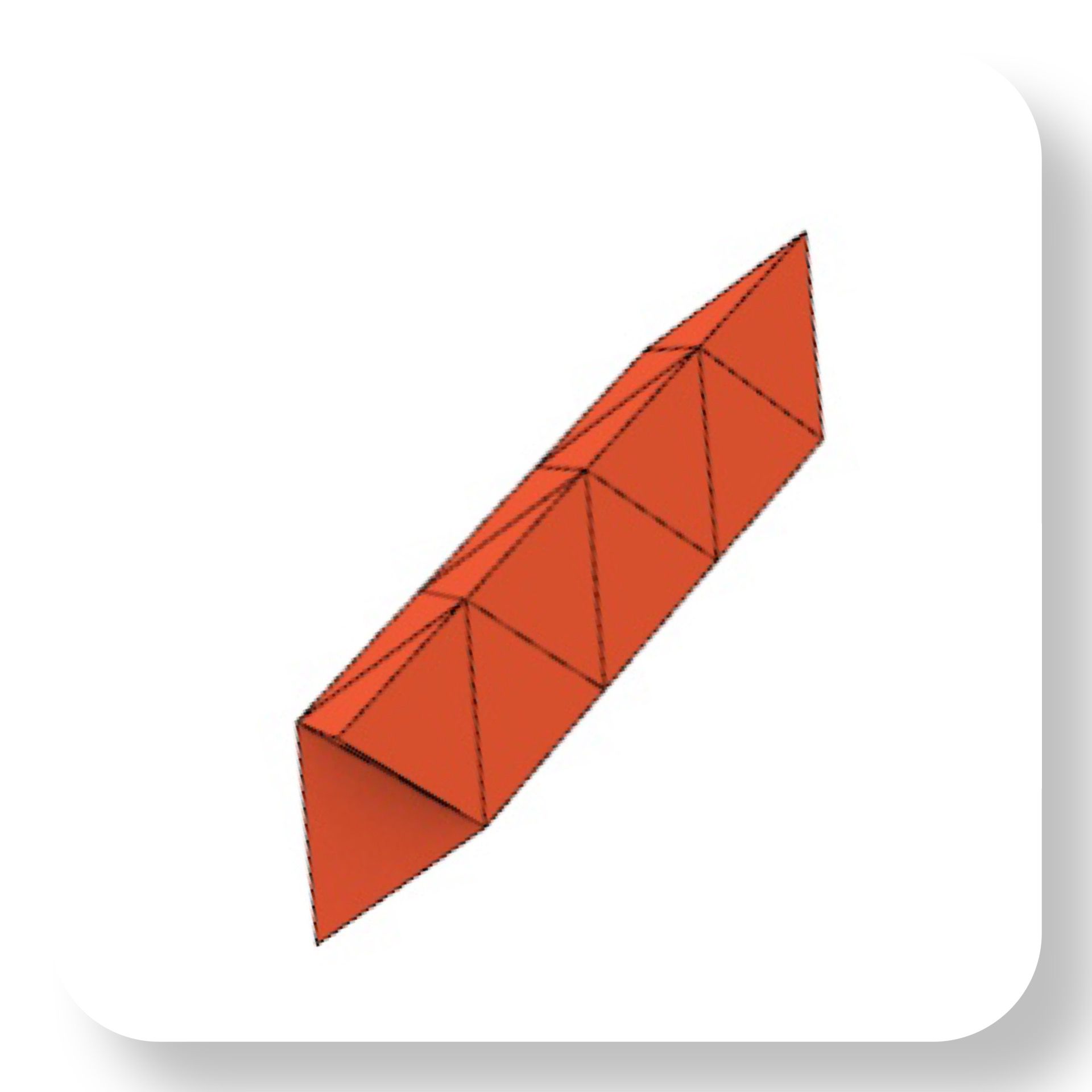}
    \label{strip21}
\end{minipage}
\begin{minipage}{0.5cm}
    
\end{minipage}
\begin{minipage}{.3\textwidth}
    \centering
    \includegraphics[height=4cm]{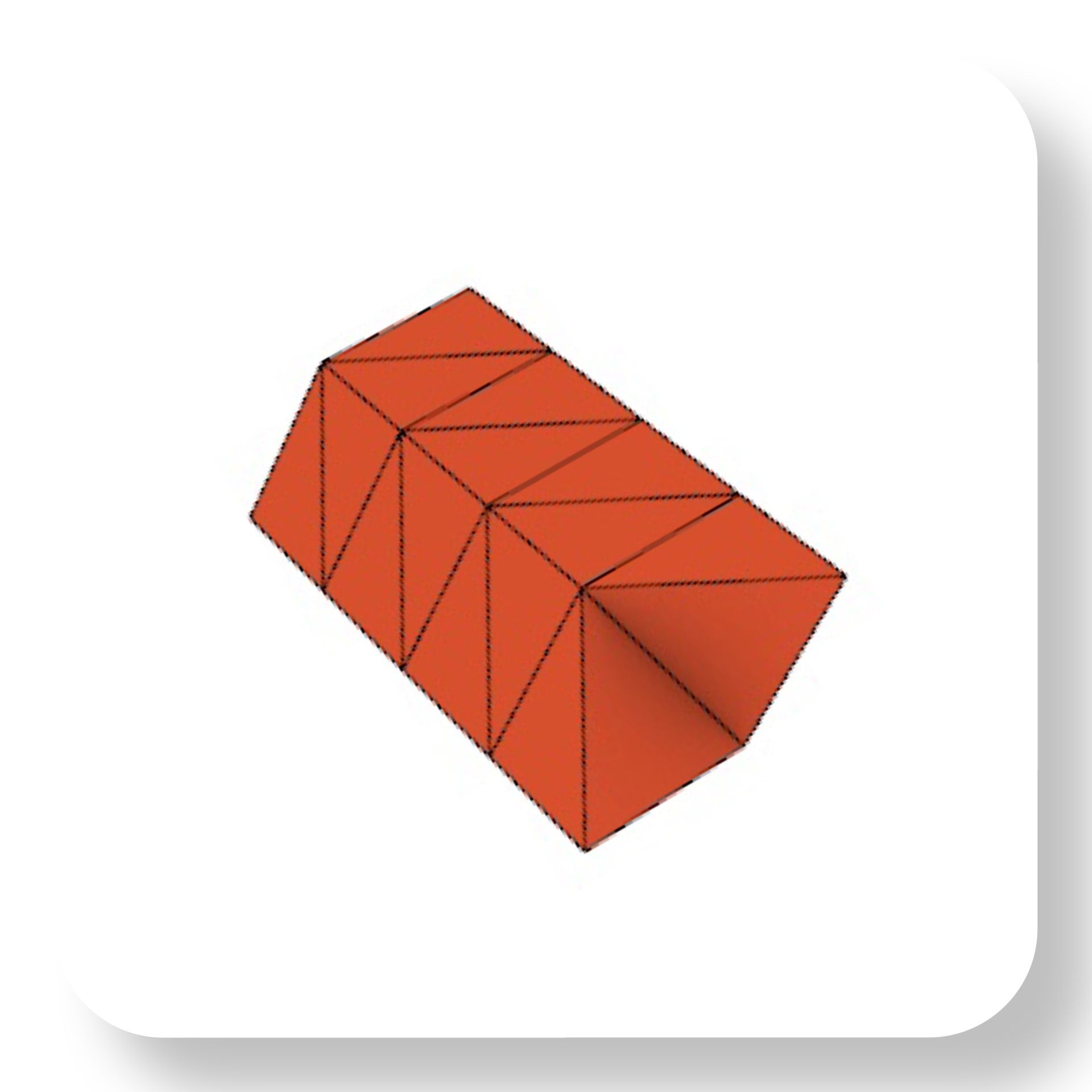}
    \label{strip22}
\end{minipage}
\begin{minipage}{0.5cm}
    
\end{minipage}
\begin{minipage}{.3\textwidth}
    \centering
    \includegraphics[height=4cm]{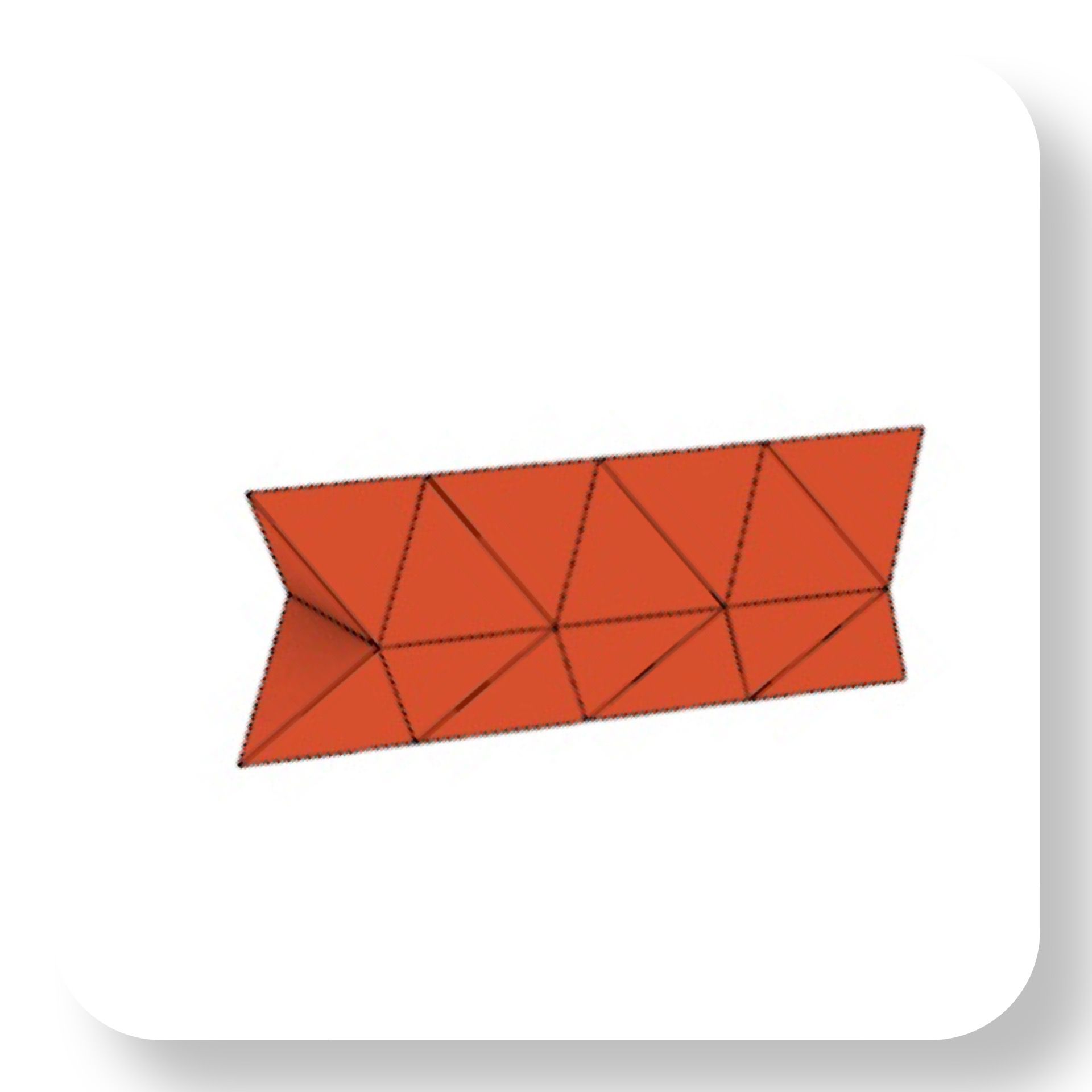}
    \label{strip32}
\end{minipage}
\caption{Various views of the 3-cushion consisting of three octahedra and $8$ tetrahedra}
\label{3cushion}
\end{figure}

This block as a rigid three-dimensional body has a symmetry group that is isomorphic to $C_2 \times C_2$. By further examining the underlying incidence structure of this block, we see that the automorphism group is isomorphic to the described symmetry group.
A folding plan of the above triangulation is presented in Figure~\ref{foldingplancushion}.
Next we formulate a generalised construction of the cushion to construct interlocking blocks, that allow similar topological interlocking assemblies. The surfaces of these blocks can be viewed as surfaces that result from a continuous deformation of the surface of the cushion. 
\begin{remark}\label{rem1}
    Let $X$ be the above cushion.
    For a positive integer $n,$ \emph{the $n$-cushion} can be defined by the following tetrahedral-octahedral~decomposition: 
    \begin{align*}
    &(T_0,\ldots,nv_1+T_1,T_2,\ldots, nv_1+T_2,O,\ldots,(n-1)v_1+O).
\end{align*}    
\end{remark}
Note, the $n$-cushion can be decomposed into $2(n+1)$ tetrahedra and $n$ octahedra.
The $3$-cushion that results from the above construction is illustrated in Figure~\ref{3cushion}.

\subsubsection{The Shuriken}
\begin{figure}[H]
\begin{minipage}{.3\textwidth}
    \centering
    \includegraphics[height=4cm]{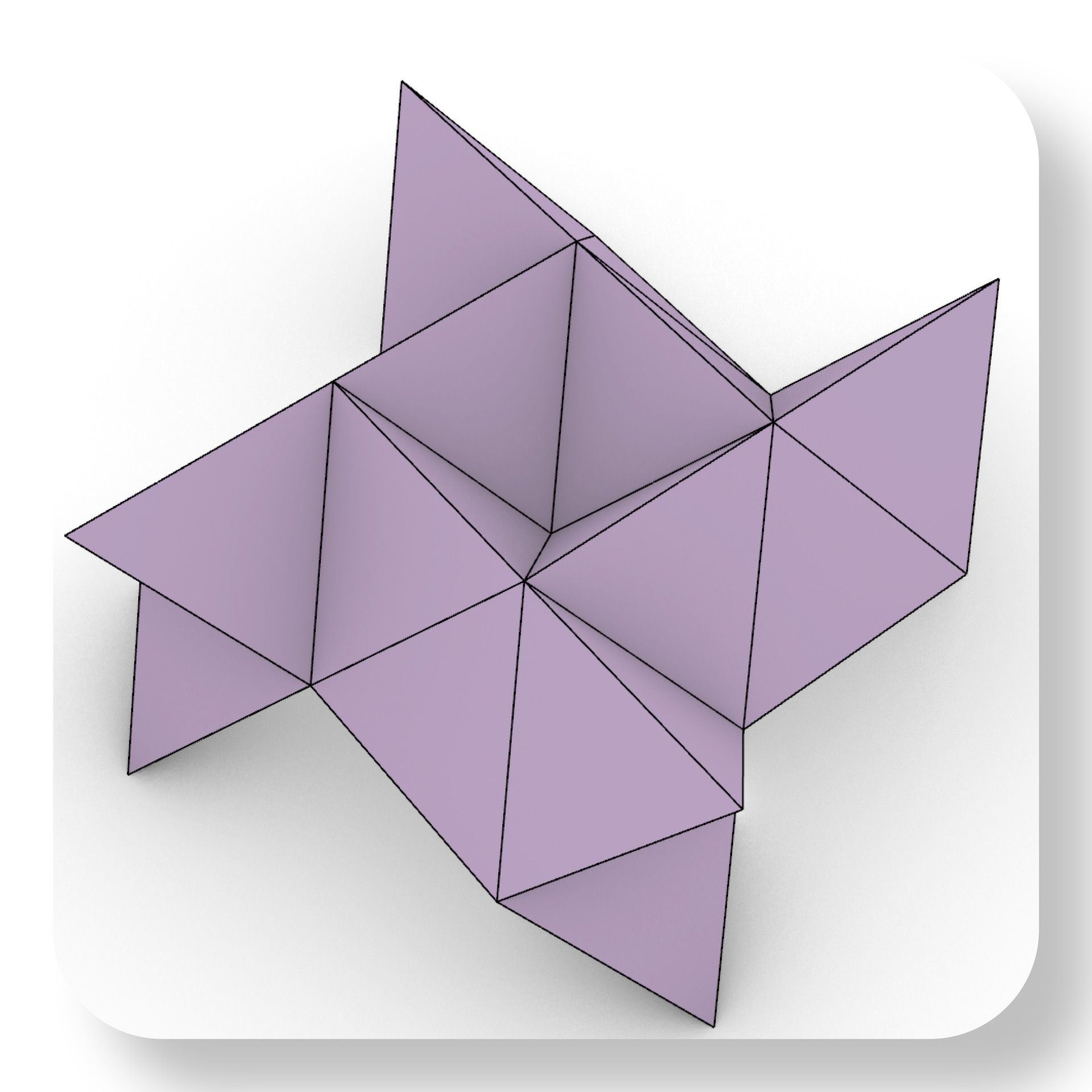}
\end{minipage}
\begin{minipage}{0.5cm}
    
\end{minipage}
\begin{minipage}{.3\textwidth}
    \centering
    \includegraphics[height=4cm]{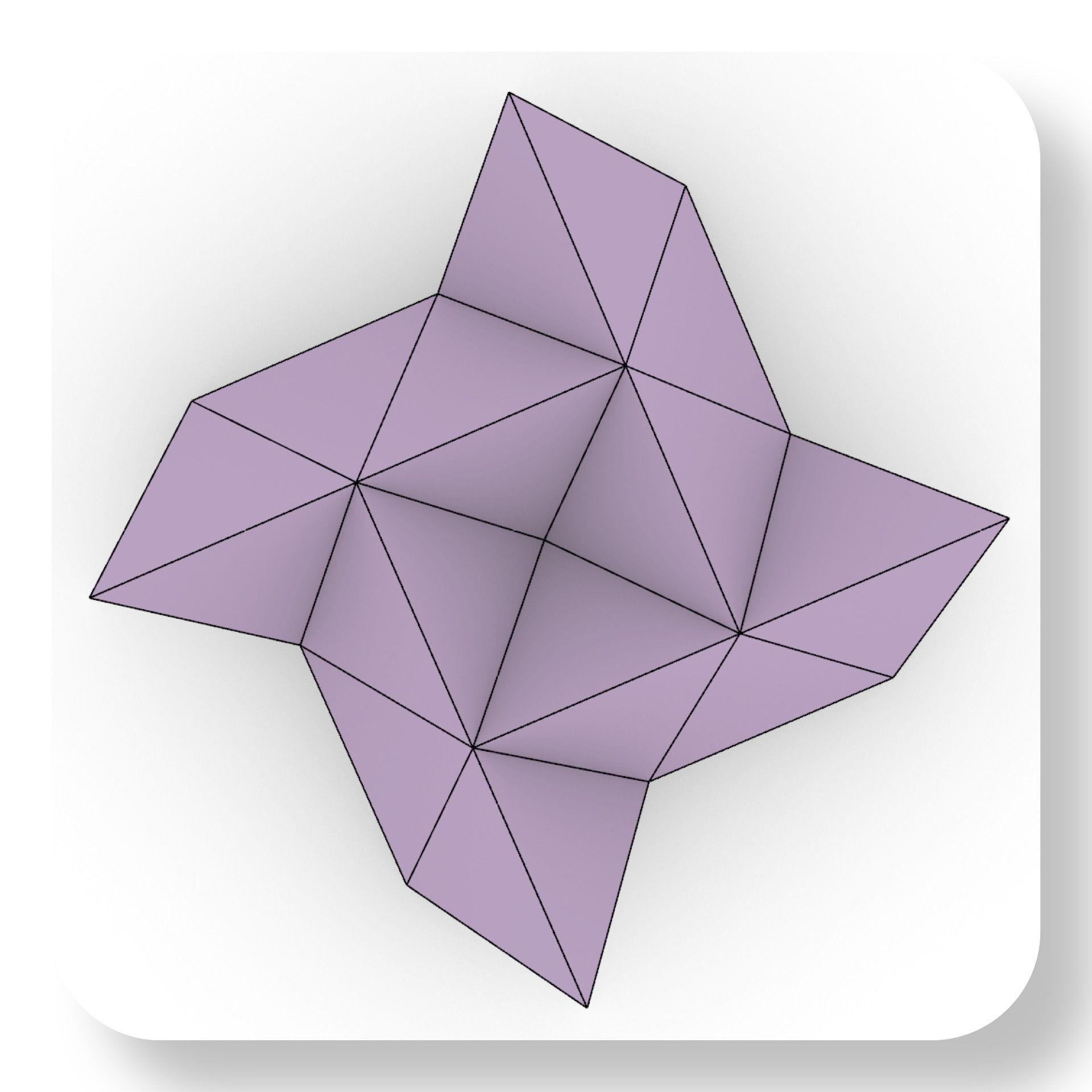}
\end{minipage}
\begin{minipage}{0.5cm}
    
\end{minipage}
\begin{minipage}{.3\textwidth}
    \centering
    \includegraphics[height=4cm]{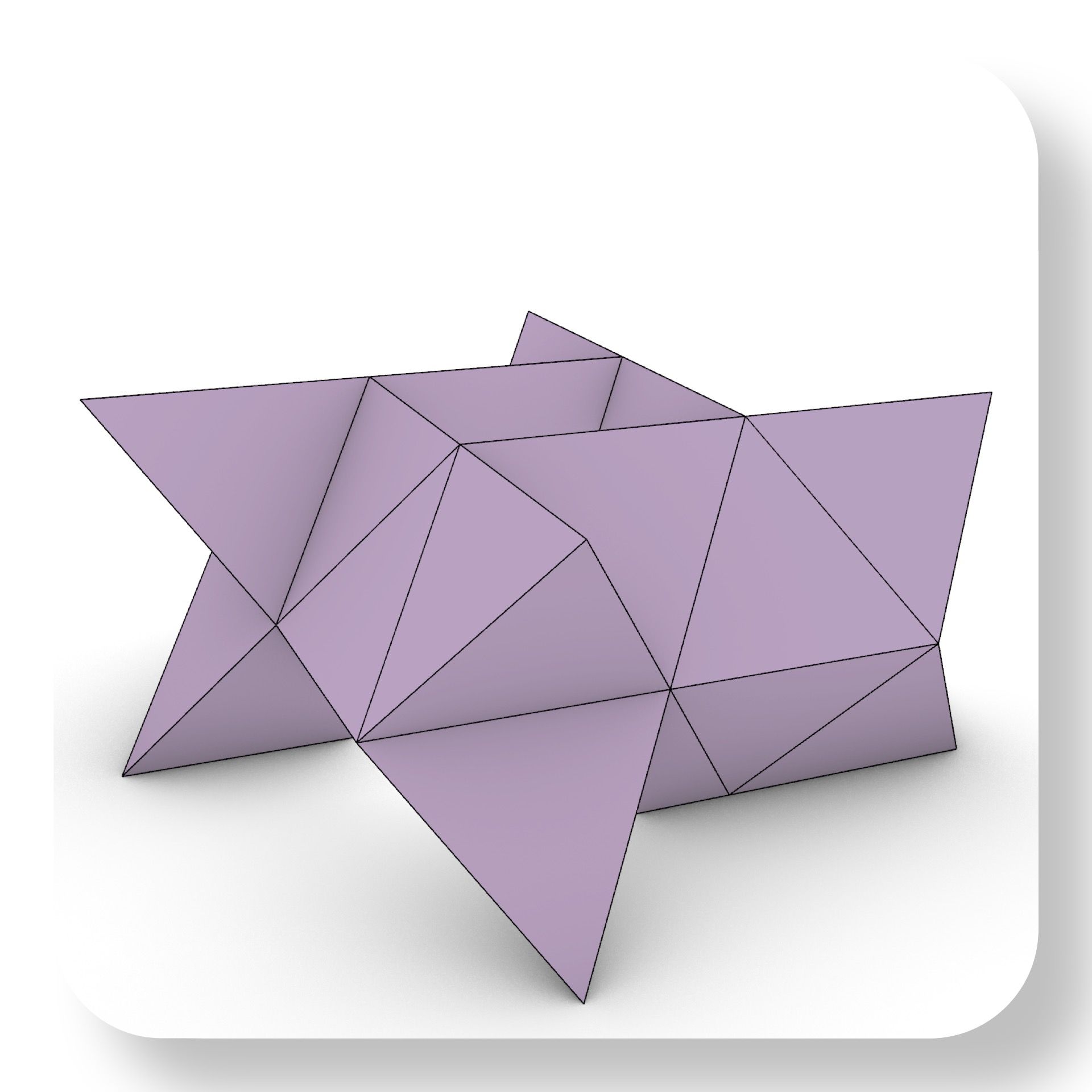}
\end{minipage}
\caption{Various views of the constructed shuriken consisting of four octahedra and $16$ tetrahedra}
\label{shuriken}
\end{figure}

In this section, we present an interlocking block that is called \emph{the shuriken}.
The shuriken is a block that can be constructed by combining 16 tetrahedra and 4 octahedra in the following way: 
Let $v_1,v_2,v_3$ be defined as in Remark~\ref{v1v2v3}. Then, we obtain the shuriken by the following tetrahedral-octahedral decomposition:
    \begin{align*}
    (&T_1,v_1+T_1,v_1+(v_2-v_3)+T_1, 2v_1+(v_2-v_3)+T_1,\\
    &T_2,v_1+T_2,v_1+(v_2-v_3)+T_2,2v_1+(v_2-v_3)+T_2,\\
    &(v_2-v_3)+T_3,2(v_2-v_3)+T_3,v_1+2(v_2-v_3)+T_3,v_1+3(v_2-v_3)+T_3\\
    &(v_2-v_3)+T_4,2(v_2-v_3)+T_4,v_1+2(v_2-v_3)+T_4,v_1+3(v_2-v_3)+T_4\\
    &O,v_1+O,(v_2-v_3)+O,v_1+(v_2-v_3)+O).
    \end{align*}
Placing the polyhedra in Euclidean space as described in the above tetrahedral-octahedral~symbol gives rise to the structure shown in Figure~\ref{shuriken}.

\begin{figure}[H]
\begin{minipage}{.3\textwidth}
    \centering
    \includegraphics[height=4cm]{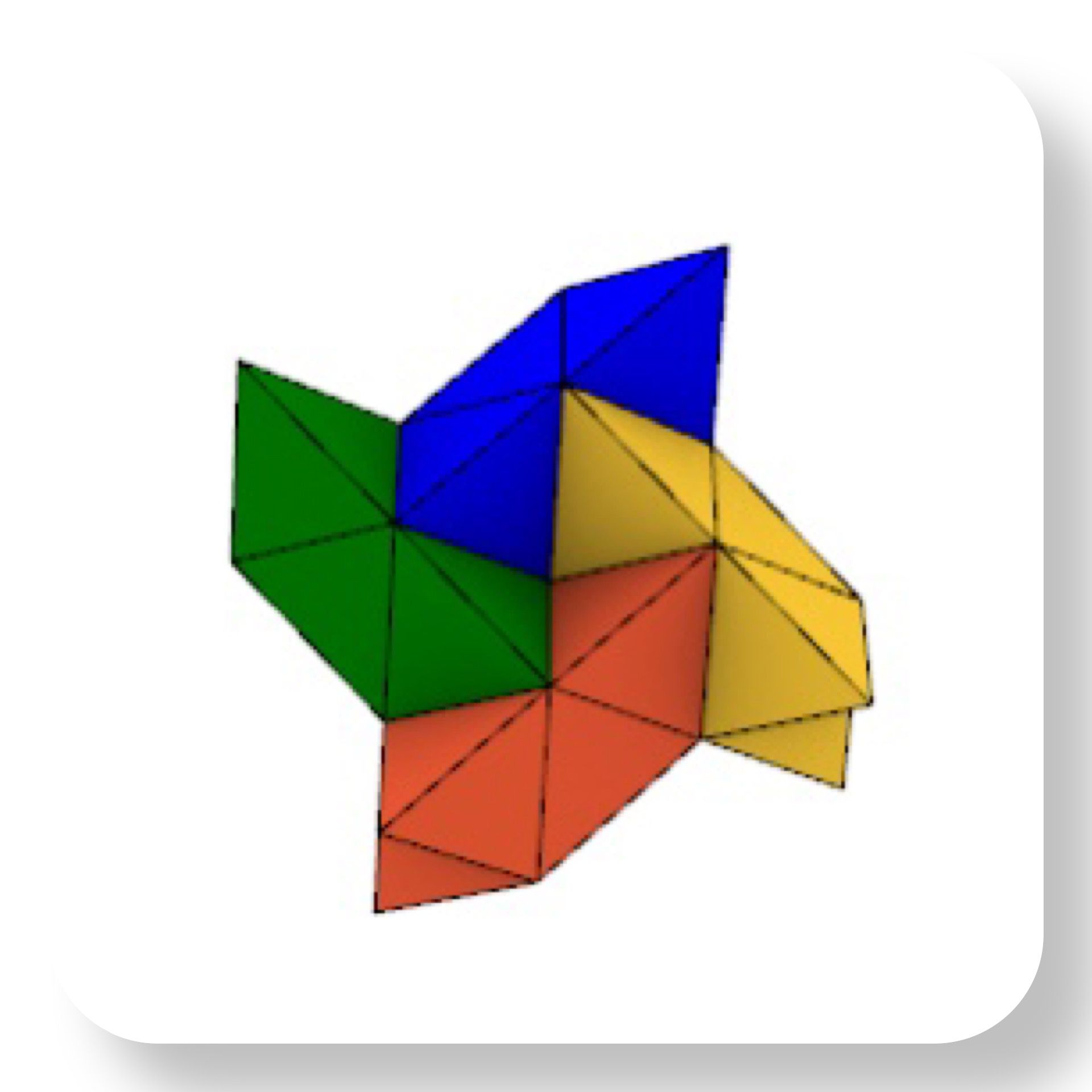}
\end{minipage}
\begin{minipage}{0.5cm}
    
\end{minipage}
\begin{minipage}{.3\textwidth}
    \centering
    \includegraphics[height=4cm]{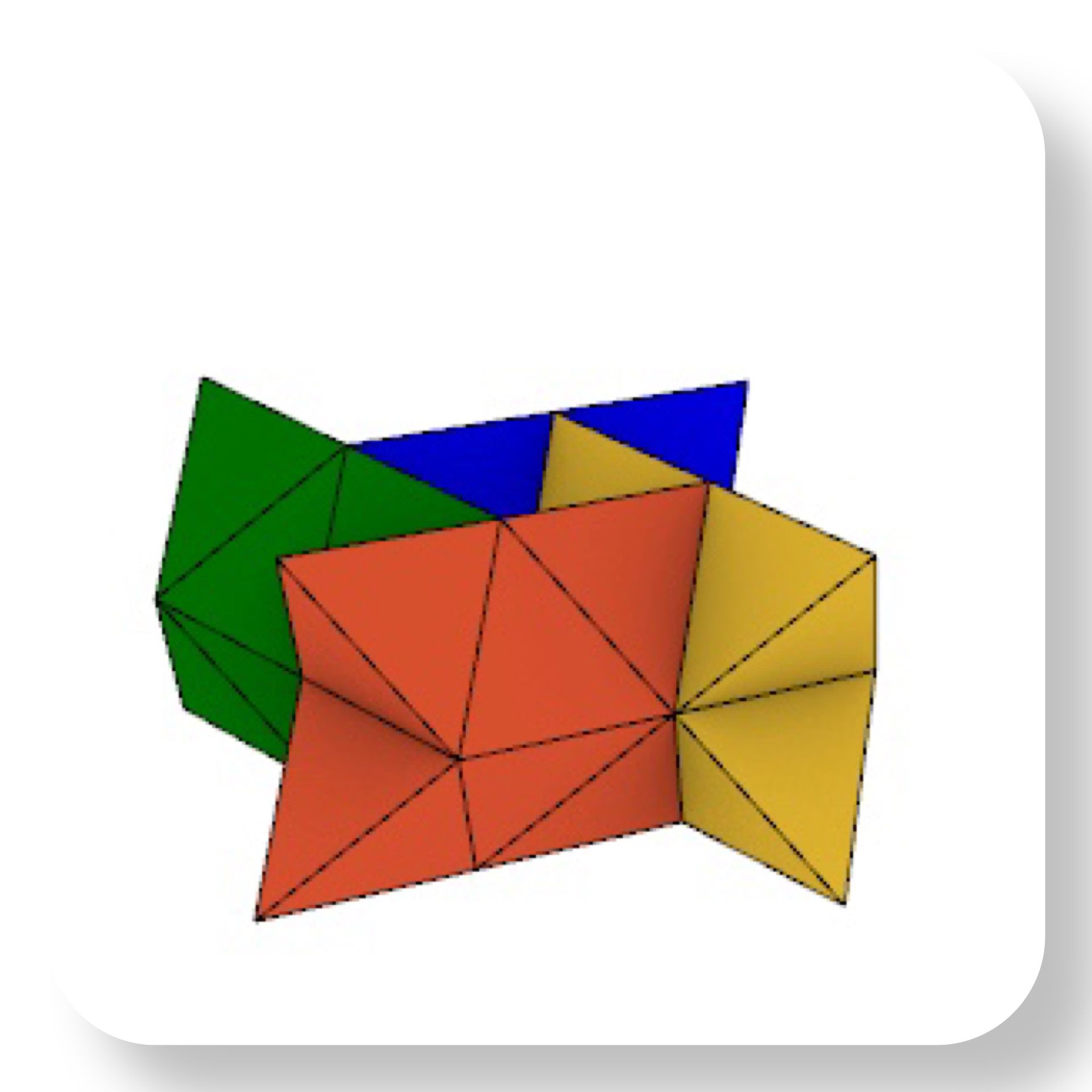}

\end{minipage}
\begin{minipage}{0.5cm}
    
\end{minipage}
\begin{minipage}{.3\textwidth}
    \centering
    \includegraphics[height=4cm]{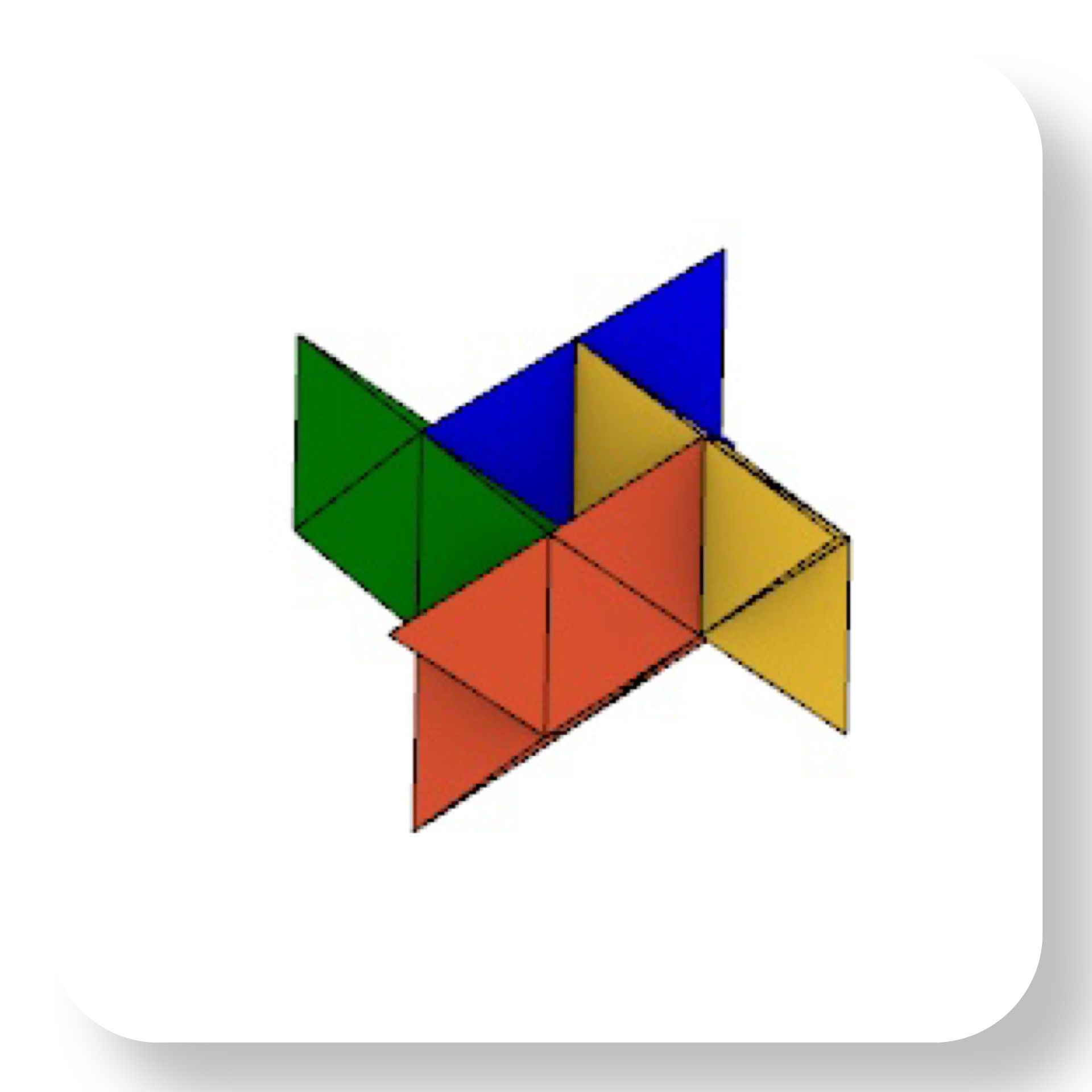}
\end{minipage}
\caption{Assembly of four cushions forming one shuriken}
\label{arrangementcushion}
\end{figure}

An alternative way to construct this block is to place 4 cushions in the Euclidean 3-space so that they all share exactly one common vertex and the resulting assembly still satisfies the rules of the tetrahedral-octahedral~honeycomb, see Figure~\ref{arrangementcushion}.

\begin{figure}[H]
\begin{minipage}{.3\textwidth}
    \centering
    \includegraphics[height=4cm]{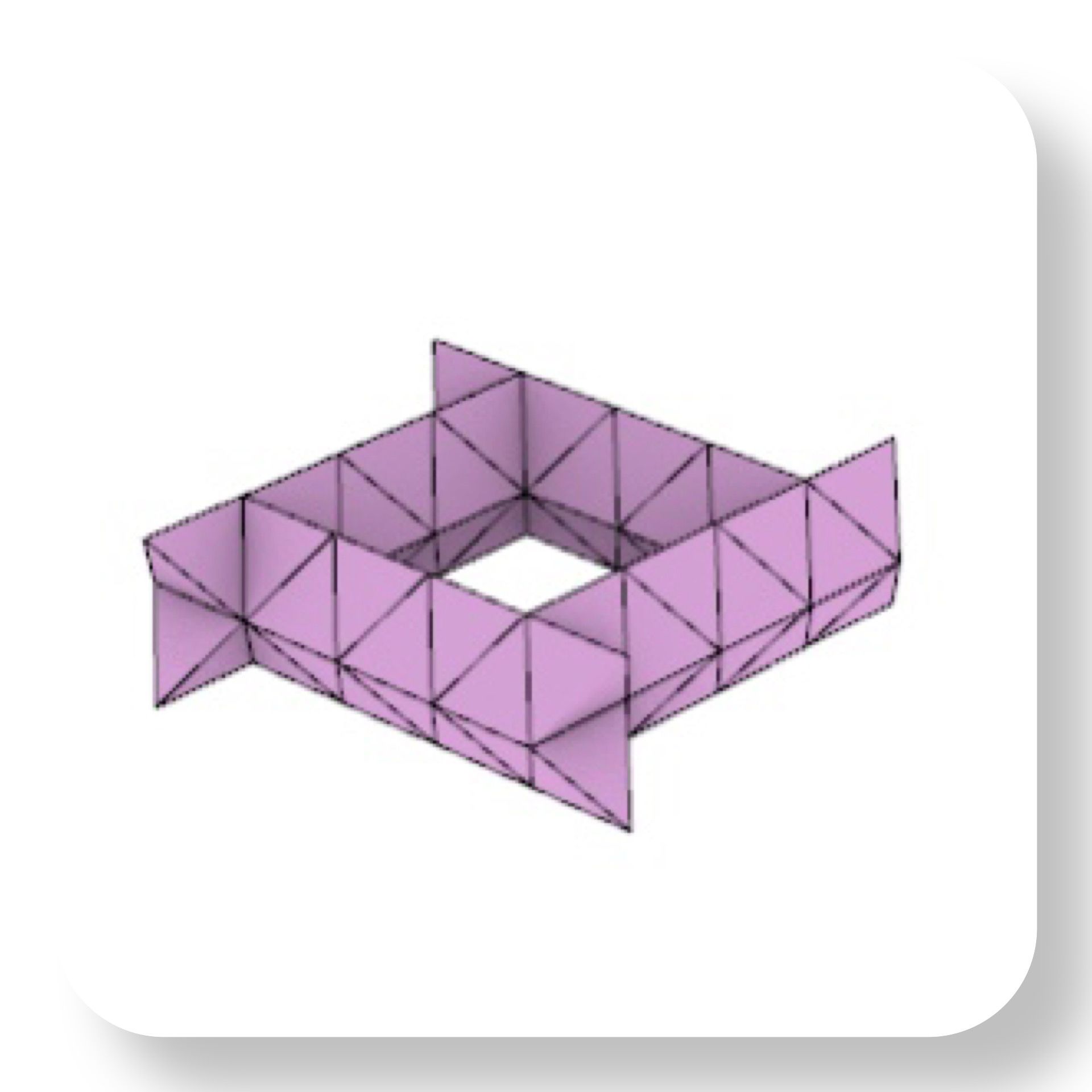}
\end{minipage}
\begin{minipage}{0.5cm}
    
\end{minipage}
\begin{minipage}{.3\textwidth}
    \centering
    \includegraphics[height=4cm]{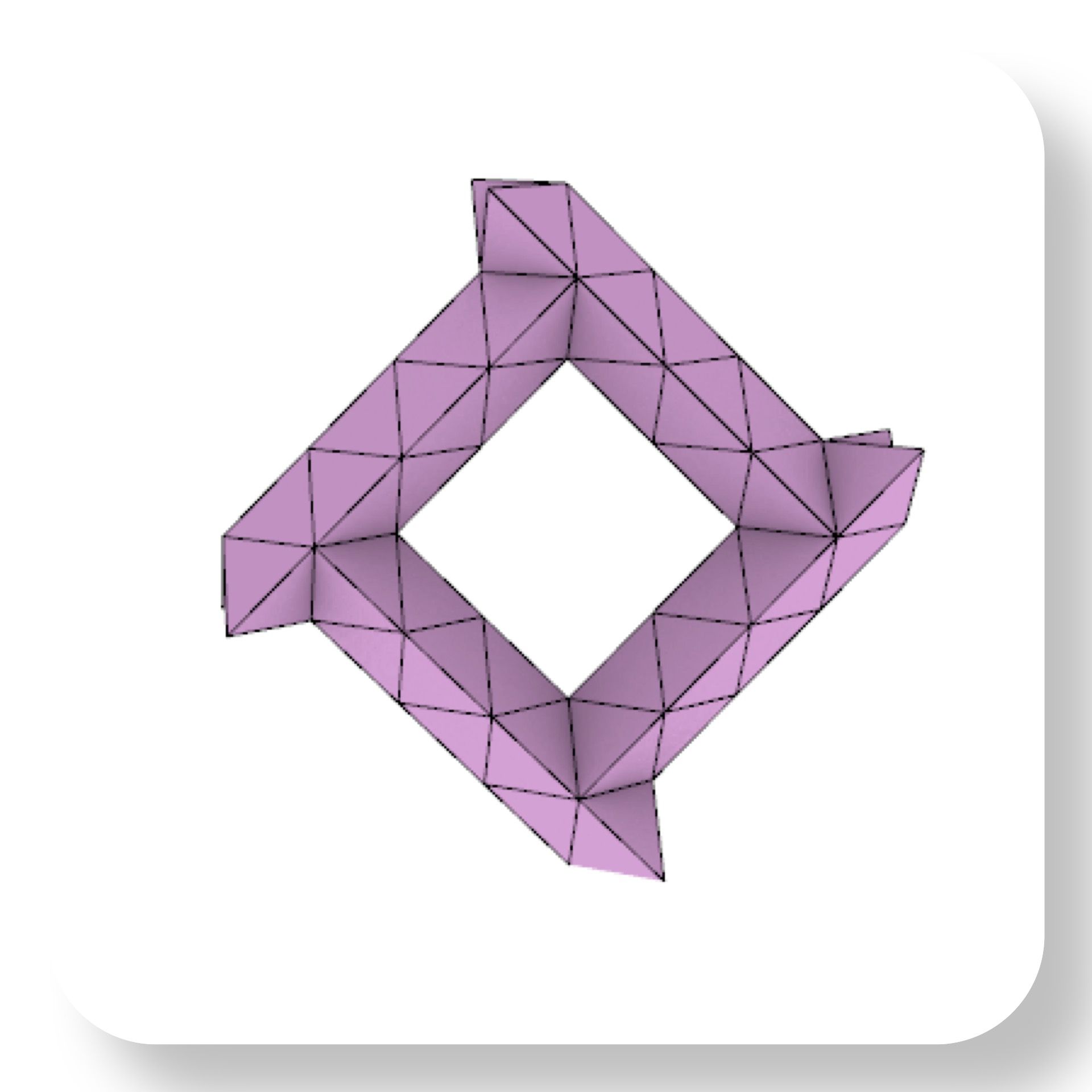}
\end{minipage}
\begin{minipage}{0.5cm}
    
\end{minipage}
\begin{minipage}{.3\textwidth}
    \centering
    \includegraphics[height=4cm]{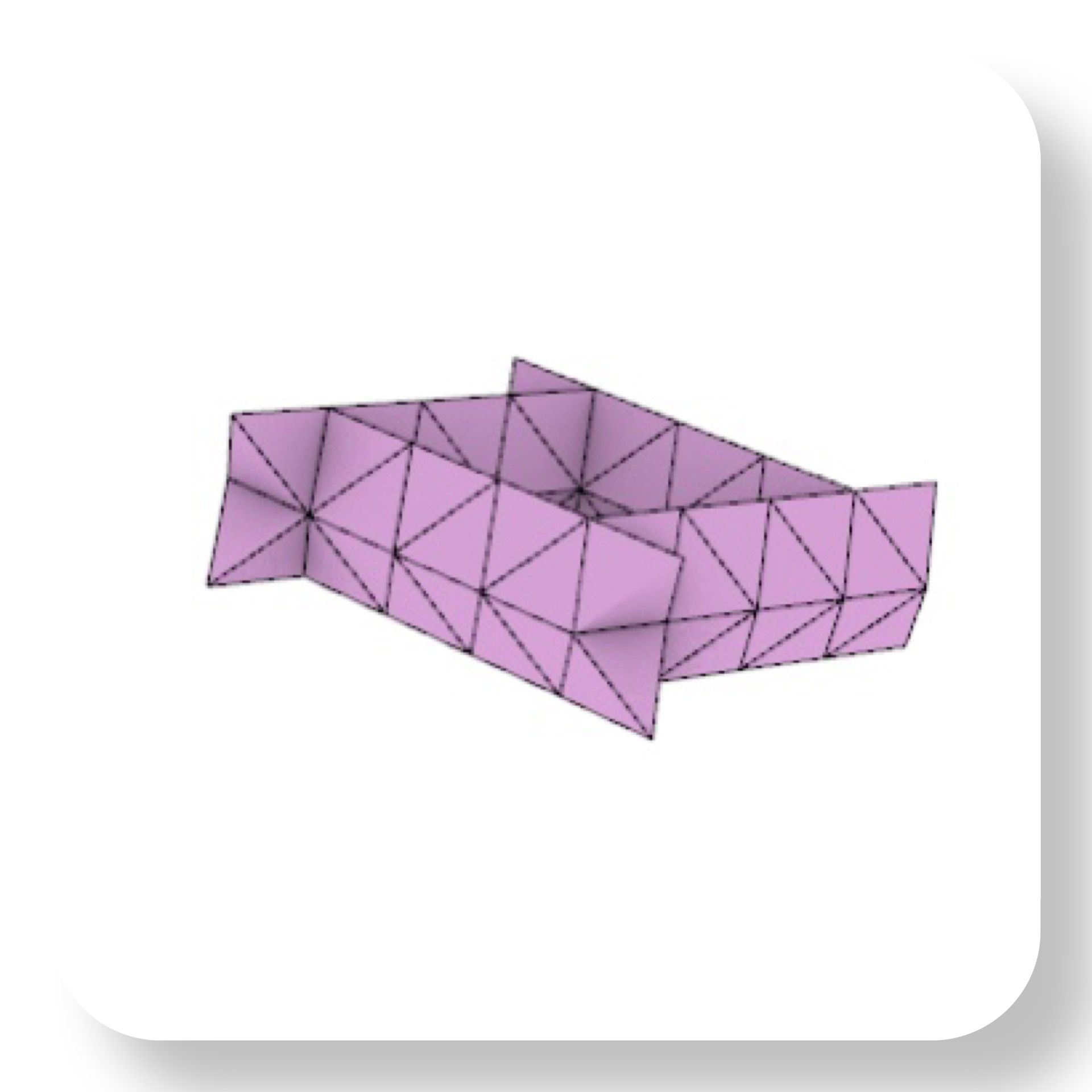}
\end{minipage}
\caption{Various views of the $(3,3)$-shuriken}
\label{33shuriken}
\end{figure}

This observation leads to a more general construction that uses $m$- and $n$-cushions to create an interlocking block. 
\begin{remark}\label{rem2}
    For positive integers $m,n$, the $(m,n)$-shuriken can be defined by the following tetrahedral-octahedral~decomposition:
    \begin{align*}
    (&T_1,\ldots,nv_1+T_1,m(v_2-v_3)+v_1+T_1,\ldots, m(v_2-v_3)+(n+1)v_1+T_1,\\
    &T_2,\ldots,nv_1+T_2,m(v_2-v_3)+v_1+T_2,\ldots, m(v_2-v_3)+(n+1)v_1+T_2,\\
    &T_2,\ldots,v_1+T_2,v_1+(v_2-v_3)+T_2,2v_1+(v_2-v_3)+T_2,\\
    &(v_2-v_3)+T_3,\ldots,m(v_2-v_3)+T_3,v_1+2(v_2-v_3)+T_3,v_1+3(v_2-v_3)+T_3\\
    &(v_2-v_3)+T_4,\ldots,m(v_2-v_3)+T_4,v_1+2(v_2-v_3)+T_4,v_1+3(v_2-v_3)+T_4\\
    &O,\ldots,nv_1+O,m(v_2-v_3)+O,\ldots,nv_1+m(v_2-v_3)+O,\\
    &(v_2-v_3)+O,\ldots,(m-1)(v_2-v_3)+O,\\
    &nv_1+(v_2-v_3)+O,\ldots,nv_1+(m-1)(v_2-v_3)+O).
    \end{align*}
\end{remark}
The tetrahedral-octahedral~decomposition of an $(m,n)$-shuriken consists of $4(n+m+2)$ tetrahedra and $2(n+m)$ octahedra.
Figure~\ref{33shuriken} shows the different views of the $(3,3)$-shuriken. 

Note, the $(m,n)-$ and the $(n,m)-$shuriken are two blocks that can be transformed into each other by using rigid motions, i.e. a 90 degree rotation. For $m,n\geq 2$ the Euler characteristic of the corresponding triangulation of the surface of this block consisting of equilateral triangles is $0$, hence the resulting blocks are tori. It can be shown that the symmetry group of $(m,n)$-shuriken is isomorphic to the dihedral group of order 4 if $m=n$ and isomorphic to $C_2\times C_2$ if $n\neq m$. 

Remarks~\ref{rem1} and \ref{rem2} demonstrate that enforcing the vertices of the tetrahedra and octahedra to be contained in the face-centred cubic lattice, simplifies the task of computing coordinates of the vertices of more complex structures: it is sufficient to apply only translations to vertex positions of the base blocks.

\subsection{Assemblies of the Blocks}
Here, we further examine the blocks that are introduced in Section~\ref{section:NewBlocks} by presenting corresponding topological interlocking assemblies. We only focus on assemblies whose blocks can be placed between the following two planes:
\[
P_1=\langle v_1,v_2-v_3 \rangle,P_2 = 2v_2+\langle v_1,v_2-v_3 \rangle.
\]
In the following, we give figures that illustrate the different interlocking assemblies. For simplicity, we colour the blocks of the corresponding assembly in three different colours. Here, we reserve one colour (black or grey) to show the frame of the assembly, whose definition is necessary to establish the interlocking property.

In order to provide a better understanding of the assemblies, we translate the assemblies into combinatorial graphs to obtain a simplified description of the assemblies. For this purpose, we need to specify the nodes and edges of the graph that is yet to be defined. This can be done in the following way:
\begin{definition}
Let $(X_i)_{i\in I}$ be an assembly of the presented blocks and $I$ a countable index set. The \emph{assembly graph}  $G=(V,E)$ of the assembly is given by the vertices $V=I$ and the edges $E\subseteq I \times I$ satisfying the following condition:

 The set $\{i,j\}$ is an edge in $G$, (that is $(\{i,j\}\in E)$) if and only if the blocks $X_i$ and $X_j$ share at least one common face in the given assembly.
\end{definition}

Translating a given assembly into an assembly graph is straight forward whereas we require additional information to reconstruct a assembly from the combinatorial structure of the assembly graph.
For instance, this reconstruction can also be achieved by denoting all possibilities of assembling any two blocks of the given assembly  to colour the edges of the assembly graph.

For simplicity, we will refrain from formal definitions of the assembly graphs of the corresponding topological interlocking assemblies and use figures to illustrate their combinatorial structures. 
\subsubsection{Assembly of the Kitten}
\begin{figure}[H]
\begin{minipage}{.3\textwidth}
    \centering
    \includegraphics[height=4cm]{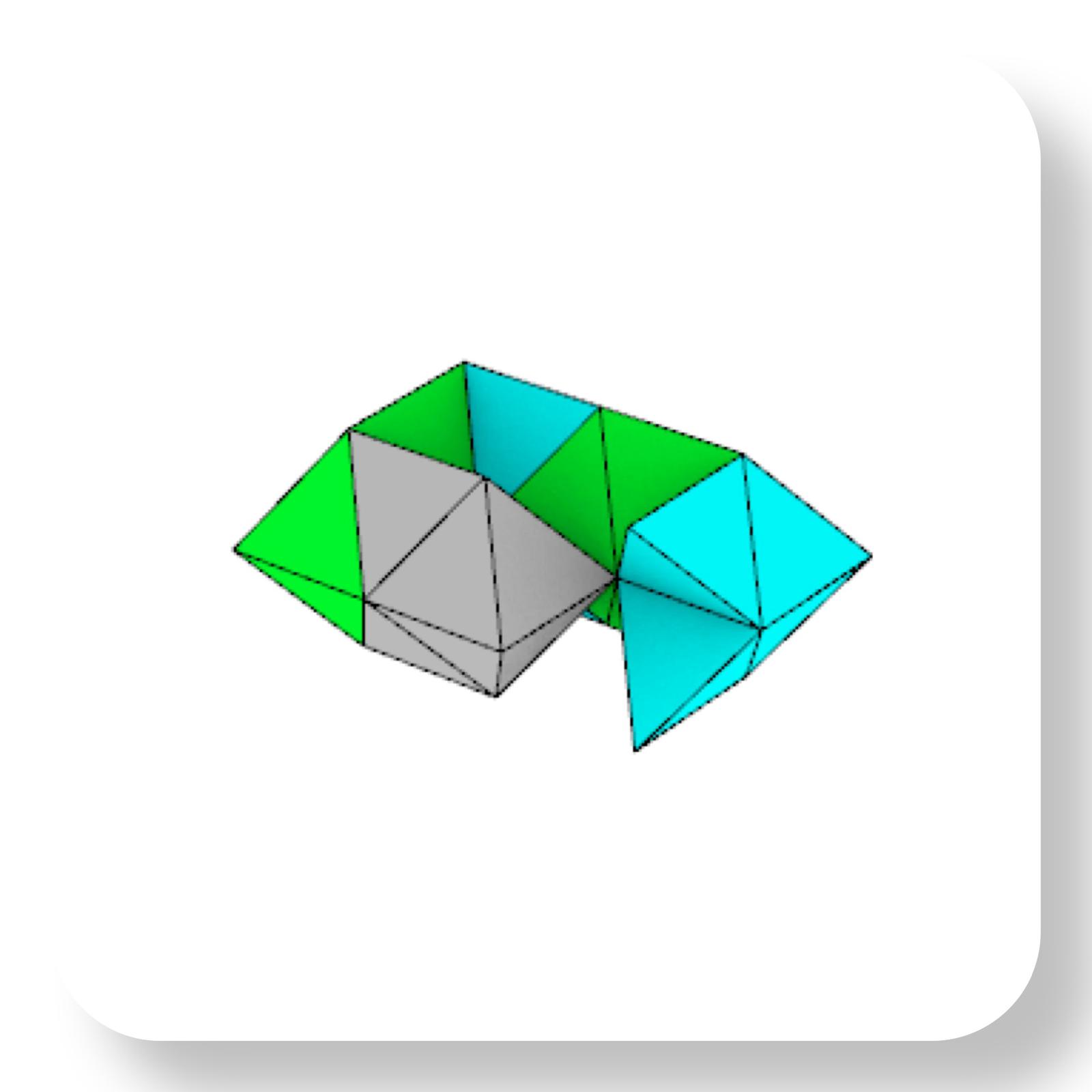}
\end{minipage}
\begin{minipage}{0.5cm}
    
\end{minipage}
\begin{minipage}{.3\textwidth}
    \centering
    \includegraphics[height=4cm]{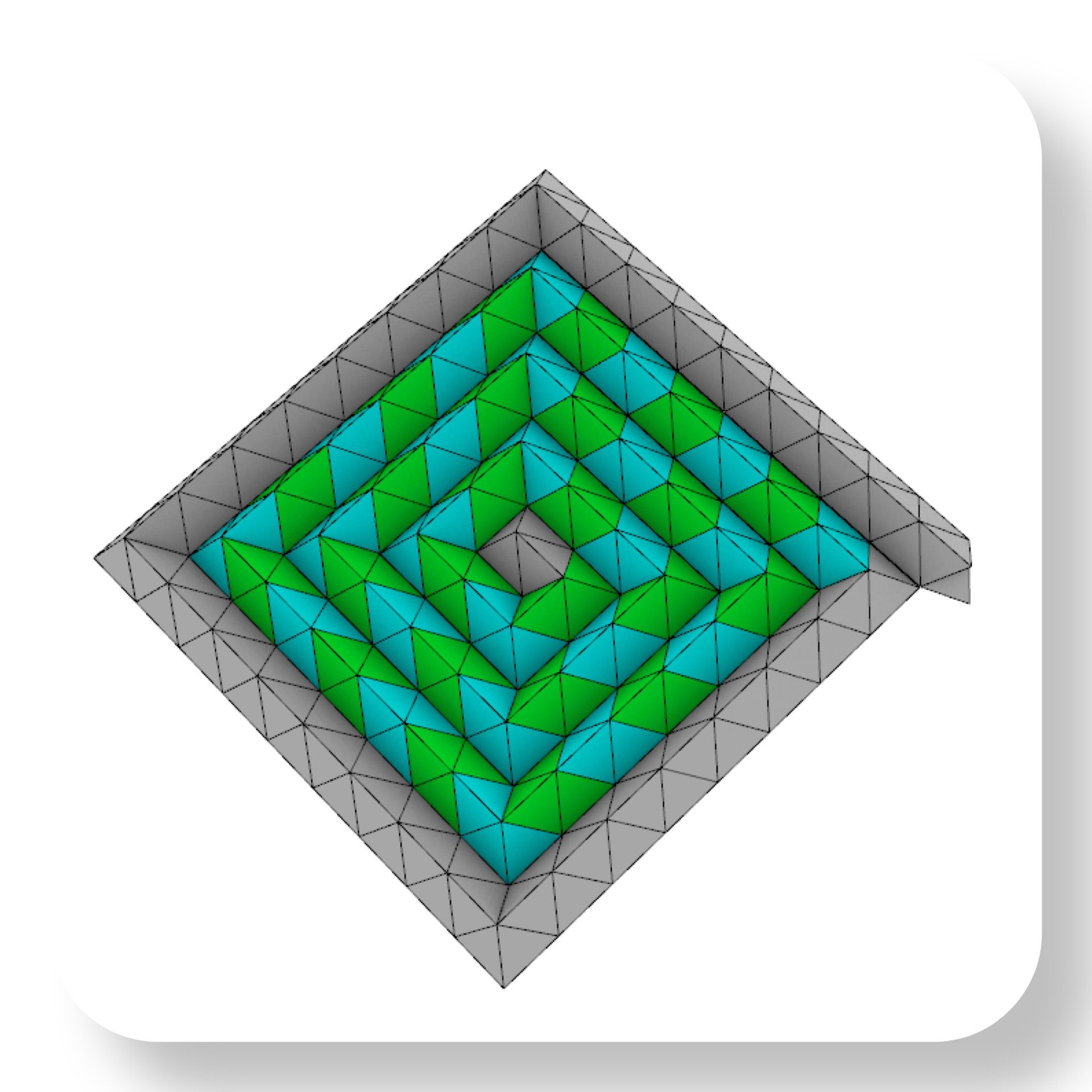}
\end{minipage}
\begin{minipage}{0.5cm}
    
\end{minipage}
\begin{minipage}{.3\textwidth}
    \centering
    \includegraphics[height=4cm]{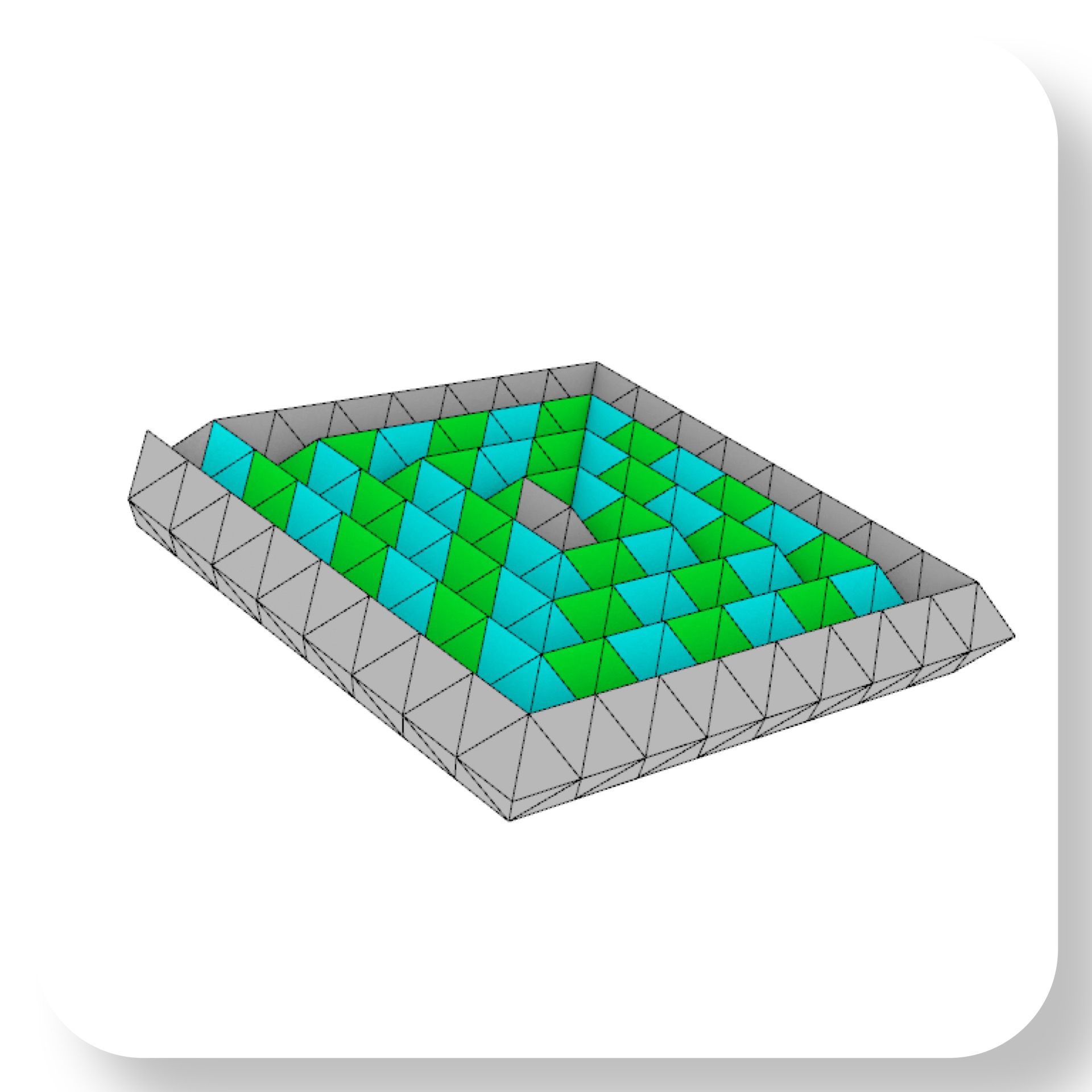}
\end{minipage}
\caption{ Assembly of the kitten in a spiral form}
\label{assemblykitten}
\end{figure}

Since the kitten is space-filling it allows various assemblies. The assembly that gives rise to a topological interlocking assembly is shown in Figure~\ref{assemblykitten}.

Here, copies of the kitten are assembled to construct a spiral form. The frame of the corresponding assembly is given by the blocks coloured in grey. 

\begin{figure}[H]
    \centering
    \includegraphics[scale=0.8,viewport=10cm 20.5cm 0cm 26.5cm]{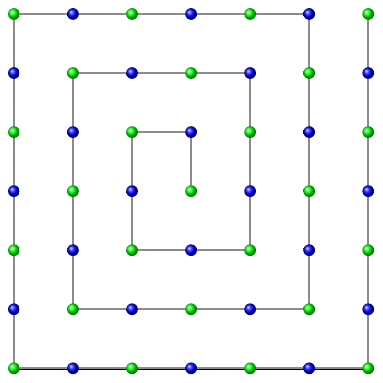}
\caption{Assembly graph of the given assembly of the kitten}
\label{assemblygraphkitten}
\end{figure}

 This assembly can be visualised by the following assembly graph, see Figure~\ref{assemblygraphkitten}. Note, the assembly graph of the kitten is isomorphic to a spanning tree.
 
\begin{figure}[H]
\begin{minipage}{.3\textwidth}
    \centering
    \includegraphics[height=4cm]{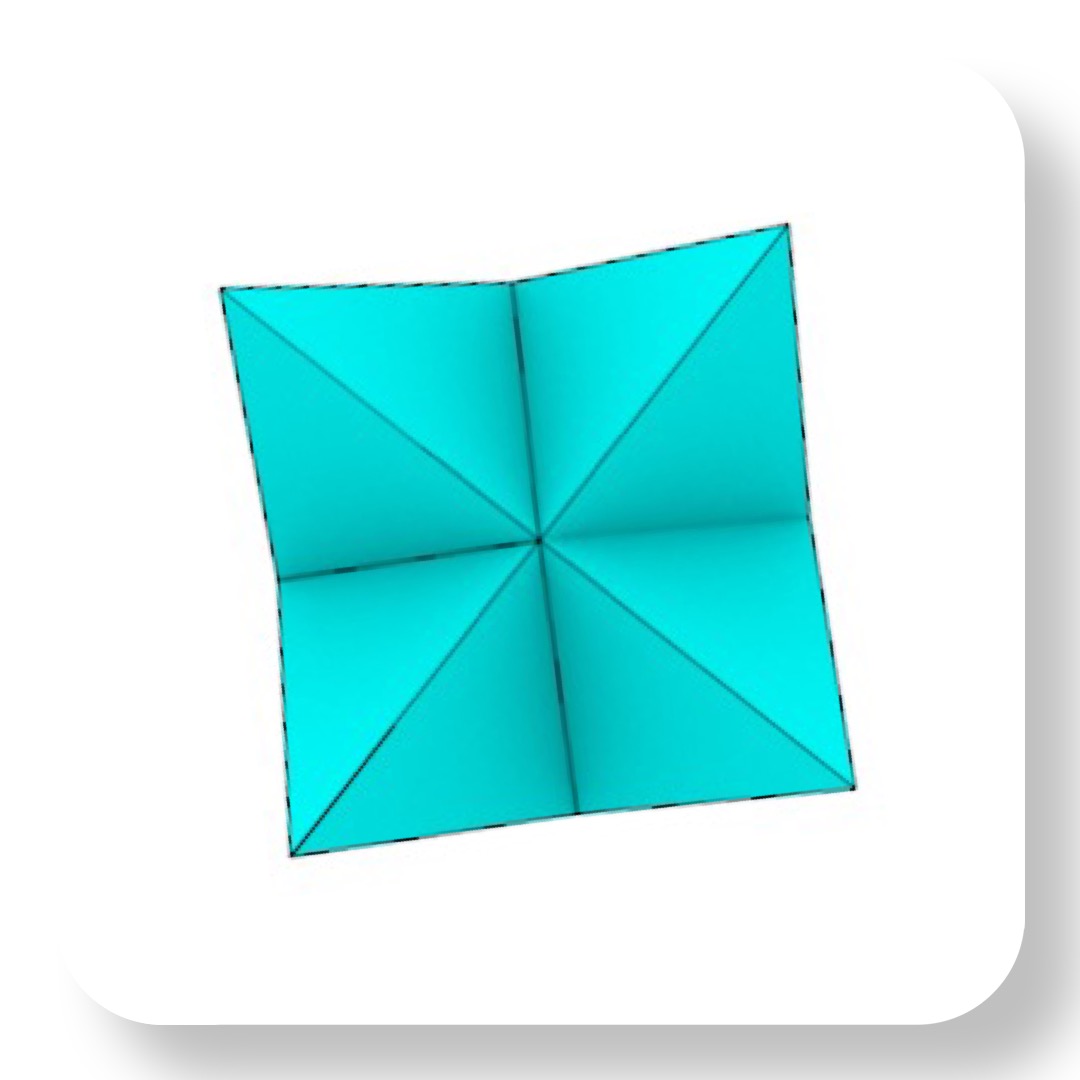}
\end{minipage}
\begin{minipage}{0.5cm}
    
\end{minipage}
\begin{minipage}{.3\textwidth}
    \centering
    \includegraphics[height=4cm]{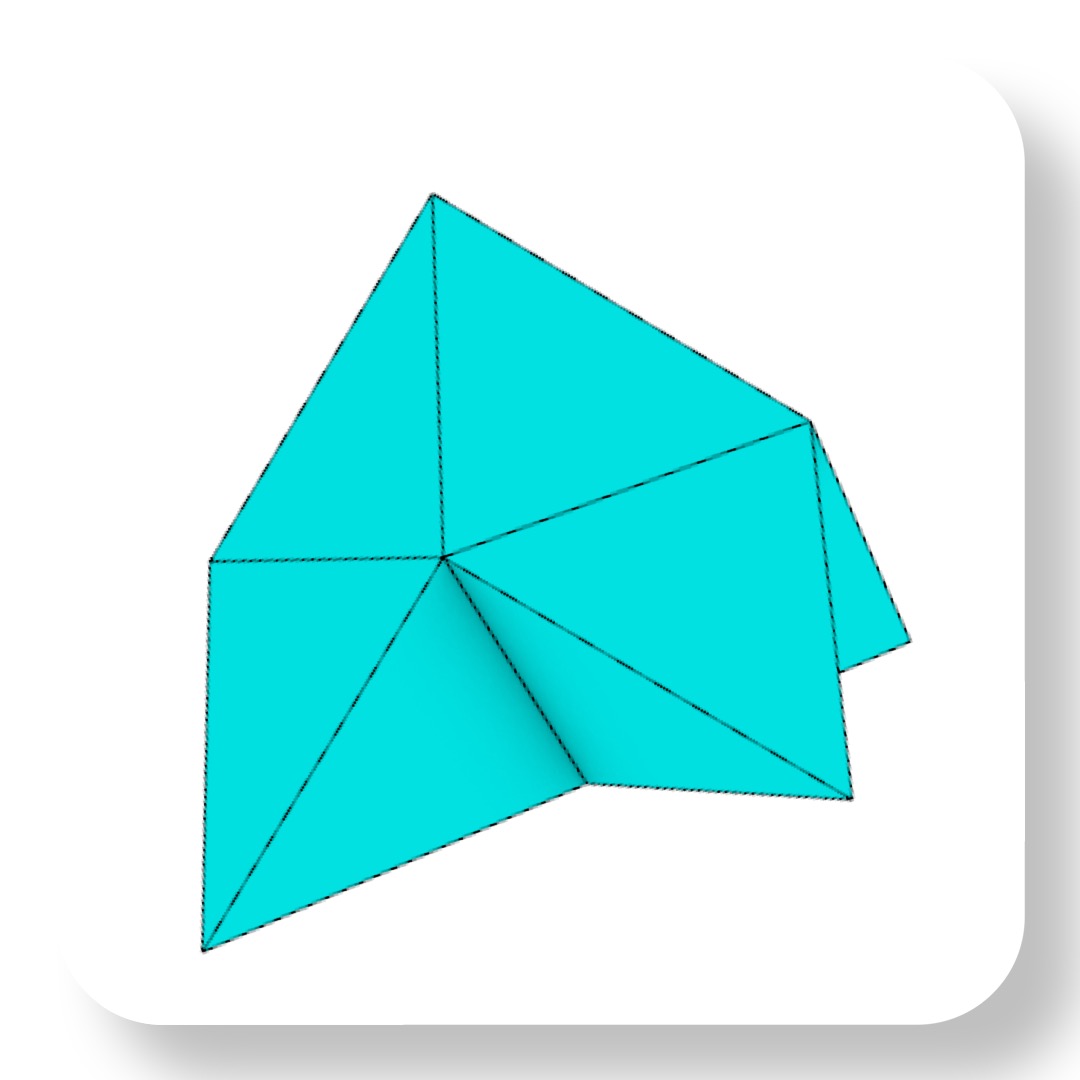}

\end{minipage}
\begin{minipage}{0.5cm}
    
\end{minipage}
\begin{minipage}{.3\textwidth}
    \centering
    \includegraphics[height=4cm]{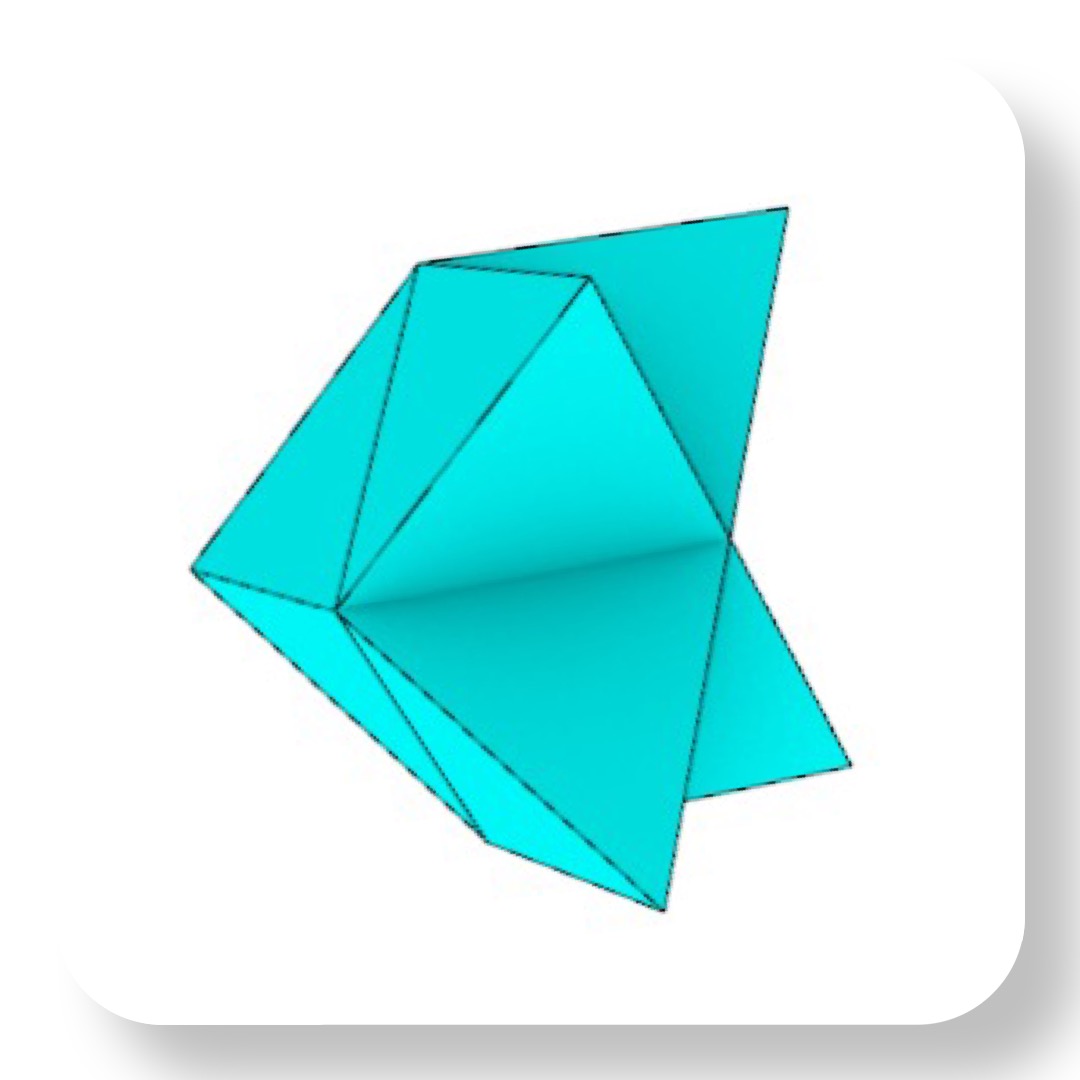}
\end{minipage}
\caption{Adding two tetrahedra to the kitten which leads to the ufo consisting of one octahedron and four tetrahedra}
\label{ufo}
\end{figure}

\begin{figure}[H]
\begin{minipage}{.3\textwidth}
    \centering
    \includegraphics[height=4cm]{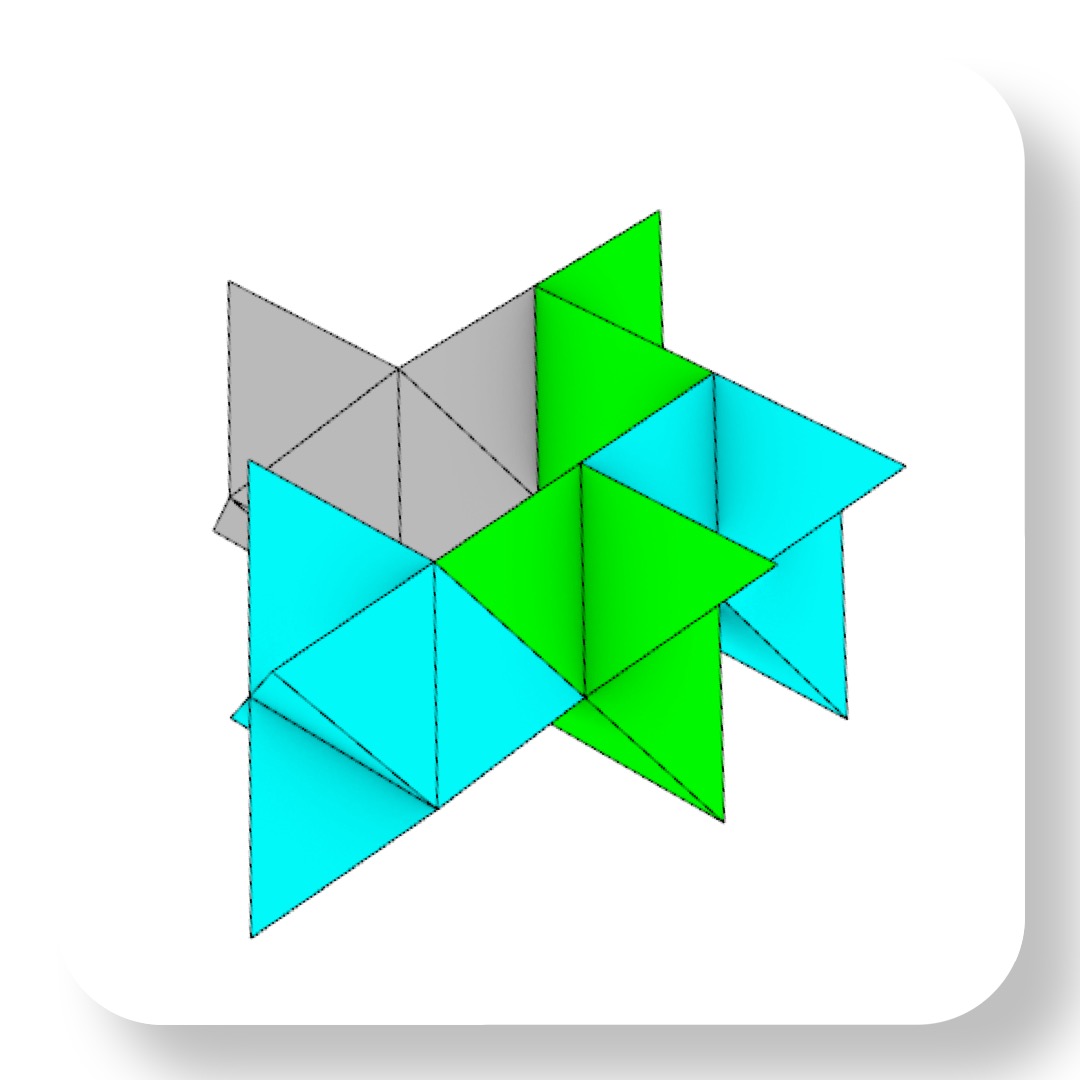}
\end{minipage}
\begin{minipage}{0.5cm}
    
\end{minipage}
\begin{minipage}{.3\textwidth}
    \centering
    \includegraphics[height=4cm]{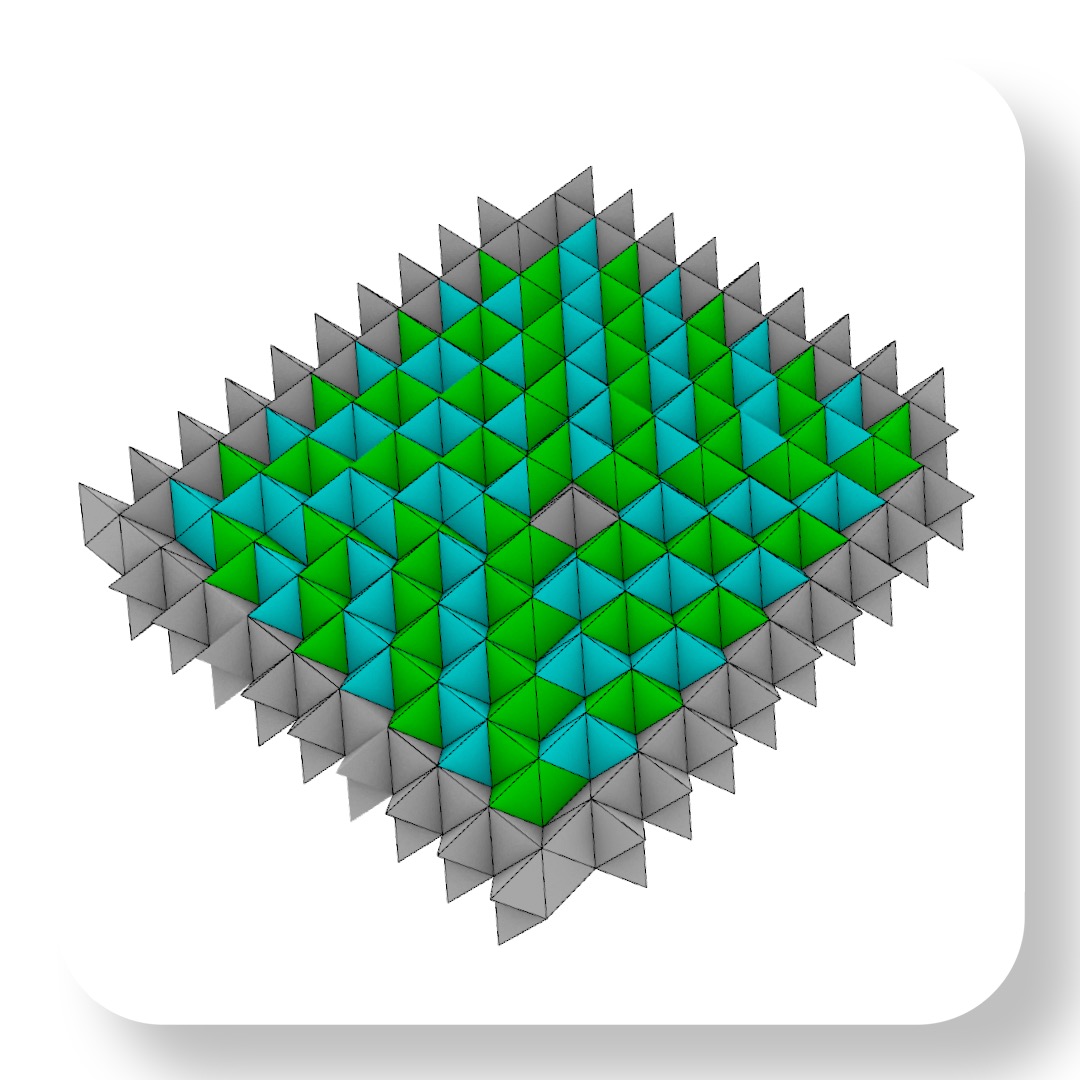}

\end{minipage}
\begin{minipage}{0.5cm}
    
\end{minipage}
\begin{minipage}{.3\textwidth}
    \centering
    \includegraphics[height=4cm]{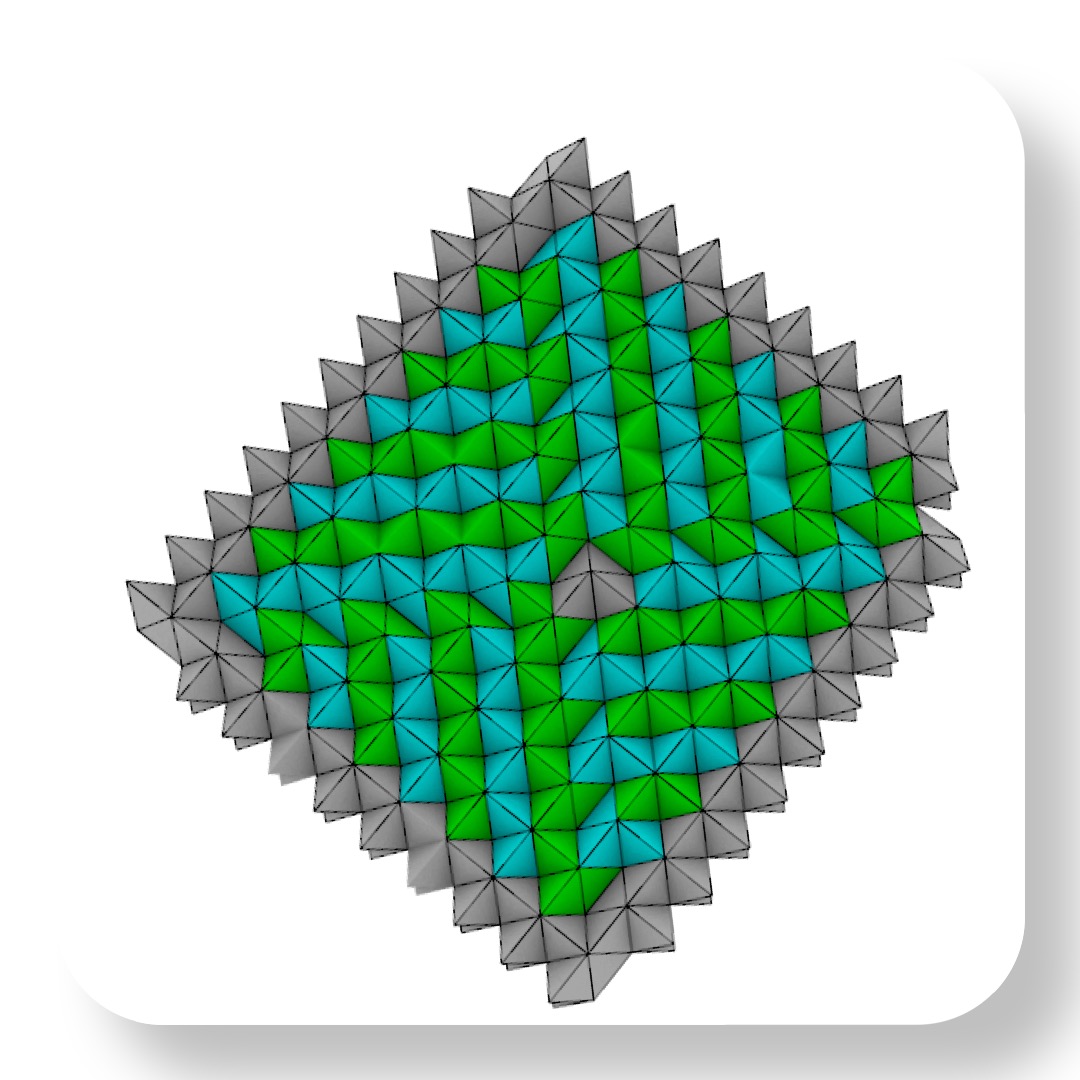}
\end{minipage}
\caption{Assembly of the newly constructed block in a spiral form}
\label{ufoassembly}
\end{figure}

Moreover, by conducting real-life experiments with 3D-printed versions of the kitten, we see that the stability of the topological interlocking assembly, illustrated in Figure \ref{assemblykitten} can be further improved. This improvement is achieved by applying the following modifications to the kitten: in our experiments we have observed that adding two additional tetrahedra to the kitten results in a block that follows the same assembly rule as the kitten and exhibits greater stability compared to the topological interlocking assembly illustrated in Figure \ref{assemblykitten}. More precisely, this modified block can be described by the tetrahedral-octahedral decomposition:
\[
(O,T_1,T_2,T_3,T_4).
\]
The block that results from this tetrahedral-octahedral decomposition is shown in Figure \ref{ufo}. We call this constructed block the \emph{ufo.}
Figure \ref{ufoassembly} illustrates the spiral form that is constructed by assembling copies of this modified block.
\subsubsection{Assembly of the cushion}
\begin{figure}[H]
\begin{minipage}{.3\textwidth}
    \centering
    \includegraphics[height=4cm]{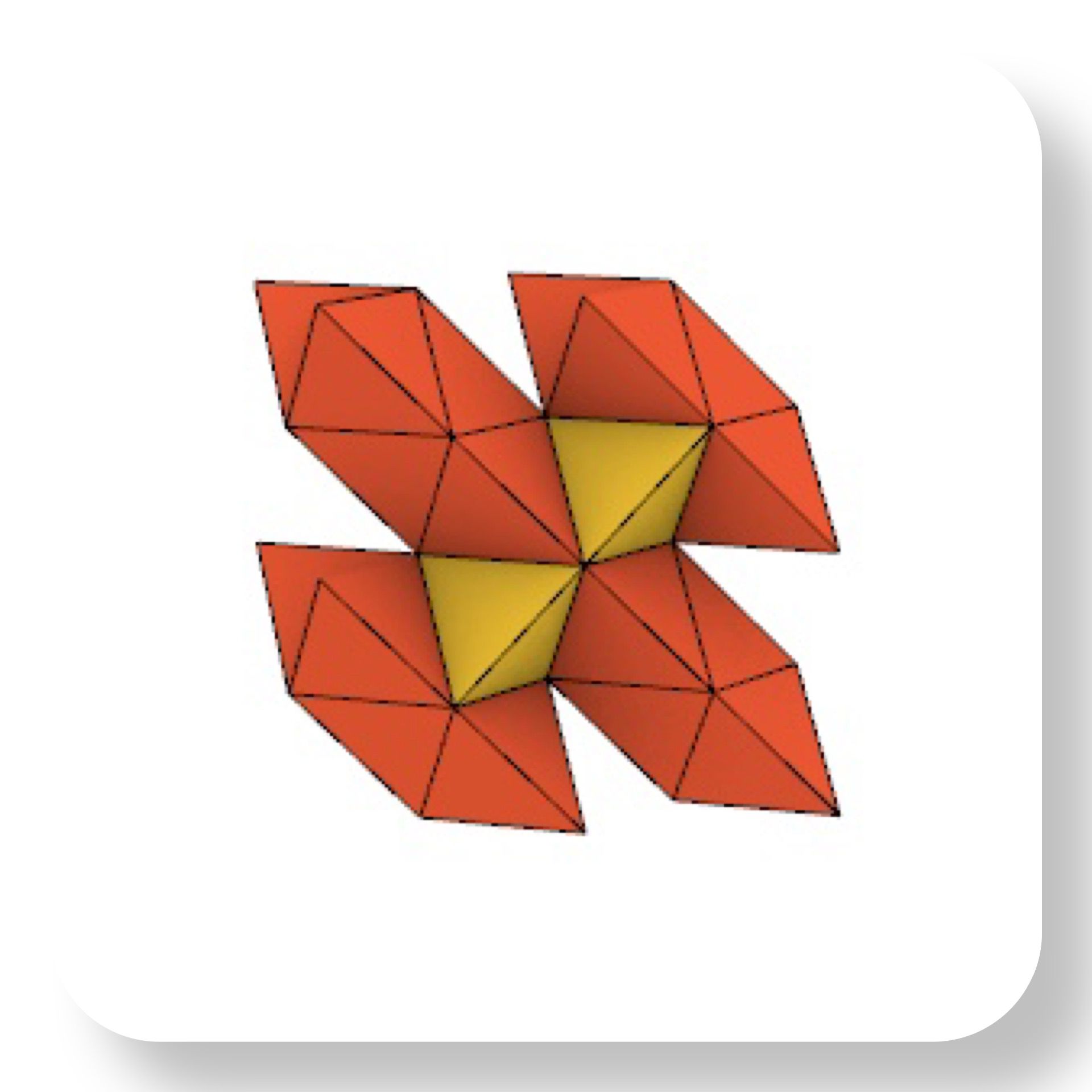}
\end{minipage}
\begin{minipage}{0.5cm}
    
\end{minipage}
\begin{minipage}{.3\textwidth}
    \centering
    \includegraphics[height=4cm]{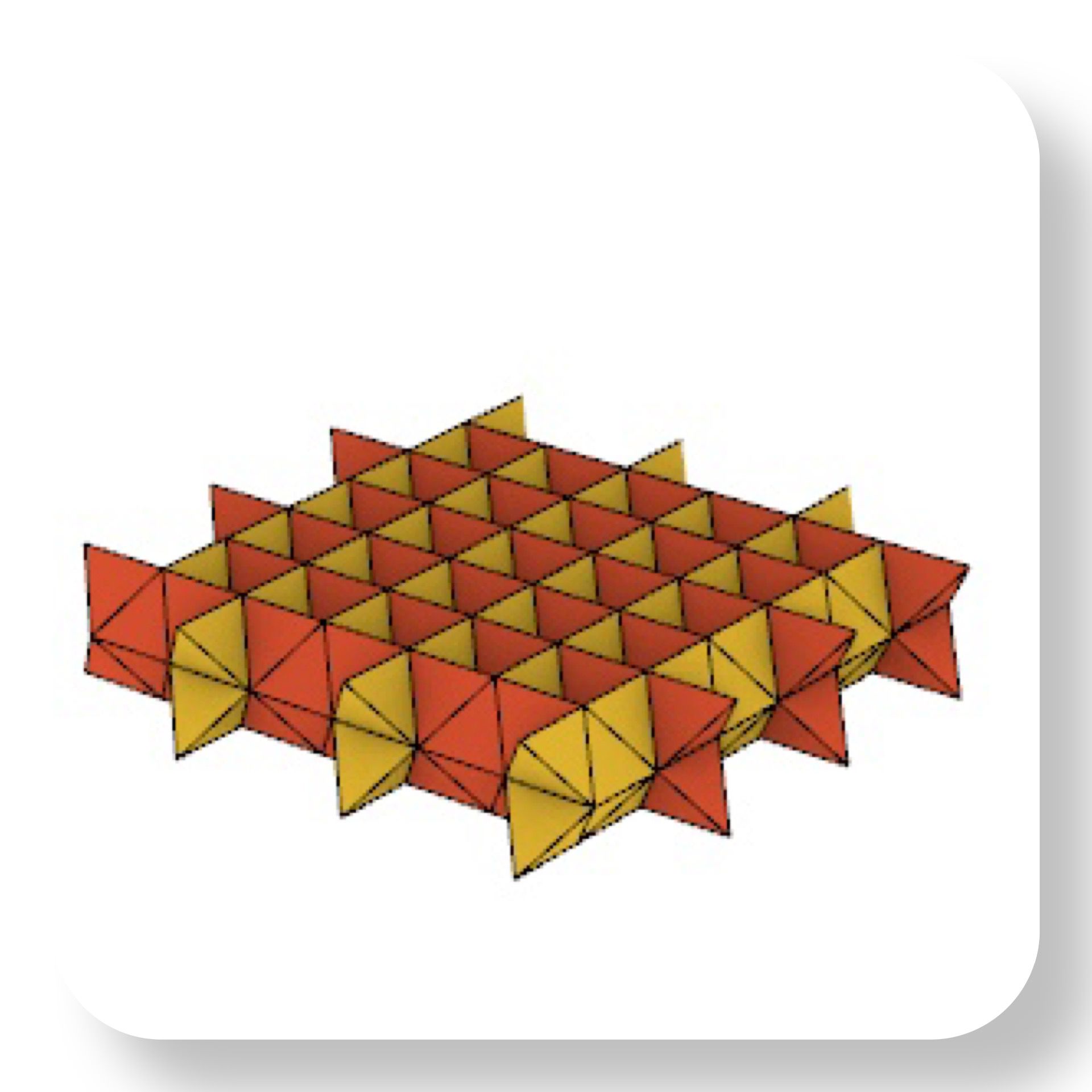}
\end{minipage}
\begin{minipage}{0.5cm}
    
\end{minipage}
\begin{minipage}{.3\textwidth}
    \centering
    \includegraphics[height=4cm]{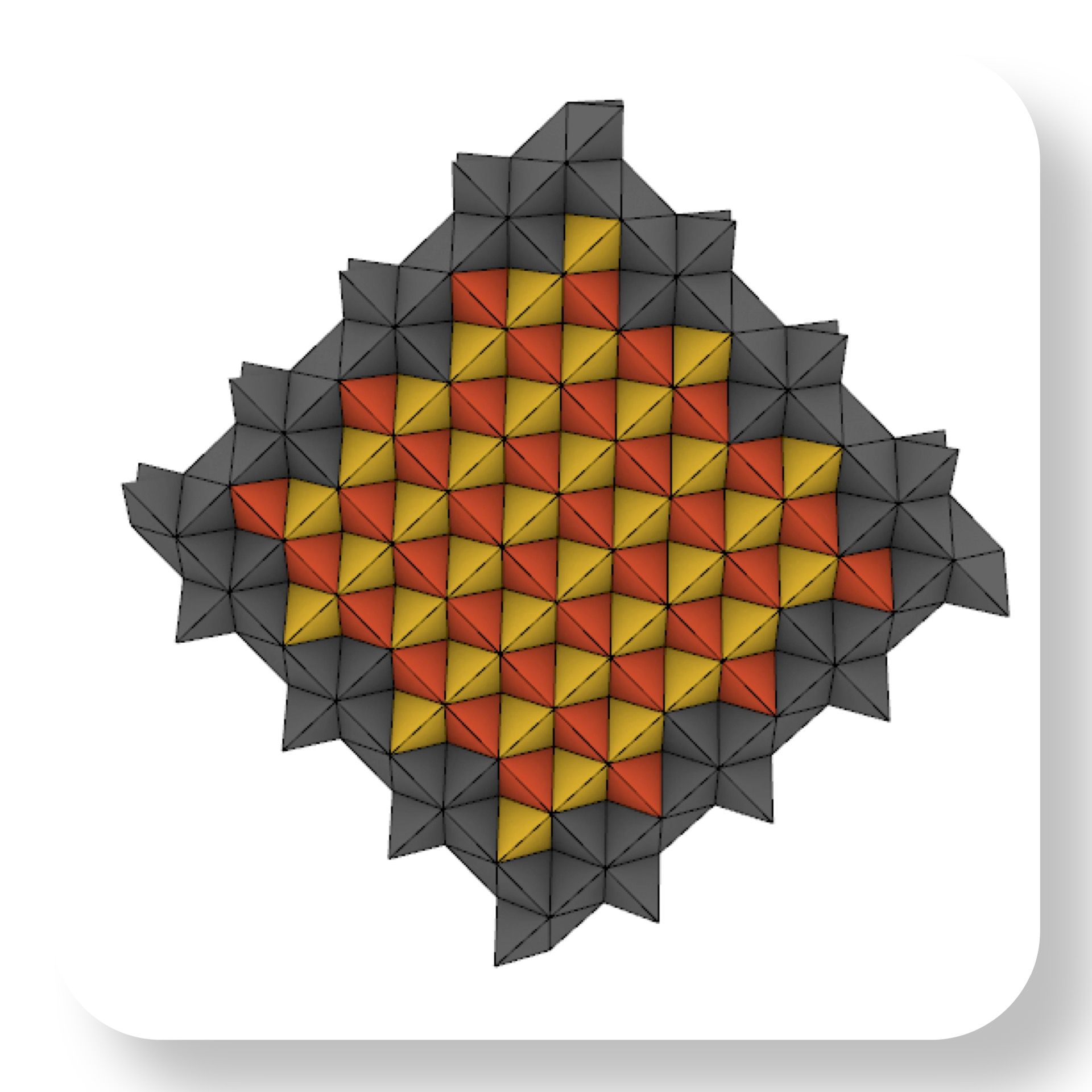}
\end{minipage}
\caption{Assembly of the cushion}
\label{assemblycushion}
\end{figure}

Next, we give topological interlocking assemblies of the different $n$-cushion. 
For $n=1$, we can realise a topological interlocking by assembling copies of the cushion in a quadratic grid so that the resulting assembly is contained in the tetroctahedrille, see Figure~\ref{assemblycushion}. Here, the frame is given by the blocks coloured in grey.

\begin{figure}[H]
\begin{minipage}{.3\textwidth}
    \centering
    \includegraphics[height=4cm]{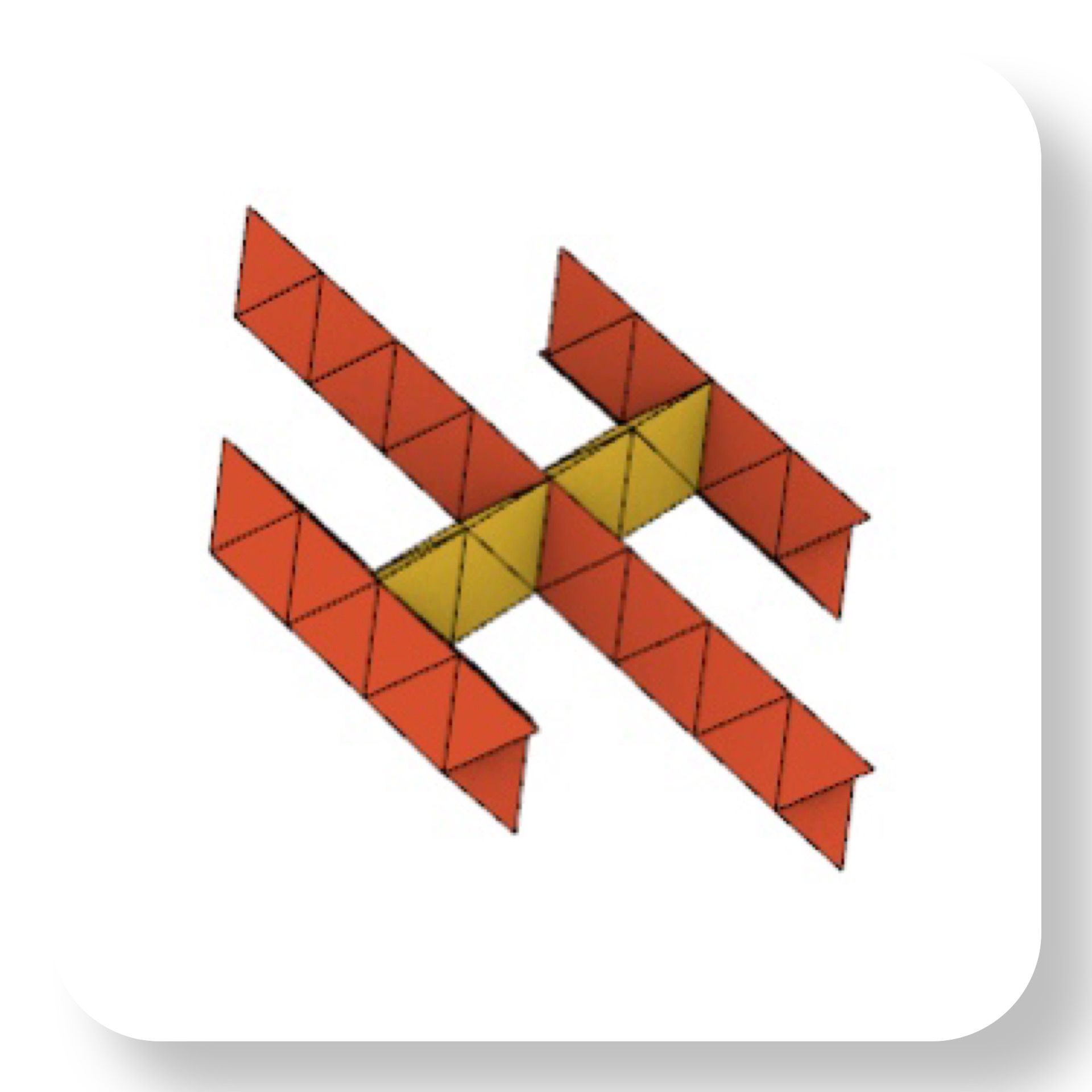}
\end{minipage}
\begin{minipage}{0.5cm}
    
\end{minipage}
\begin{minipage}{.3\textwidth}
    \centering
    \includegraphics[height=4cm]{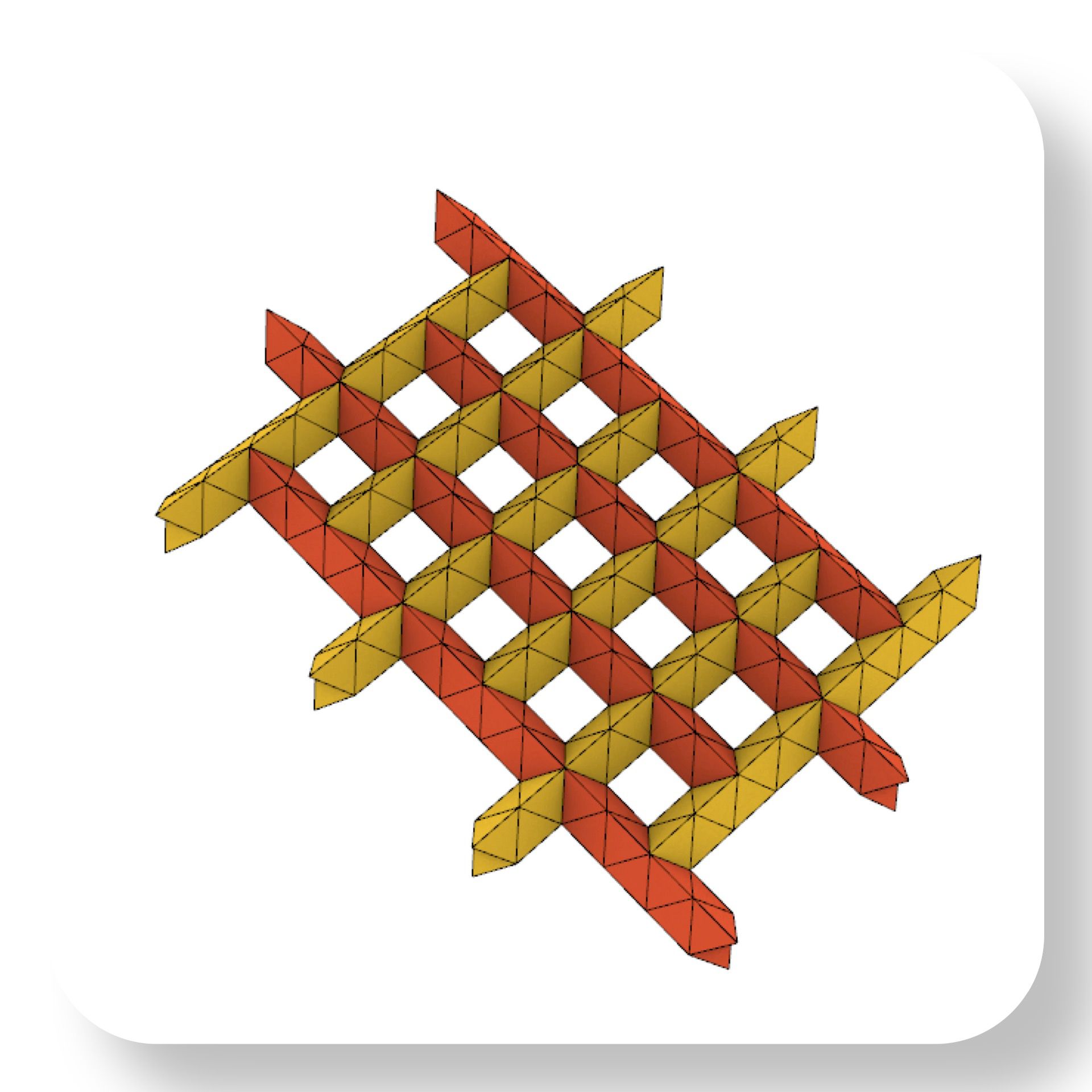}
\end{minipage}
\begin{minipage}{0.5cm}
    
\end{minipage}
\begin{minipage}{.3\textwidth}
    \centering
    \includegraphics[height=4cm]{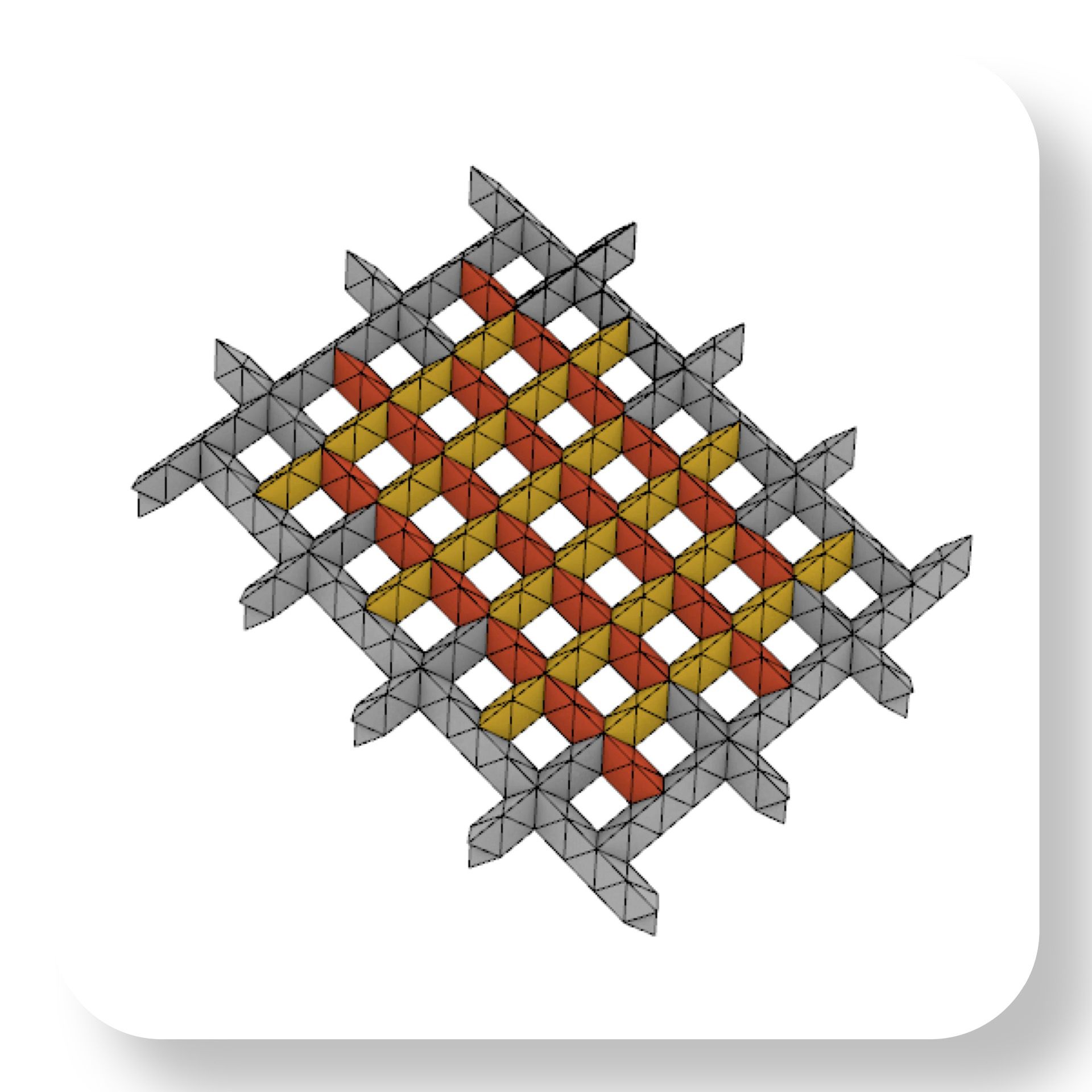}
\end{minipage}
\caption{Assembly of the 3-cushion}
\label{assembly3cushion}
\end{figure}

The surface of this assembly looks like the surface of the tetrahedra-interlocking illustrated in Figure \ref{tetra_oct_interlocking}.
We can apply the same logic for the general case and assemble $n$-cushions in a grid satisfying the rules of the tetroctahedrille to construct topological interlocking assemblies, see Figure~\ref{assembly3cushion}. In this figure the frame of the assembly is coloured in grey.

\begin{figure}[H]
    \centering    
    \includegraphics[scale=1,viewport=10cm 21.5cm 0cm 26.5cm]{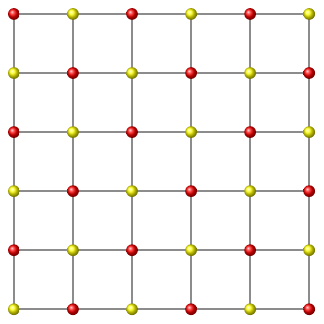}
\caption{Assembly graph of the given assembly of the kitten}
\label{assemblygraphcushion}
\end{figure}

Constructing a topological interlocking assembly by following the introduced logic of forming a grid with copies of the $n$-cushion for arbitrary $n$ gives rise to a graph that can be embedded into the square lattice, see Figure~\ref{assemblygraphcushion}. Hence, these assemblies give rise to isomorphic graphs.

\subsubsection{Assembly of the Shuriken}
\begin{figure}[H]
\begin{minipage}{.3\textwidth}
    \centering
    \includegraphics[height=4cm]{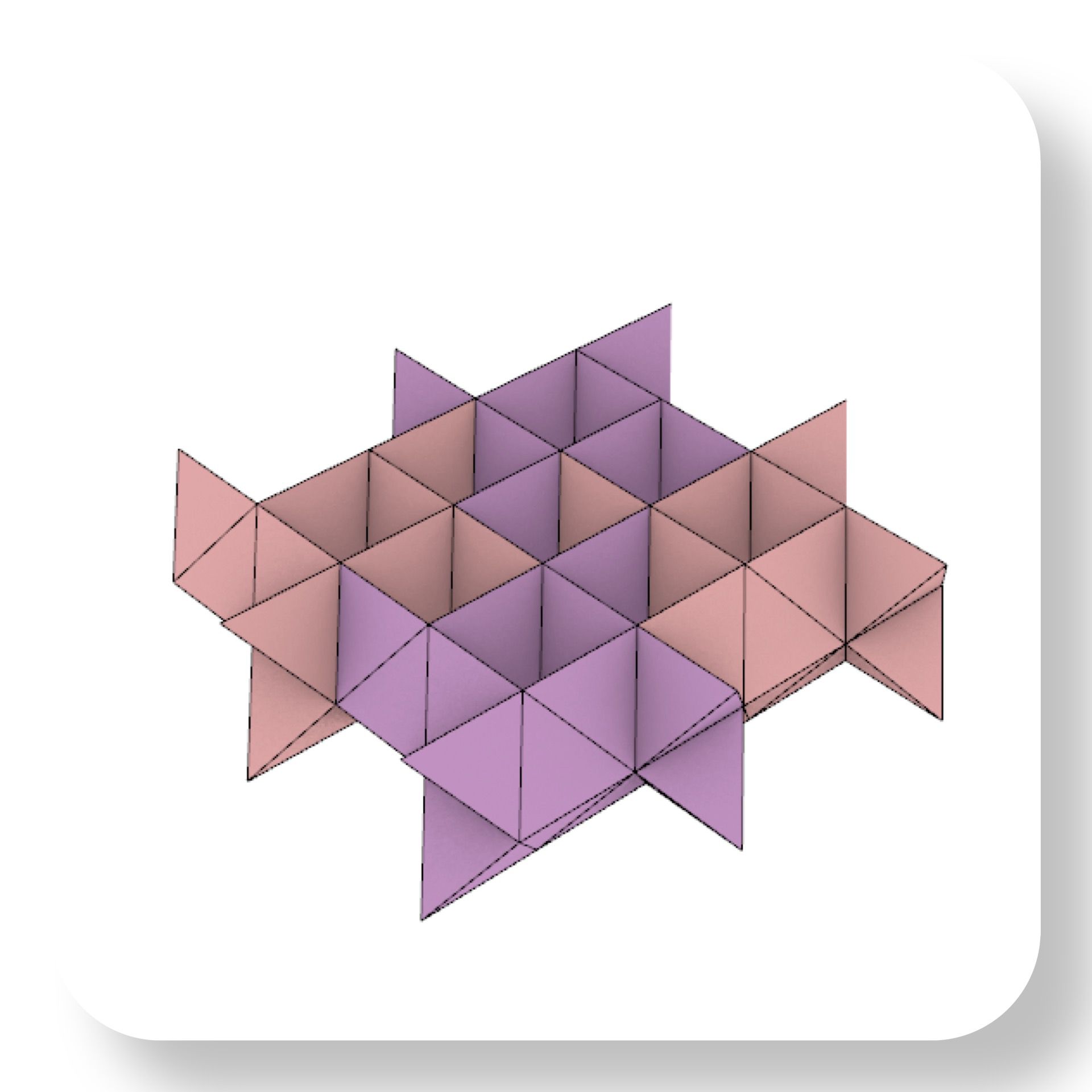}
\end{minipage}
\begin{minipage}{0.5cm}
    
\end{minipage}
\begin{minipage}{.3\textwidth}
    \centering
    \includegraphics[height=4cm]{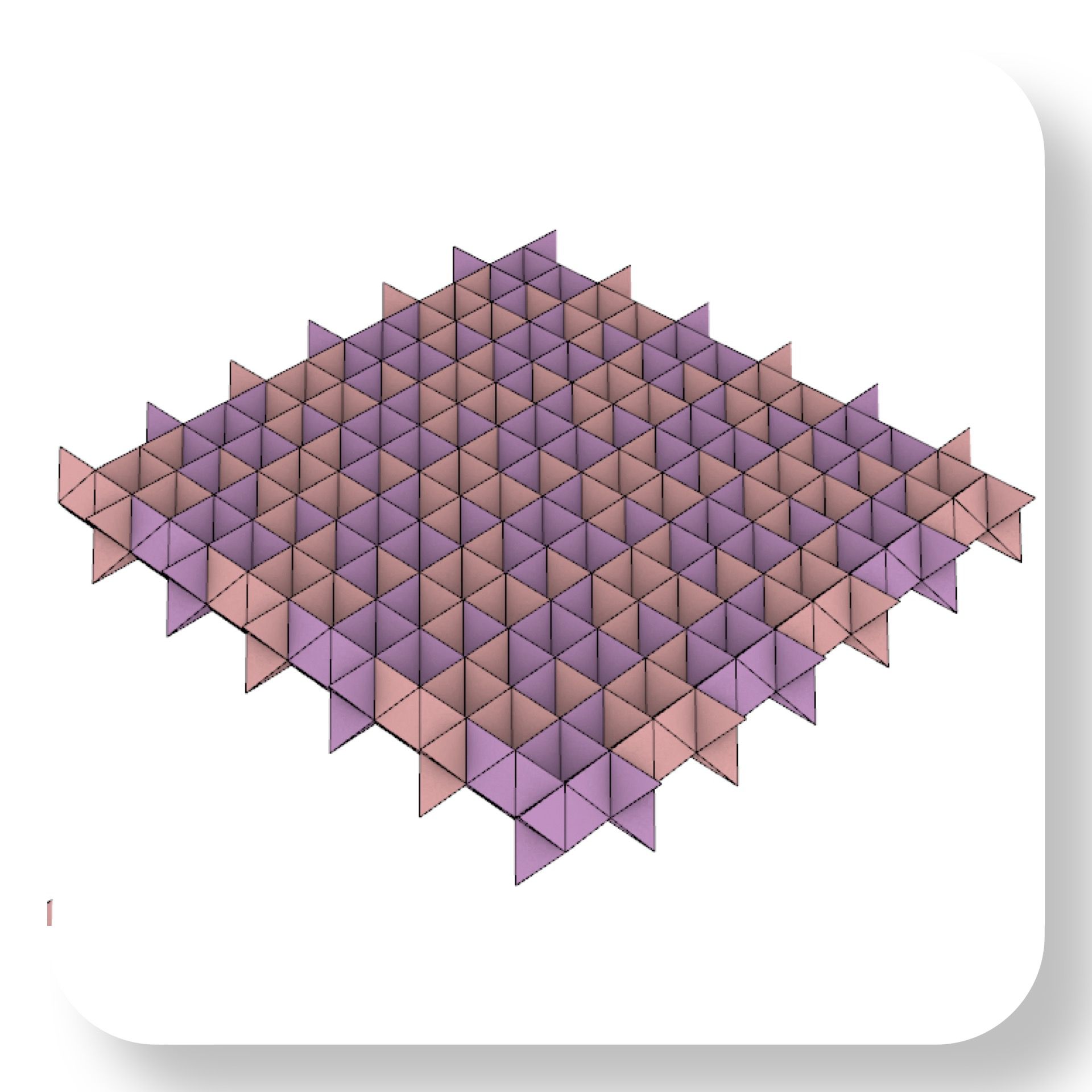}
\end{minipage}
\begin{minipage}{0.5cm}
    
\end{minipage}
\begin{minipage}{.3\textwidth}
    \centering
    \includegraphics[height=4cm]{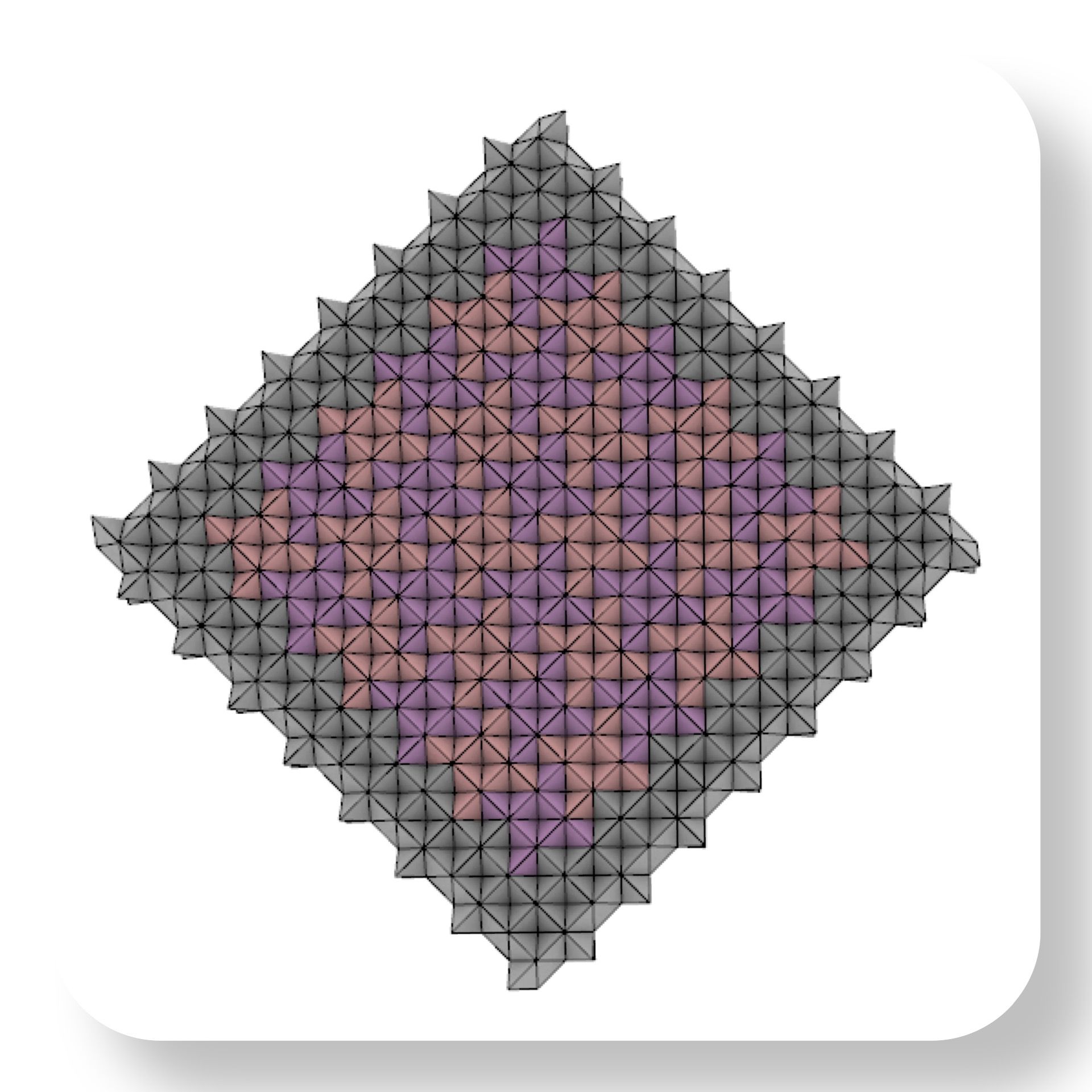}
\end{minipage}
\caption{Assembly of the shuriken}
\label{assemblyShuriken}
\end{figure}

\begin{figure}[H]
\begin{minipage}{.3\textwidth}
    \centering
    \includegraphics[height=4cm]{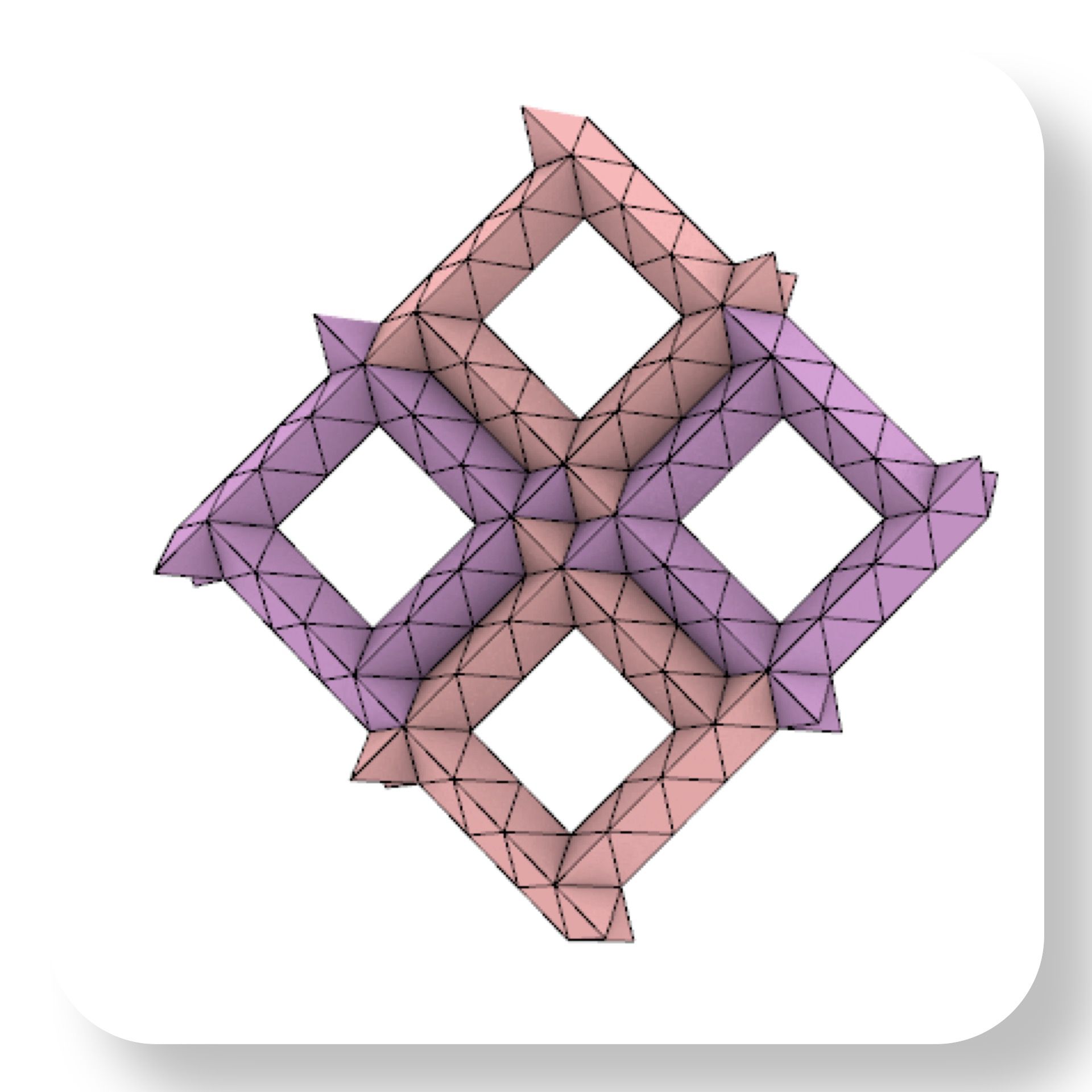}
\end{minipage}
\begin{minipage}{0.5cm}
    
\end{minipage}
\begin{minipage}{.3\textwidth}
    \centering
    \includegraphics[height=4cm]{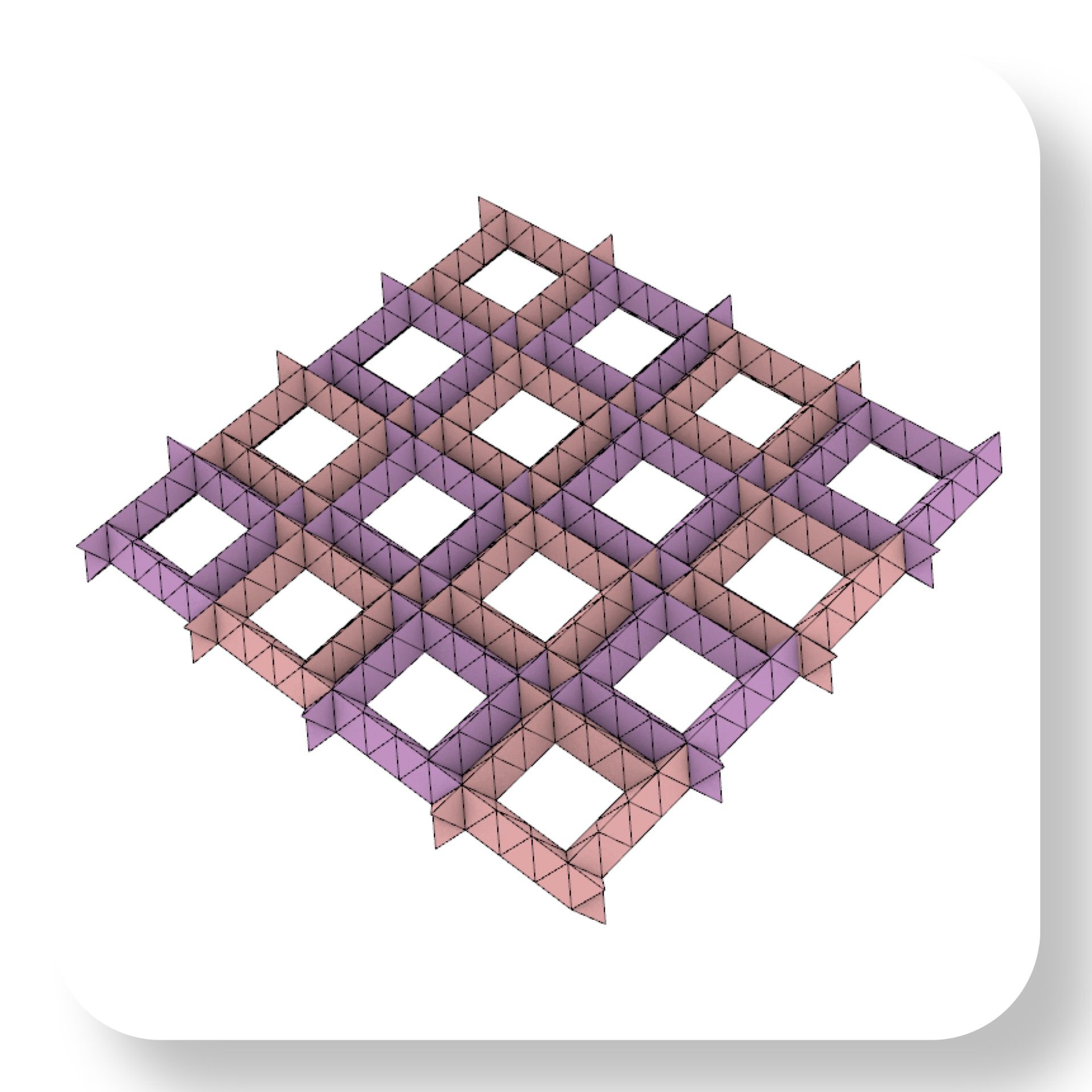}
\end{minipage}
\begin{minipage}{0.5cm}
    
\end{minipage}
\begin{minipage}{.3\textwidth}
    \centering
    \includegraphics[height=4cm]{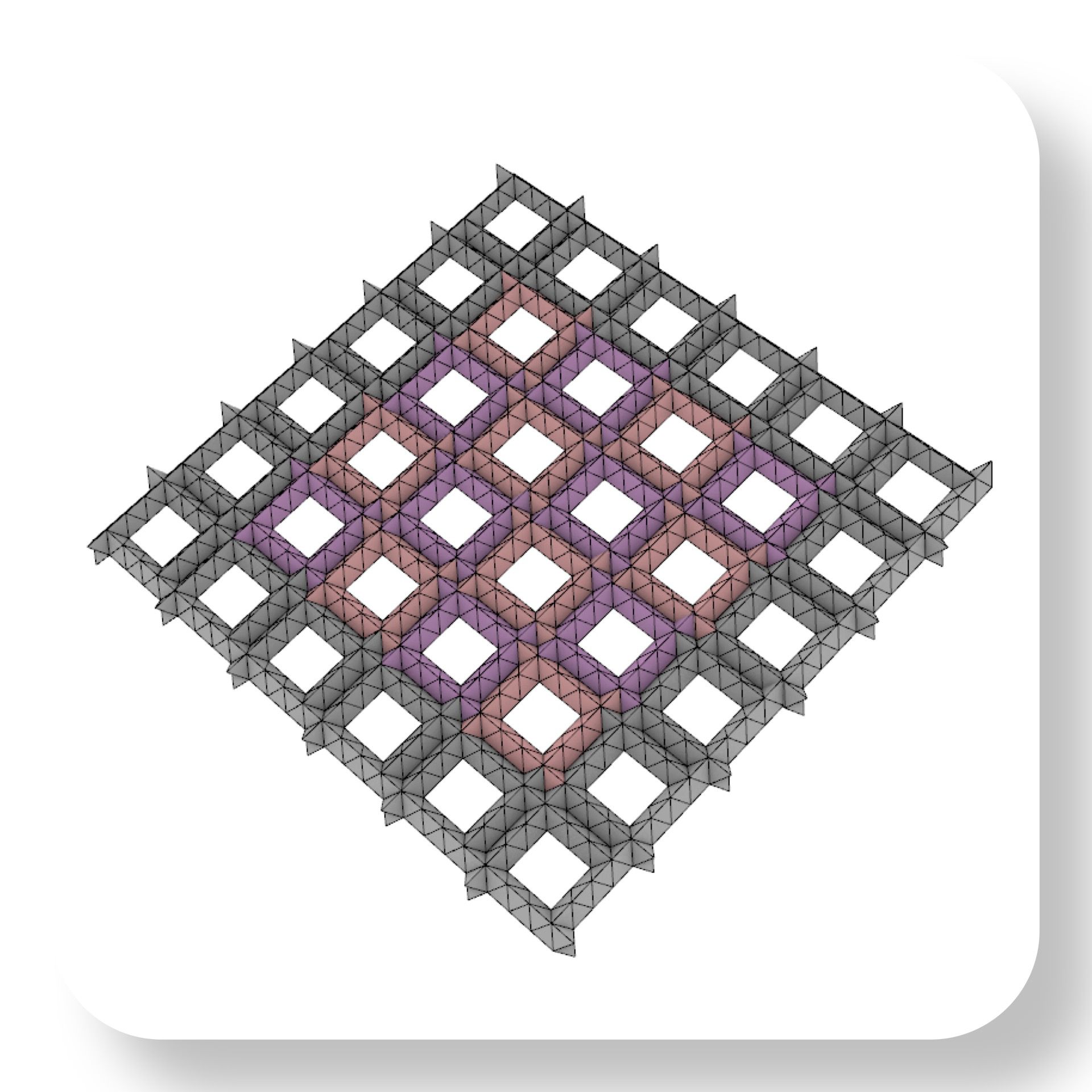}
\end{minipage}
\caption{Assembly of the $(3,3)$-shuriken}
\label{assembly33Shurkien}
\end{figure}

We provide topological interlocking assemblies of the $(m,n)$-shuriken by translating copies of the $(m,n)$-cushion by multiples of the vectors
\[
(n+1)v_1,(m+1)(v_2-v_3).
\]
Figure~\ref{assemblyShuriken} shows the described assembly of the shuriken and Figure~\ref{assembly33Shurkien} presents the assembly of $(3,3)$-shuriken. In both figures the frame is illustrated by the grey blocks.

\begin{figure}[H]
    \centering        \includegraphics[scale=1,viewport=10cm 21.5cm 0cm 26.5cm]{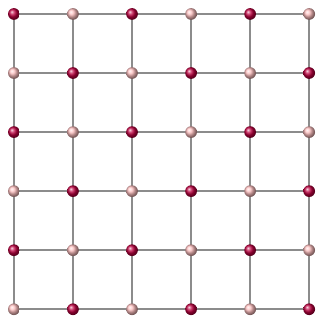}
\caption{Assembly graph of the given assembly of the shuriken}
\label{assemblygraphshuriken}
\end{figure}

These assemblies of the $(m,n)$-shuriken give rise to isomorphic assembly graphs. This graph is shown in Figure~\ref{assemblygraphshuriken}.

\section{Modifying Blocks}\label{section:Modyfying}
In this section, we present modifications of some blocks introduced in Section~\ref{section:NewBlocks}. More precisely, some of these blocks can be continuously deformed or truncated so that copies of the modified blocks can be arranged to form topological interlocking assemblies. In particular, we discuss the tetrahedron, the kitten and the $n$-cushion as examples of blocks that can be modified to obtain new blocks with interlocking properties.

\subsection{Truncation}
The \emph{truncation} of a polyhedron is the process of cutting off edges or vertices of the given polyhedron and therefore creating a new polyhedron with an increased number of faces. Here, we illustrate truncations of some of the presented blocks that yield modified blocks with interlocking properties. It is well known that truncating interlocking blocks, is a method to create new interlocking blocks. For instance, truncating a tetrahedron results in the Abeille block \citep{gallon_machines_1735,glickman_g-block_1984,dyskin_new_2001} and truncating a cube can lead to a compressed octahedron \citep{kanel-belov_interlocking_2010,spherical}.
\subsubsection{Tetrahedron}
\begin{figure}[H]
\begin{minipage}{.3\textwidth}
    \centering
    \includegraphics[height=4cm]{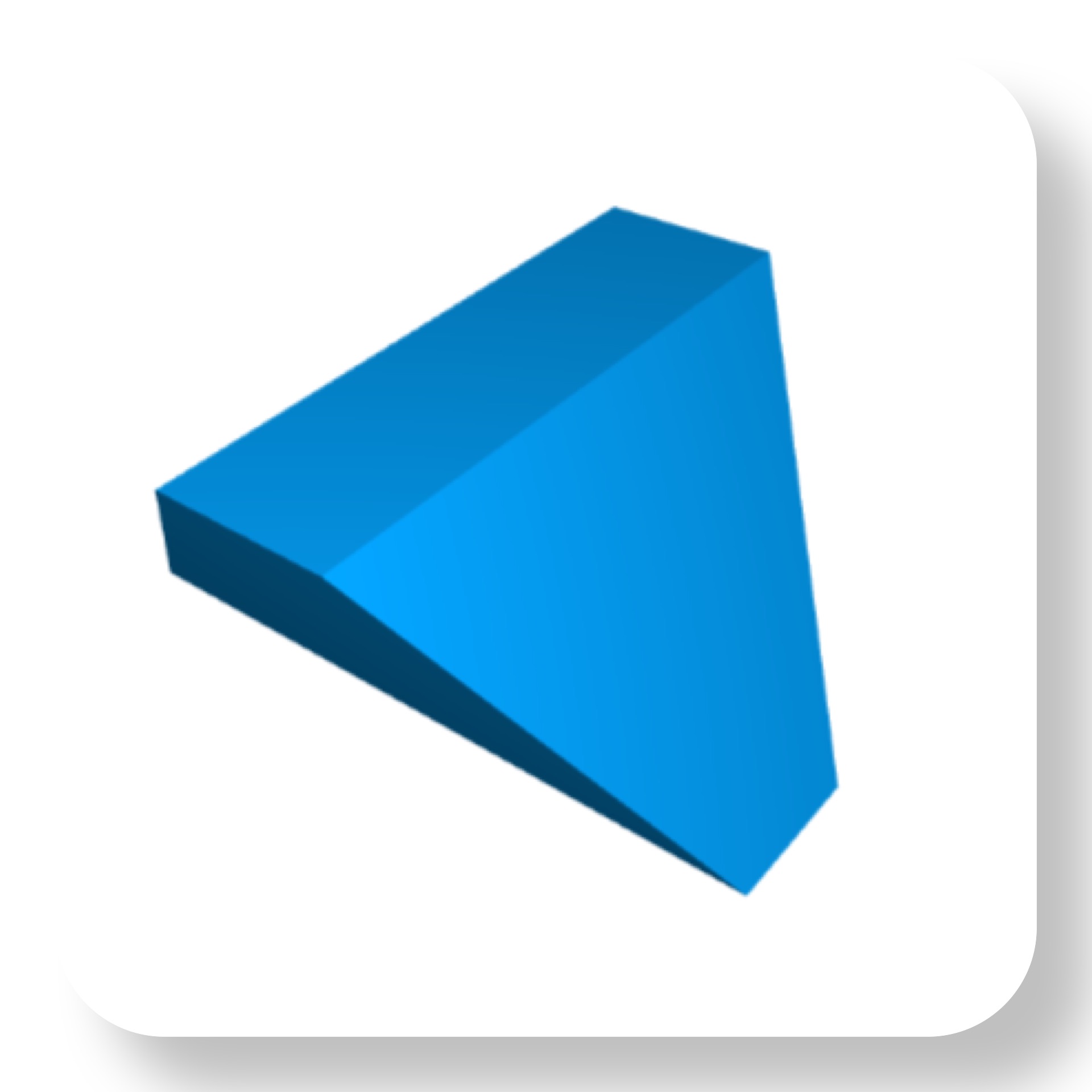}
\end{minipage}
\begin{minipage}{0.5cm}
    
\end{minipage}
\begin{minipage}{.3\textwidth}
    \centering
    \includegraphics[height=4cm]{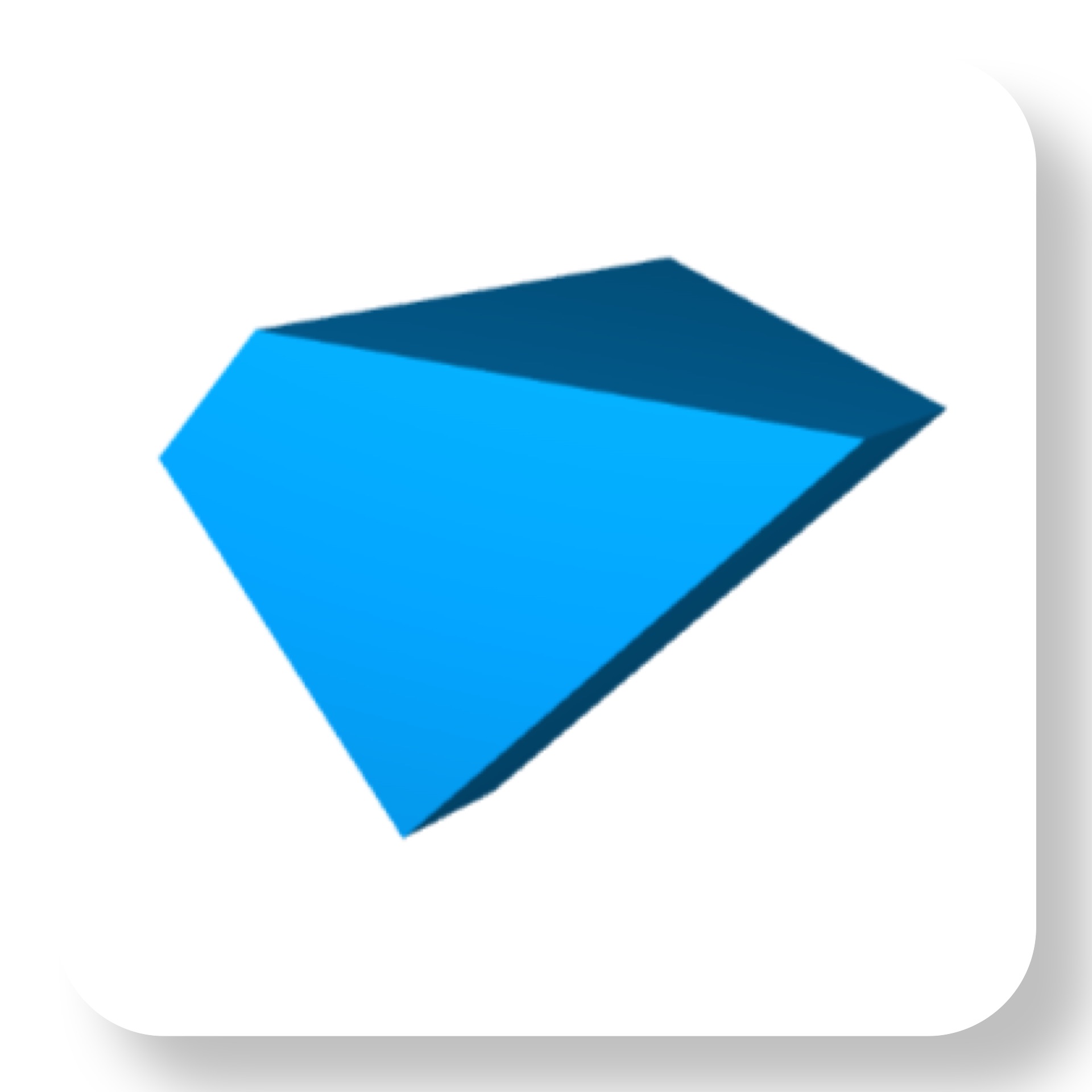}
\end{minipage}
\begin{minipage}{0.5cm}
    
\end{minipage}
\begin{minipage}{.3\textwidth}
    \centering
    \includegraphics[height=4cm]{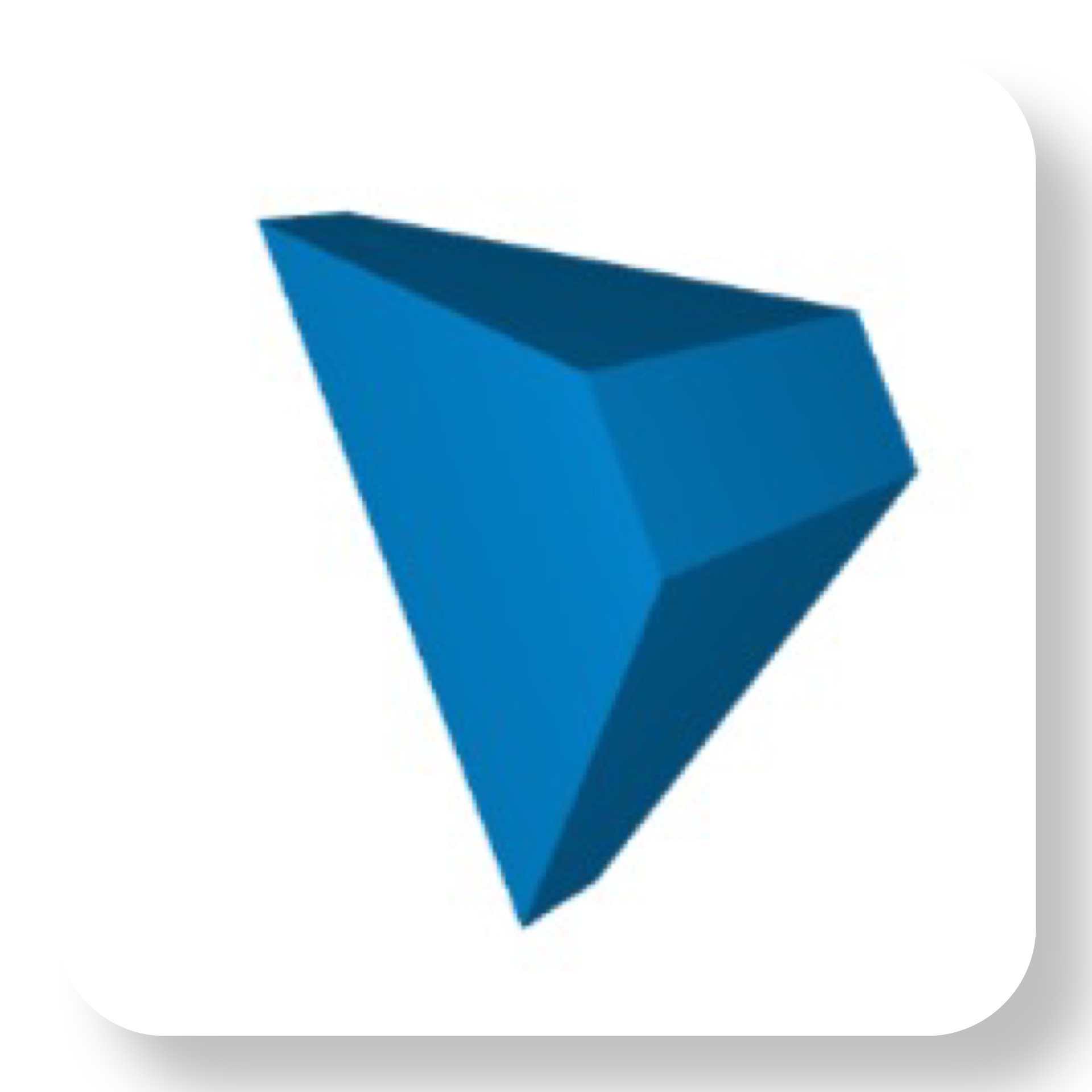}
\end{minipage}
\caption{ Truncating two opposite edges of a tetrahedron leads to the Abeille block, see \citep{gallon_machines_1735}}
\label{truncatedtetrahedron}
\end{figure}

Cutting off two opposite edges of a tetrahedron results in the block illustrated in Figure \ref{truncatedtetrahedron}.

\begin{figure}[H]
\begin{minipage}{.5\textwidth}
    \centering
    \includegraphics[height=5cm]{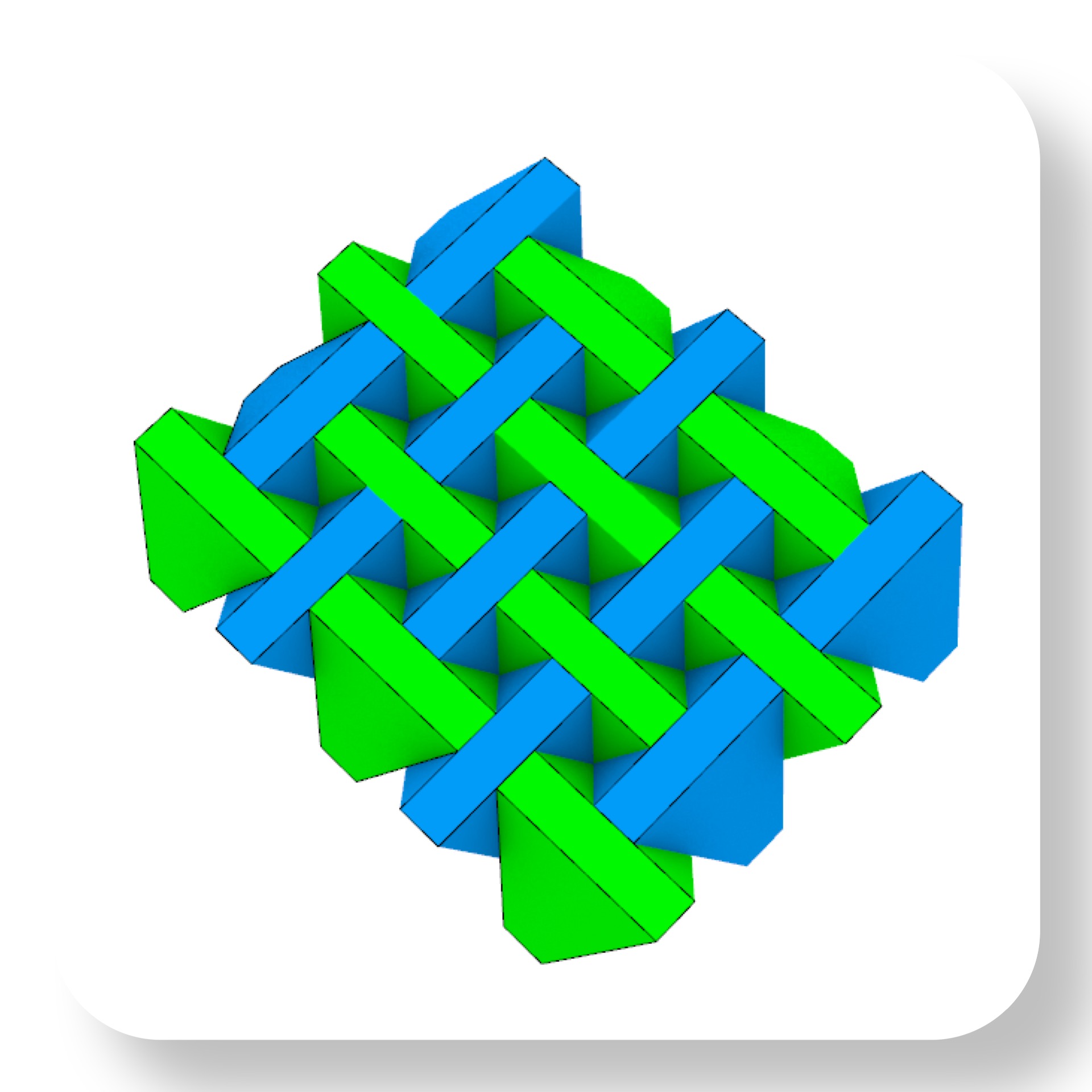}
\end{minipage}
\begin{minipage}{0.5cm}
    
\end{minipage}
\begin{minipage}{.5\textwidth}
    \centering
    \includegraphics[height=5cm]{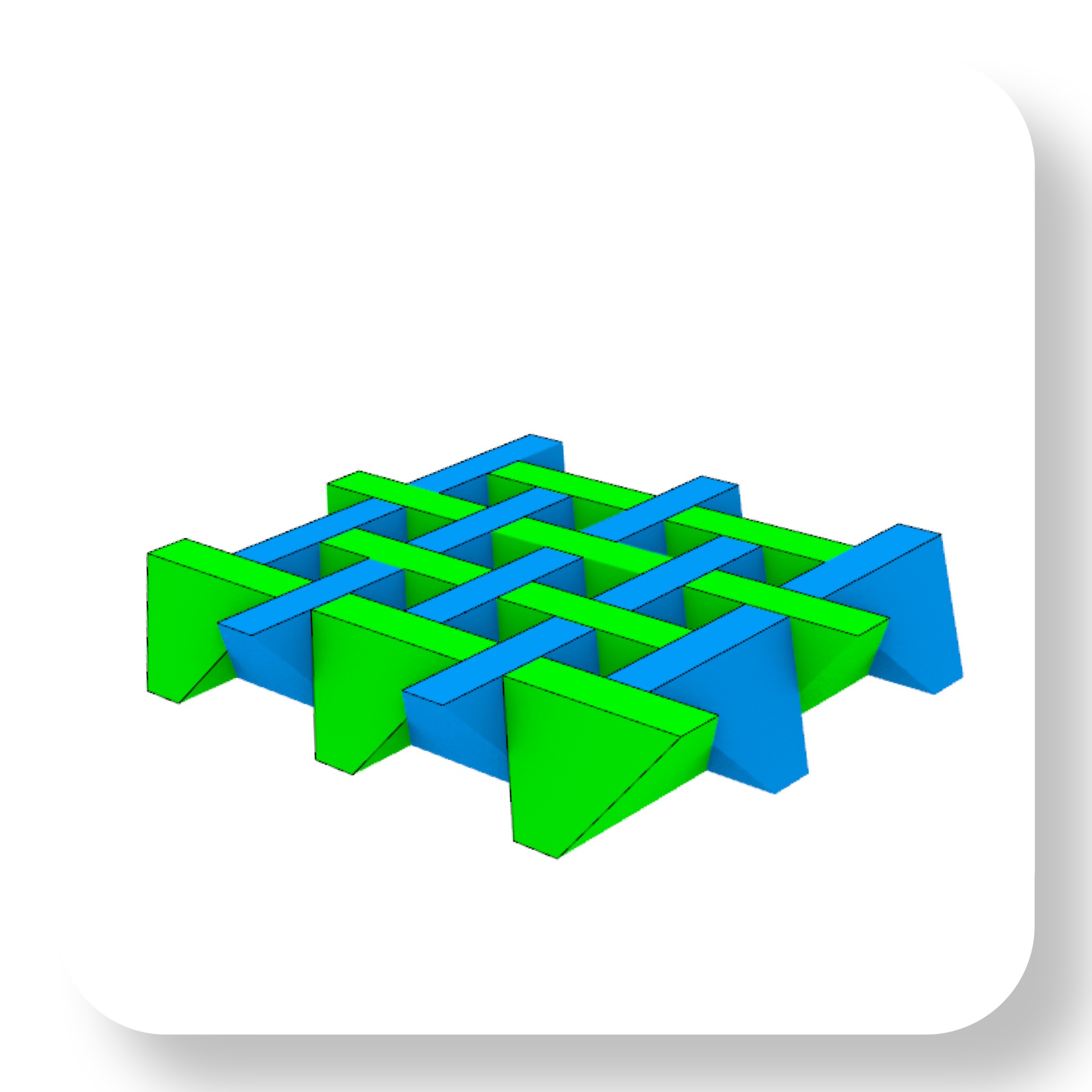}
\end{minipage}
\caption{ Various views of the assembly of the truncated tetrahedron also called flat Abeille vault in literature}
\label{assemblytruncatedtetrahedra}
\end{figure}

Introducing this modification to the tetrahedra in the assembly in Figure \ref{tetra_oct_interlocking} results in the interlocking shown in Figure \ref{assemblytruncatedtetrahedra}. 
Note, these modifications of the blocks do not decrease the number of contact faces of different blocks. More precisely, two tetrahedra in the given assembly share a contact face if and only if the truncated blocks share a contact face. Moreover, the modified blocks in the assembly can be seen as subsets of the tetrahedra that form the tetrahedra-interlocking. Thus, it can be concluded that the interlocking property also holds for the assembly of the modified blocks. 

\subsubsection{Kitten}
\begin{figure}[H]
\begin{minipage}{.3\textwidth}
    \centering
    \includegraphics[height=4cm]{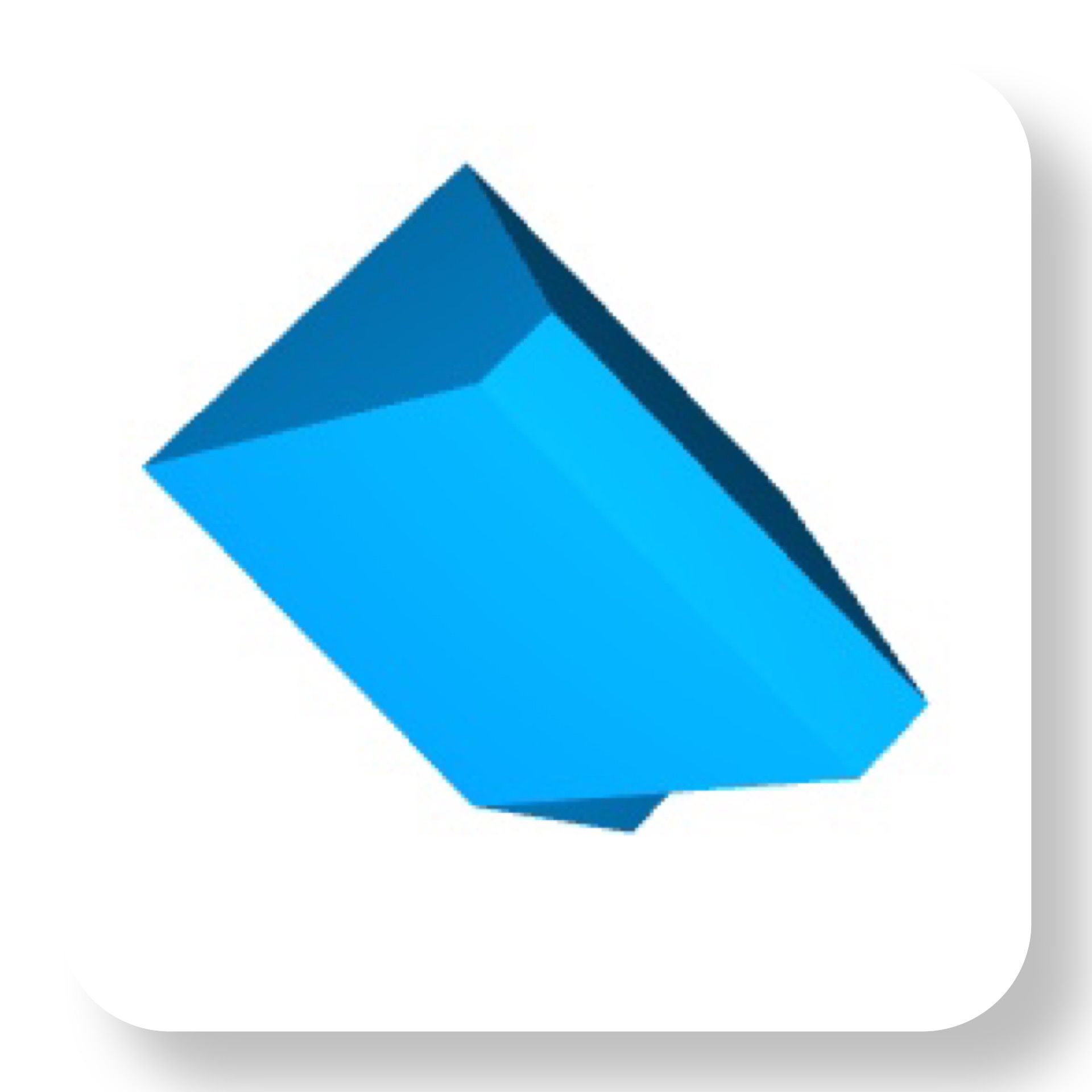}
\end{minipage}
\begin{minipage}{0.5cm}
    
\end{minipage}
\begin{minipage}{.3\textwidth}
    \centering
    \includegraphics[height=4cm]{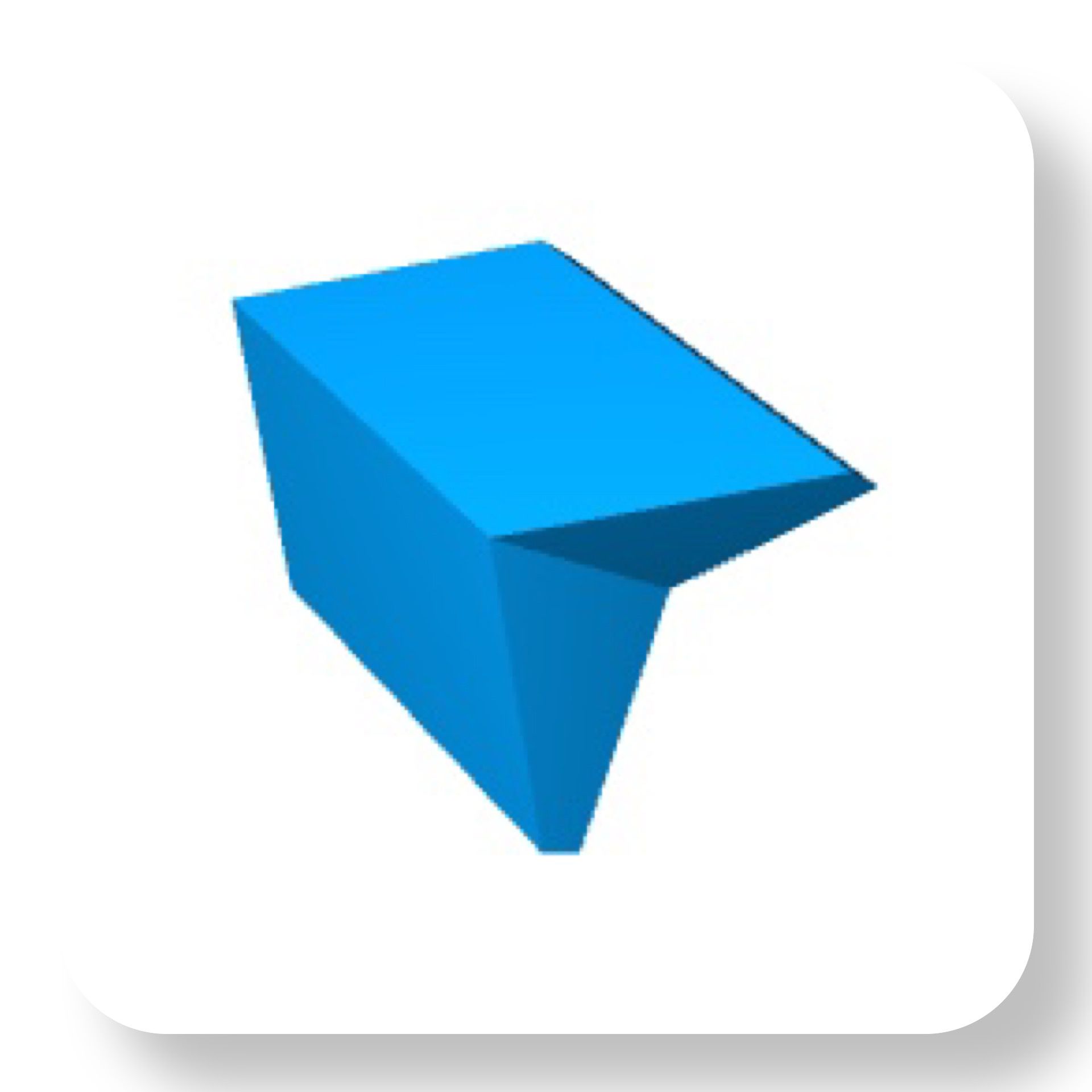}
\end{minipage}
\begin{minipage}{0.5cm}
    
\end{minipage}
\begin{minipage}{.3\textwidth}
    \centering
    \includegraphics[height=4cm]{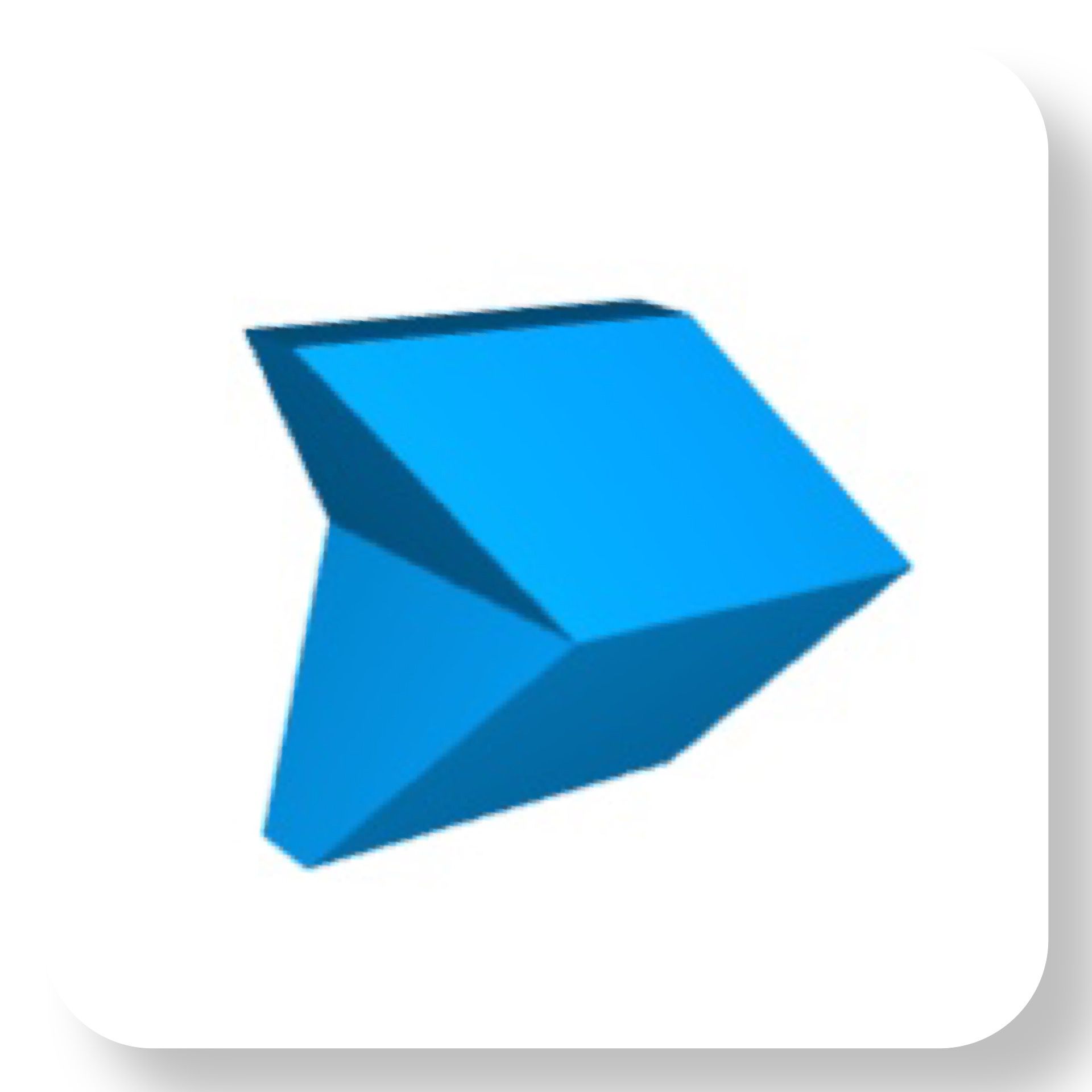}
\end{minipage}
\caption{ Truncating two opposite edges of the $n$-kitten }
\label{deformedkitten}
\end{figure}

Next we describe a truncation of the kitten and a corresponding interlocking assembly. There exists a pair of edges of the kitten that results in the block shown in Figure \ref{deformedkitten}.

\begin{figure}[H]
\begin{minipage}{.5\textwidth}
    \centering
    \includegraphics[height=5cm]{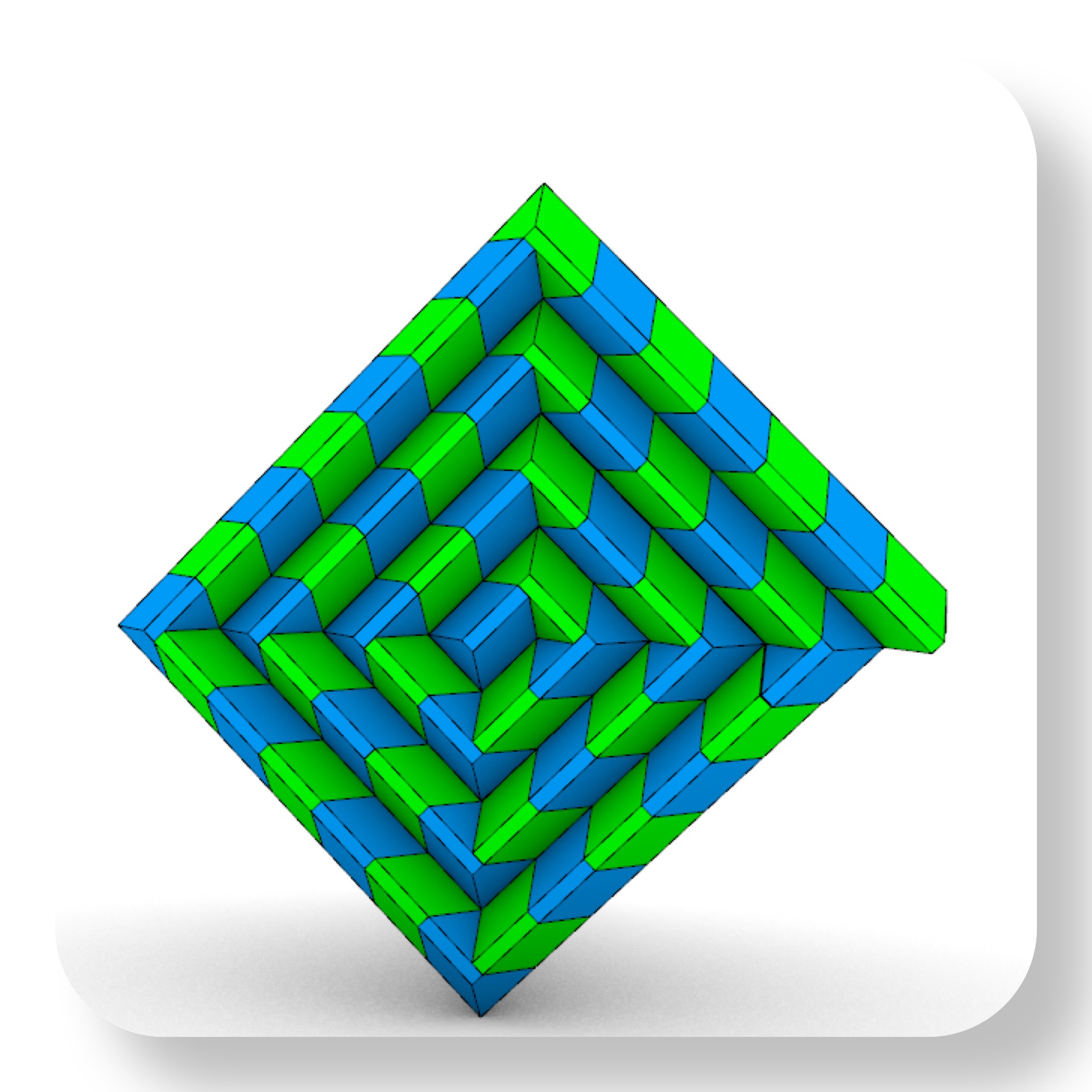}
\end{minipage}
\begin{minipage}{0.5cm}
    
\end{minipage}
\begin{minipage}{.5\textwidth}
    \centering
    \includegraphics[height=5cm]{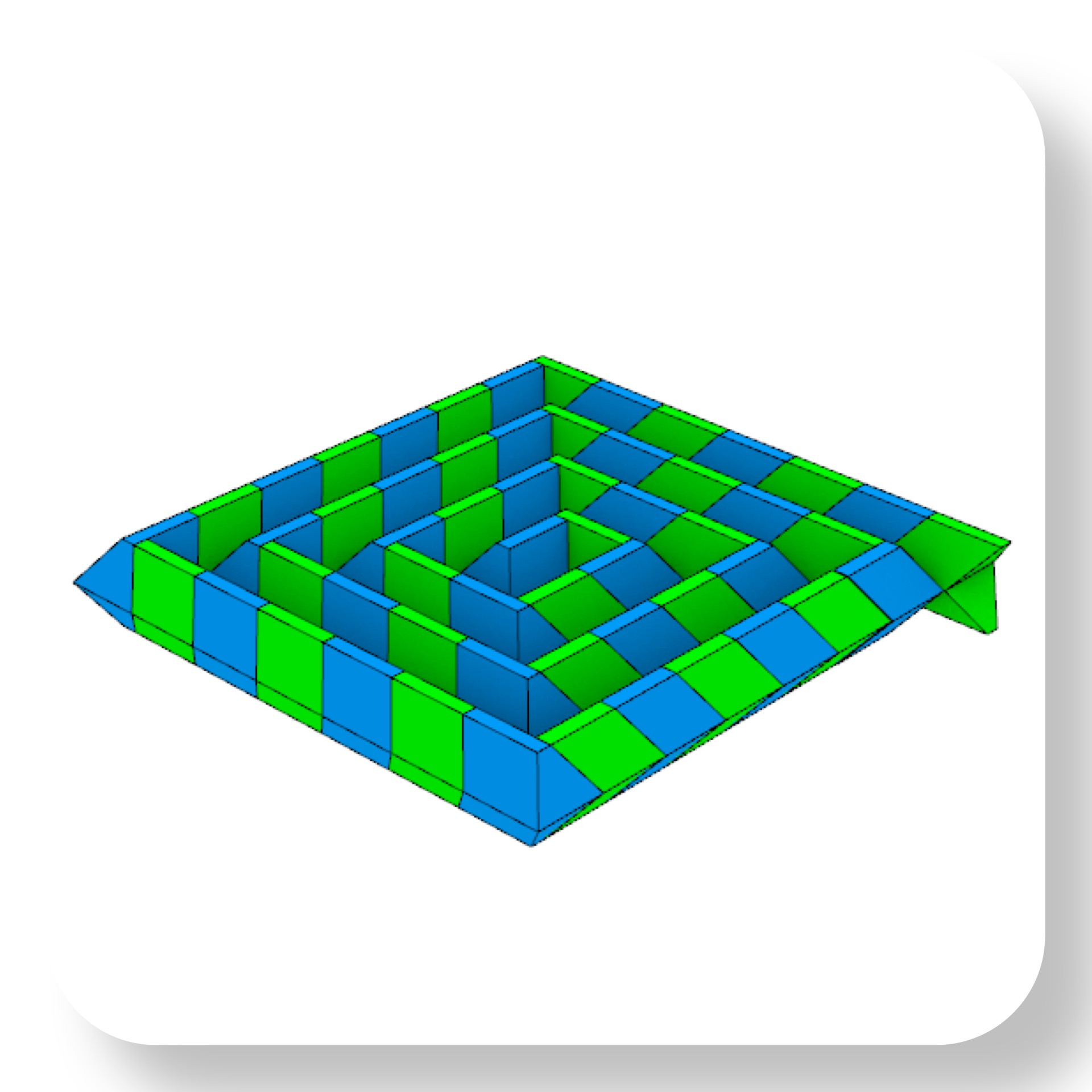}
\end{minipage}
\caption{ Various views of the assembly of the truncated kitten}
\label{assemblytruncatedkitten}
\end{figure}

In this case, applying this deformation to the copies of the kitten in the corresponding interlocking assembly does not increase the number of contact faces and due to the fact that the truncated blocks can be seen as a subset of the blocks in Figure \ref{assemblytruncatedkitten}, the interlocking property of the modified assembly can be established.

\subsubsection{Cushion}
\begin{figure}[H]
\begin{minipage}{.3\textwidth}
    \centering
    \includegraphics[height=4cm]{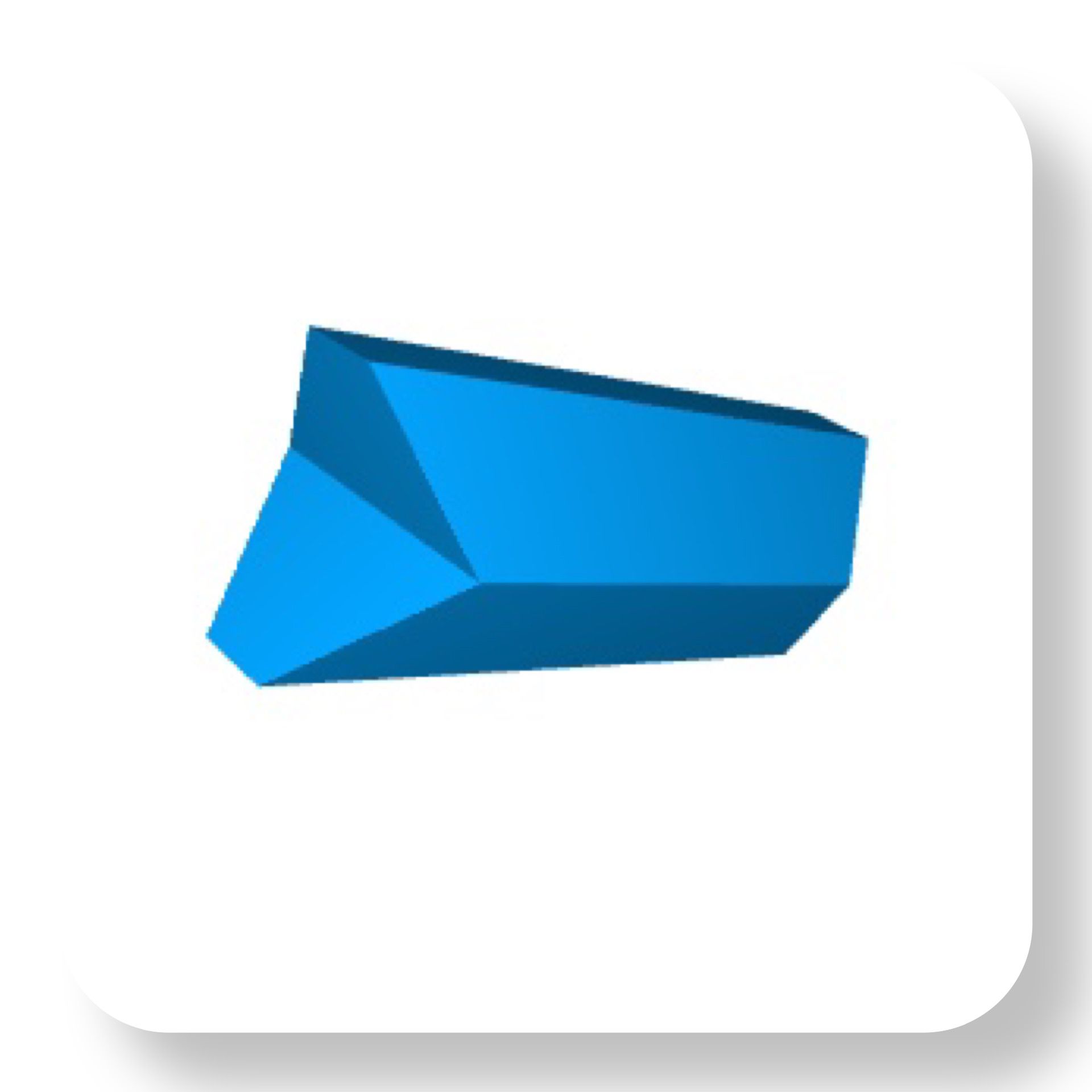}
\end{minipage}
\begin{minipage}{0.5cm}
    
\end{minipage}
\begin{minipage}{.3\textwidth}
    \centering
    \includegraphics[height=4cm]{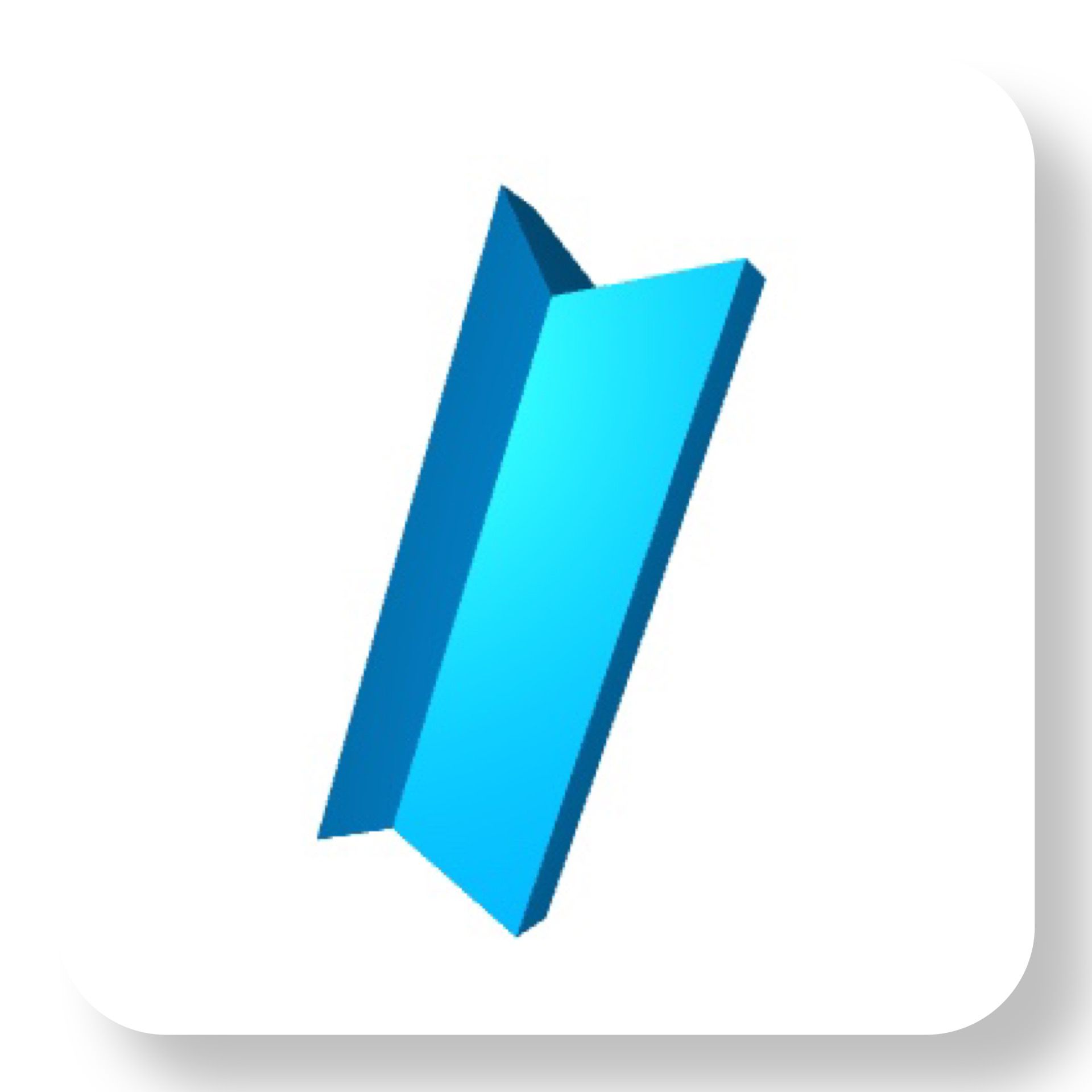}
\end{minipage}
\begin{minipage}{0.5cm}
    
\end{minipage}
\begin{minipage}{.3\textwidth}
    \centering
    \includegraphics[height=4cm]{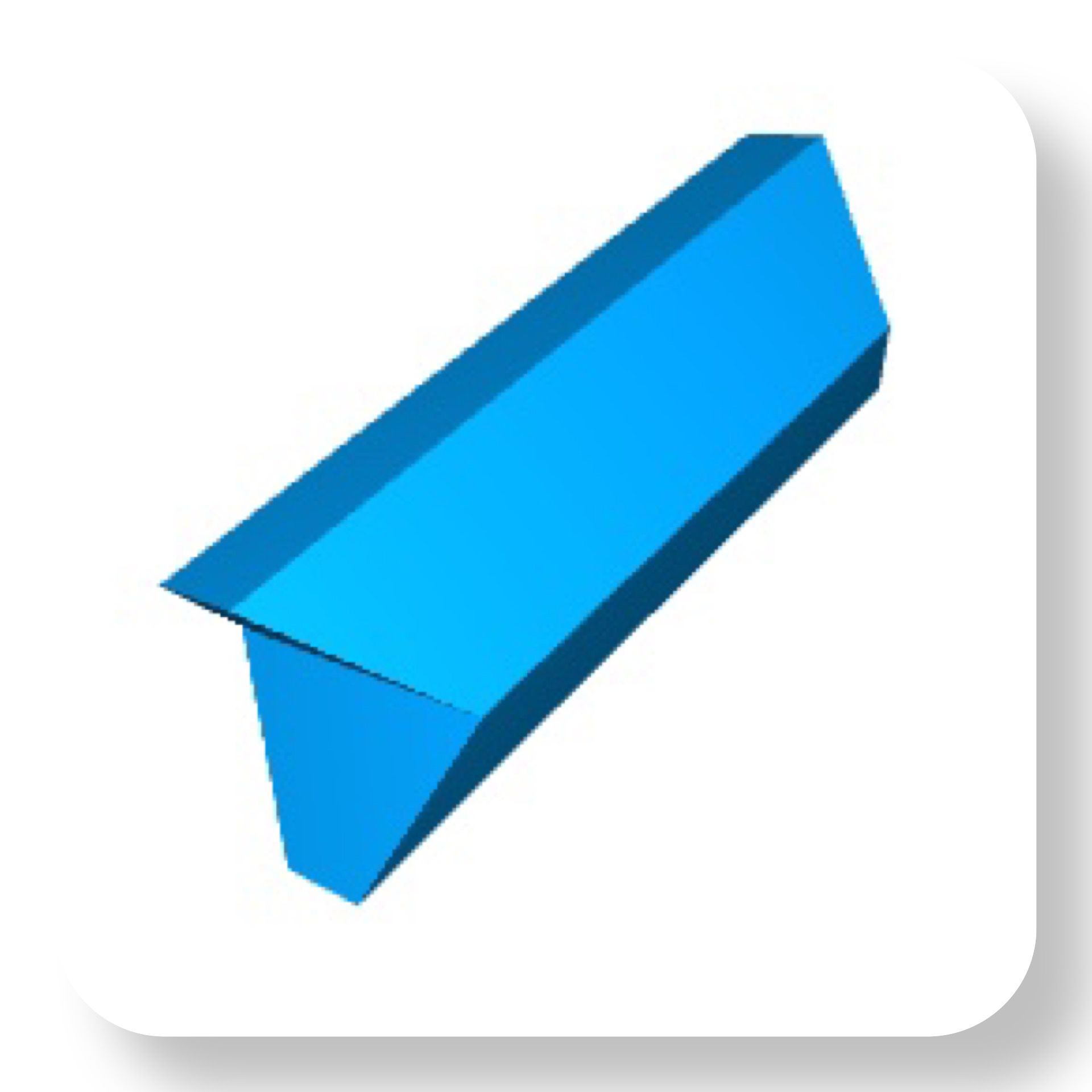}
\end{minipage}
\caption{Truncating two opposite edges of the $n$-cushion}
\label{truncatedcushion}
\end{figure}

\begin{figure}[H]
\begin{minipage}{.5\textwidth}
    \centering
    \includegraphics[height=5cm]{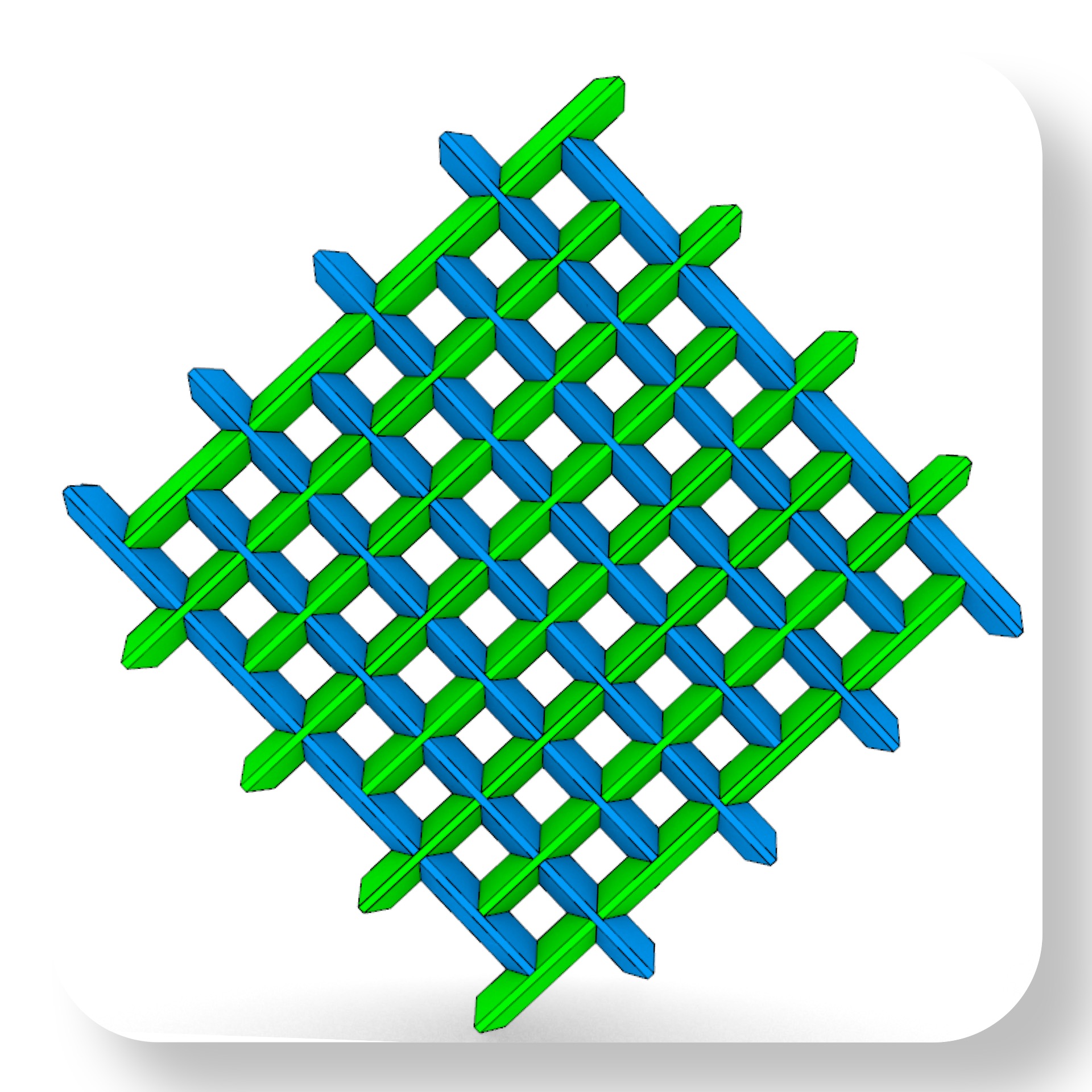}
\end{minipage}
\begin{minipage}{2cm}
    
\end{minipage}
\begin{minipage}{.5\textwidth}
    \centering
    \includegraphics[height=5cm]{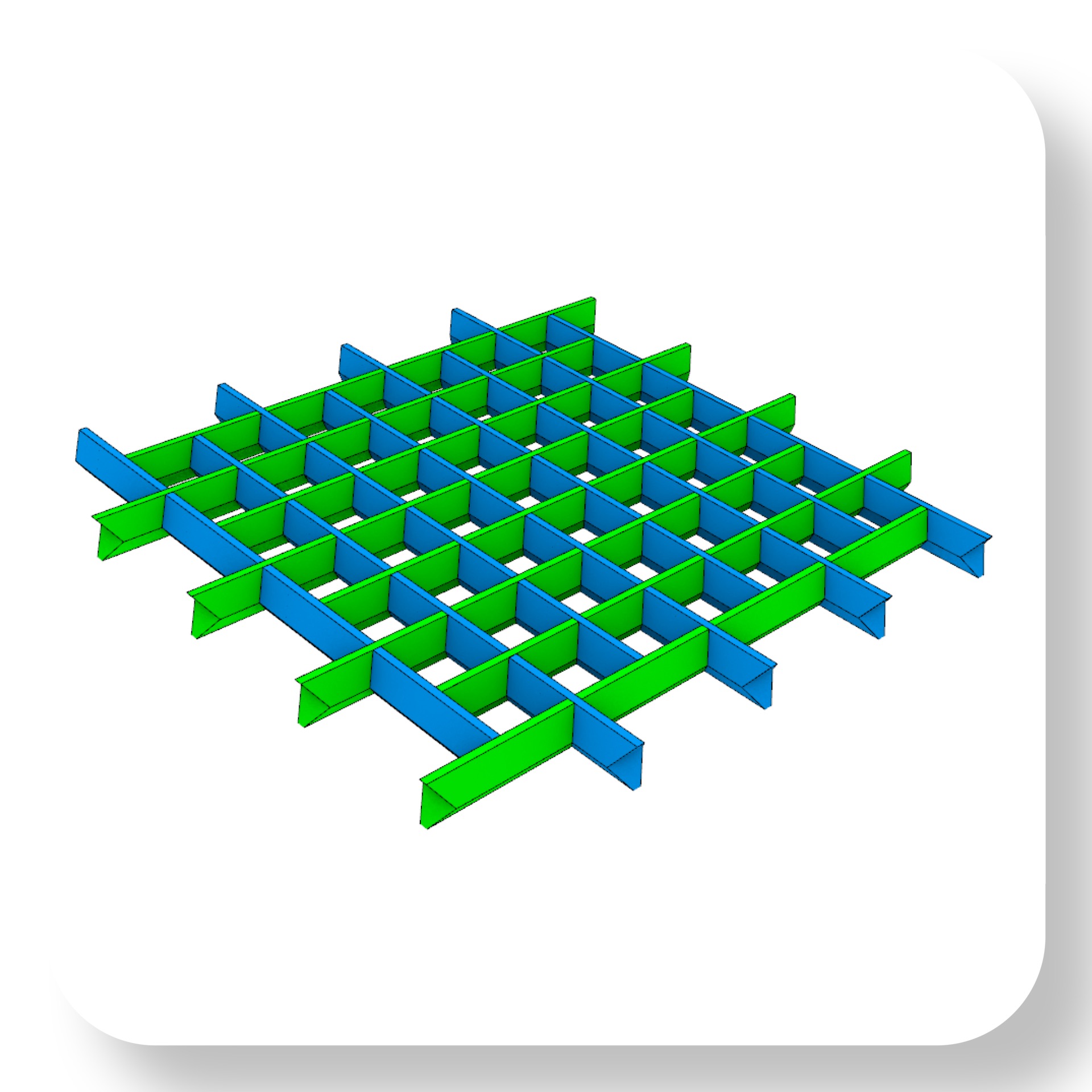}
\end{minipage}
\caption{ Various views of the assembly of the truncated kitten}
\label{assemblytruncatedcushion}
\end{figure}

Truncations can be applied to a pair of edges of the $n$-cushion to construct the block illustrated in Figure \ref{truncatedcushion}.
If copies of the modified block are assembled similarly as the assembly of the cushion in Figure \ref{assemblycushion}, the interlocking property can be proven by using similar arguments as in the case of the assembly of the modified kitten or the modified tetrahedron. This assembly is illustrated in Figure \ref{assemblytruncatedcushion}.

\subsection{Continuous Deformation}
\begin{figure}[H]
\begin{minipage}{.5\textwidth}
    \centering
    \includegraphics[height=5cm]{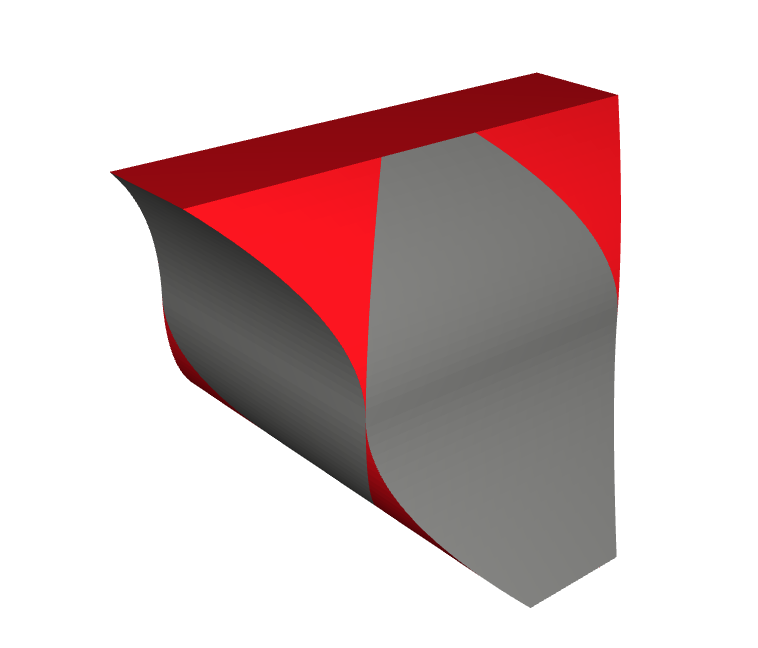}
\end{minipage}
\begin{minipage}{2cm}
    
\end{minipage}
\begin{minipage}{.5\textwidth}
    \centering
    \includegraphics[height=5cm]{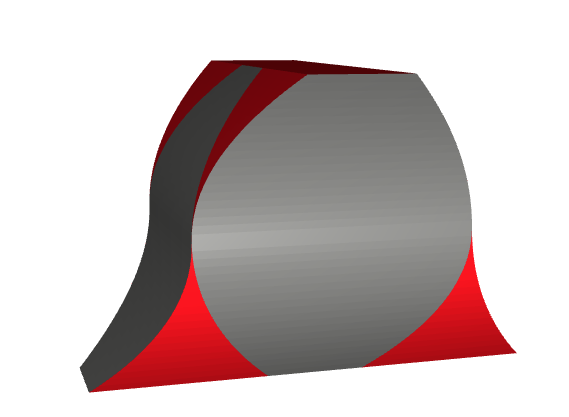}
\end{minipage}
\caption{ Various views of the deformed Abeille block}
\label{deformedTetrahedron}
\end{figure}

A \emph{continuous deformation}, also known as homeomorphism, is the gradual transition of an object or shape into another, while preserving specific properties or features consistently throughout the transformation. Here, a continuous deformation corresponds to a continuous map $\mathbb{R}^3\to \mathbb{R}^3$. For instance, we can continuously deform a unit sphere into a cube or a cup into a torus. Here, we illustrate interlocking assemblies that arise from deforming blocks of already existing topological interlocking assemblies. As an example, we consider the truncated tetrahedron, which can be continuously deformed into the block shown in Figure \ref{deformedTetrahedron}.
 
\begin{figure}[H]
\begin{minipage}{.5\textwidth}
    \centering
    \includegraphics[height=4.cm]{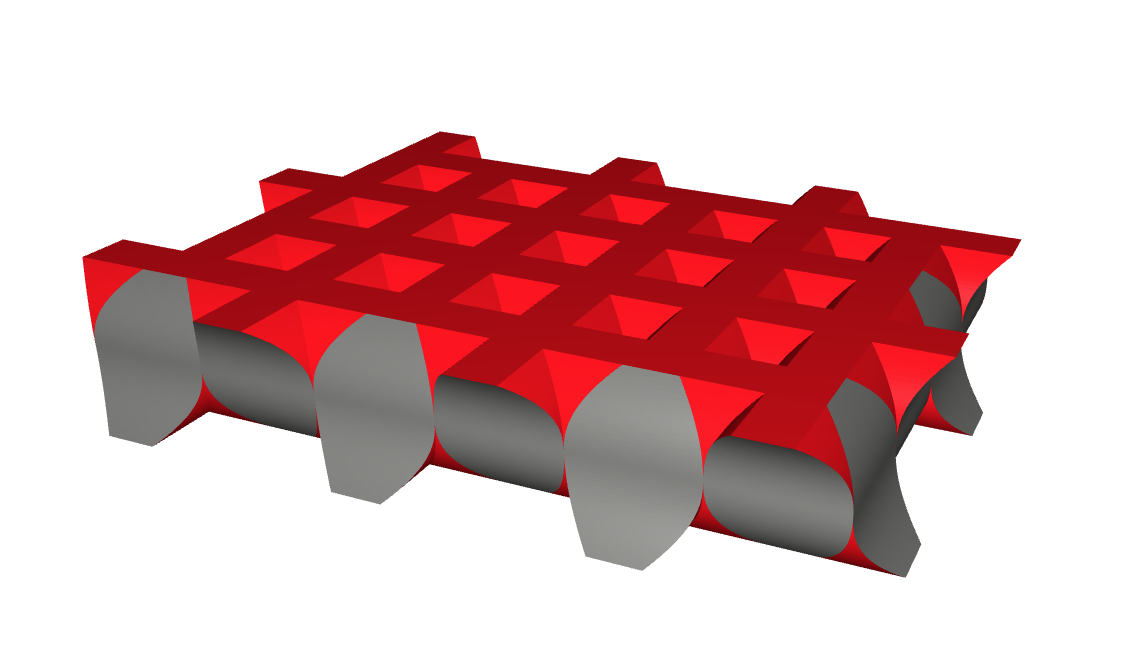}
\end{minipage}
\begin{minipage}{2cm}
    
\end{minipage}
\begin{minipage}{.5\textwidth}
    \centering
    \includegraphics[height=4.cm]{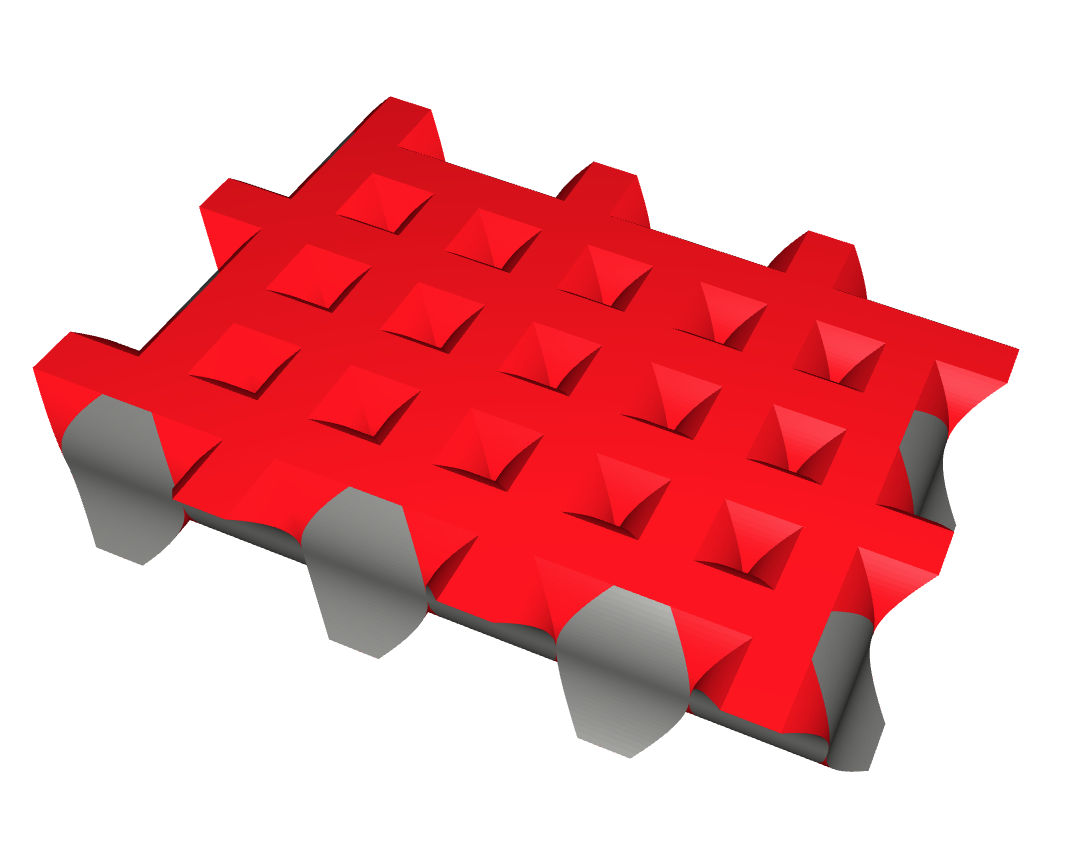}
\end{minipage}
\caption{ Various views of the assembly of the deformed block}
\label{assemblydeformedTetrahedron}
\end{figure}

Note, copies of this new block can be arranged similarly as copies of the tetrahedron in the tetrahedra-interlocking shown in Figure \ref{TetrahedraInterlocking}. This assembly of the modified blocks is shown in Figure \ref{assemblydeformedTetrahedron}. 
The relative contact area (coloured in grey in Figure \ref{assemblydeformedTetrahedron}) increases with the deformation. It can be shown that the interlocking property of the new assembly follows from the interlocking property of the truncated tetrahedron assembly.

\section{Conclusion}
This paper shows that the blocks constructed by combining tetrahedra and octahedra, so that the resulting body is contained in the tetrahedral-octahedral~honeycomb, allow many different planar and even three-dimensional assemblies. 

It might be interesting to investigate our proposed interlocking assemblies with the infinitesimal interlocking criteria given by \cite{wang_design_2019}.  
It is well-known that cubes can be assembled to obtain a TIA.
Since the vertices of a cube can be chosen, so that they are contained in the primitive cubic lattice, i.e.\ the lattice spanned by the unit vectors in $\mathbb{R}^3$, examining other three-dimensional lattices might give rise to interesting examples of interlocking blocks.

Moreover, deformations of the proposed blocks can be analysed to investigate whether it is possible to realise structures such as vaults and sphere as a topological interlocking assembly created out of the deformed blocks.

\section*{Declarations}

\subsection*{Funding}
The authors gratefully acknowledge the funding by the Deutsche Forschungsgemeinschaft (DFG, German Research Foundation) in the framework of the Collaborative Research Centre CRC/TRR 280 “Design Strategies for Material-Minimized Carbon Reinforced Concrete Structures – Principles of a New Approach to Construction” (project ID 417002380). 

Tom Goertzen was partially supported by the FY2022 JSPS Postdoctoral Fellowship for Research in
Japan (Short-term), ID PE22747.

\subsection*{Conflicts of Interest}
The authors declare that they have no conflicts of interest.

\subsection*{Availability of Data and Material}
The data and material used in this research are available upon request.

\subsection*{Code Availability}
The custom code developed for this research can be obtained by contacting the corresponding author.

\bibliography{Interlocking,other_citations}

\end{document}